\documentclass[reqno]{amsart}

\newtheorem{Theorem}{Theorem}[section]
\newtheorem{Corollary}[Theorem]{Corollary}
\newtheorem{Lemma}[Theorem]{Lemma}
\newtheorem{Proposition}[Theorem]{Proposition}
\theoremstyle{definition}
\newtheorem{Definition}[Theorem]{Definition}
\theoremstyle{remark}
\newtheorem{Remark}[Theorem]{Remark}
\numberwithin{equation}{section}

\newcommand{\diag}{\operatorname{diag}}
\newcommand{\Image}{\operatorname{Im}}
\newcommand{\hess}{\operatorname{Hess}}

\newcommand{\delt}{\mbox{\boldmath$\delta$}}
\newcommand{\grad}{\operatorname{grad}}
\newcommand{\ind}{\operatorname{Ind}}
\newcommand{\vol}{\operatorname{Vol}}
\newcommand{\dist}{\operatorname{dist}}
\newcommand{\divv}{\overset v {\operatorname{div}}\,}
\newcommand{\divh}{\overset h {\operatorname{div}}\,}
\newcommand{\divm}{\overset {m} {\operatorname{div}}\,}
\newcommand{\de}[2]{\frac{\partial #1}{\partial #2}}

\renewcommand{\d}{d^s}

\begin{document}
\title{The boundary rigidity problem in the presence
of~a~magnetic field}
\author[N.S. Dairbekov]{Nurlan S. Dairbekov}
\address{Laboratory of Mathematics,
Kazakh British Technical University,
Tole bi 59, 050000 Almaty, Kazakhstan }
\email{Nurlan.Dairbekov@gmail.com}
\thanks{First, third, and forth authors
partly supported by CRDF Grant KAM1-2851-AL-07}

\author[G.P. Paternain]{Gabriel P. Paternain}
 \address{Department of Pure Mathematics and Mathematical Statistics,
University of Cambridge,
Cambridge CB3 0WB, England}
 \email {g.p.paternain@dpmms.cam.ac.uk}

\author[P. Stefanov]{Plamen Stefanov}
\address{ Department of Mathematics,
Purdue University,
West Lafayette, IN 47907, USA}
\email {stefanov@math.purdue.edu}
\thanks{Third author partly supported by
NSF Grant DMS-0400869}

\author[G. Uhlmann]{Gunther Uhlmann}
\address{Department of Mathematics,
University of Washington,
Seattle, WA 98195-4350, USA}
\email{gunther@math.washington.edu}
\thanks{Fourth author partly supported by
NSF and a Walker Family Endowed Professorship}

\begin{abstract}
For a compact Riemannian manifold with boundary, endowed with a magnetic potential $\alpha$, we consider the problem of restoring the  metric $g$ and  the magnetic
potential $\alpha$ from the values of the Ma\~n\'e action potential
between boundary points and the associated linearized problem. We study simple magnetic systems. In this case, knowledge of the Ma\~n\'e action potential is equivalent to knowledge of the scattering relation on the boundary which maps a starting point and a direction of a magnetic geodesic into its end point and direction. This problem can only be solved up to an isometry and a gauge transformation of $\alpha$.

For the linearized problem, we show injectivity, up to the natural obstruction,  under explicit bounds on the curvature and on $\alpha$. We also show injectivity and stability for $g$ and $\alpha$ in a generic class $\mathcal{G}$    including real analytic ones.

For the nonlinear problem, we show rigidity for real analytic simple ($g, \alpha)$,
rigidity for metrics in a given conformal class, and locally, near any $(g,\alpha)\in \mathcal{G}$.
We also show that simple magnetic systems on two-dimensional manifolds are always rigid.
\end{abstract}

\maketitle
\tableofcontents

\section{Introduction}
\subsection{Statement of the problem}
Let $M$ be a compact manifold with boundary,
endowed with a Riemannian metric~$g$, and let
$\pi:TM\to M$ denote the canonical projection,
$\pi:(x,\xi)\mapsto x$,
$x\in M$, $\xi\in T_xM$.

Denote by $\omega_0$ the canonical symplectic form on $TM$, which
is obtained by pulling back the canonical symplectic form of $T^*M$
via the Riemannian metric.
Let $H:TM\to\mathbb R$ be defined by
$$
H(v)=\frac 12|v|^2_g,\quad v\in TM.
$$
The Hamiltonian flow of $H$ w.r.t.\ $\omega_0$
gives rise to the geodesic flow of $M$.
Let $\Omega$ be a closed $2$-form on $M$ and consider
the new symplectic form $\omega$ defined as
$$
\omega=\omega_0+\pi^*\Omega.
$$
The Hamiltonian flow of $H$ w.r.t.\ $\omega$ gives rise to
the {\em magnetic} (or {\em twisted geodesic}) {\em flow}
$\psi^t:TM\to TM$.
This flow models the motion of
a unit charge of unit mass in a magnetic field
whose Lorentz force $Y:TM\to TM$
is the bundle map uniquely determined by
\begin{equation}\label{lorentz}
\Omega_x(\xi,\eta)=\langle Y_x(\xi),\eta\rangle
\end{equation}
for all $x\in M$ and $\xi,\eta\in T_xM$.

Magnetic flows were first considered by V. I. Arnold in \cite{Ar}
and by D. V. Anosov and Y. G. Sinai in \cite{AS}.
As shown in \cite{ArG, N1, N2, NS, Ko, PP}, they are closely related to other problems of classical mechanics, mathematical physics,
symplectic geometry, and dynamical systems.

It is not hard to show that the generator $\mathbf G_\mu$
of the magnetic flow is
$$
\mathbf G_\mu(x,\xi)=\mathbf G(x,\xi)+Y^j_i(x)\xi^i\de{}{\xi^j},
$$
where $\mathbf G$ is the generator of the geodesic flow,
and that every trajectory of the magnetic flow
is a curve of the form $t\mapsto (\gamma(t),\dot\gamma(t))$,
where $\gamma$ is a curve on $M$ which satisfies the equation
\begin{equation}\label{geodesic}
\nabla_{\dot\gamma}\dot\gamma=Y(\dot\gamma),
\end{equation}
which is nothing but Newton's law of motion.
Such a curve $\gamma$ is called a {\em magnetic geodesic}.
Note that time is not reversible
on the magnetic geodesics, unless $\Omega=0$.
If $\Omega=0$ we recover the ordinary geodesic flow
and ordinary geodesics.

When $\Omega$ is exact, i.e. $\Omega=d\alpha$ for some magnetic potential $\alpha$,
the magnetic flow also arises as the Hamiltonian flow of
$H(x,p) = \frac12 (p+\alpha)^2_g$ with respect to the standard symplectic form
of $T^*M$. The function $H$ is the symbol of the semiclassical magnetic Schr\"odinger operator.

Since the magnetic flow preserves the level sets of
the Hamiltonian function $H$,
every magnetic geodesic has constant speed.

Unlike the geodesic flow, where the flow is the same (up to time scale)
on any energy levels, the magnetic flow depends essentially on the choice
of the energy level. We fix the energy level $H^{-1}(1/2)$,
thus considering the magnetic flow on the unit sphere
bundle $SM$ of $M$, in consequence we consider only the unit speed
magnetic geodesics. Note that fixing the energy level to be
$SM$ is no restriction at all,
since one can obtain the behavior in any energy level
by considering the flow on $SM$
upon changing $\Omega$ by $\lambda\,\Omega$, where $\lambda\in\mathbb R$.

We define the action $\mathbb A(x,y)$ between boundary points as
a minimizer of the appropriate action functional, see (\ref{axy})
and Appendix~A. In the case $\Omega=0$,
$\mathbb A(x,y)$ coincides with the boundary distance function $d_g(x,y)$.
In this case, we cannot recover $g$ from $d_g$ up to isometry,
unless some additional assumptions are imposed on $g$, see, e.g.,   \cite{Cr1}. One such assumption is the simplicity of the metric, see, e.g., \cite{M, Sh1, SU, SU1}. We consider below the analog of simplicity for magnetic systems.

Let $\Lambda$ stand for the second fundamental form of $\partial M$
and $\nu(x)$ for the inward unit normal to $\partial M$ at $x$.
We say that $\partial M$ is {\em strictly magnetic convex}
if
\begin{equation}\label{strict-conv}
\Lambda(x,\xi)>\langle Y_x(\xi),\nu(x)\rangle
\end{equation}
for all $(x,\xi)\in S(\partial M)$
(see Appendix A for an explanation).
Note that if we replace $\xi$ by $-\xi$, we can put an absolute
value in the right-hand side of \eqref{strict-conv}.
In particular, magnetic convexity is stronger
than Riemannian one.

For $x\in M$, we define the {\em magnetic
exponential map} at $x$ to be the partial map
$\exp^\mu_x:T_xM\to M$ given by
$$
\exp^\mu_x(t\xi)=\pi\circ\psi^t(\xi),\quad t\ge 0,\ \xi\in S_xM.
$$
It is not hard to show that, for every $x\in M$,
$\exp^\mu_x$  is a $C^1$-smooth partial map on $T_xM$
which is $C^\infty$-smooth on $T_xM\setminus\{0\}$
(see Appendix A).

\begin{Definition}\label{simple}
We say that $M$ is {\em simple}
(w.r.t.\ $(g,\Omega)$) if $\partial M$ is strictly magnetic convex
and the magnetic exponential map
$\exp^\mu_x:(\exp^\mu_x)^{-1}(M)\to M$
is a diffeomorphism for every $x\in M$
(cf. the definition of a simple Riemannian manifold \cite{PU}).
\end{Definition}

In this case,
$M$ is diffeomorphic to the unit ball of $\mathbb R^n$
(so we can assume that $M$ is this ball); therefore,
$\Omega$ is exact, and we let $\alpha$ be
a {\em magnetic potential},
i.e., $\alpha$ is a $1$-form on $M$ such that
\begin{equation}\label{d-alpha}
d\alpha=\Omega.
\end{equation}

Henceforth we call $(g,\alpha)$
a~{\em simple magnetic system} on $M$. We will also say
that $(M,g,\alpha)$ is a~simple magnetic system.

Given $x,y\in M$, let
$$
\mathcal C(x,y)=\{\gamma:[0,T]\to M: T>0,\, \gamma(0)=x,\gamma(T)=y, \gamma\
\text{is absolutely continuous}\}.
$$
The {\em time-free action} of a curve $\gamma\in\mathcal C(x,y)$
w.r.t.\ $(g,\alpha)$ is defined as
\begin{equation}\label{action}
\mathbb A(\gamma)=\mathbb A_{g,\alpha}(\gamma)
=\frac12\int_0^T |\dot\gamma(t)|_g^2\,dt +\frac12 T-\int_\gamma\alpha.
\end{equation}

For a simple magnetic system,
unit speed magnetic geodesics
minimize the time-free action (Lemma \ref{minimize} in Appendix A).
More precisely,
\begin{equation}\label{axy}
\mathbb A(x,y):=\inf_{\gamma\in \mathcal C(x,y)}\mathbb A(\gamma)
=\mathbb A(\gamma_{x,y})
=T_{x,y}-\int_{\gamma_{x,y}}\alpha,
\end{equation}
where $\gamma_{x,y}:[0,T_{x,y}]\to M$ is
the unique unit speed magnetic geodesic from $x$ to $y$.

The function $\mathbb A(x,y)$ is referred to as
{Ma\~n\'e's action potential}
(of energy $1/2$), and we call
the restriction $\mathbb A|_{\partial M\times\partial M}$
the {\em boundary action function}.

Now, we state the {\em boundary rigidity problem in the presence
of a magnetic field} as follows: To which extent is a magnetic system
$(g,\alpha)$ determined by the boundary action function?

To be more precise, we say that
two simple magnetic systems $(g,\alpha)$ and $(g',\alpha')$ are
{\em gauge equivalent} if there are a diffeomorphism
$f:M\to M$, which is the identity on the boundary,
and a function $\varphi:M\to\mathbb R$, vanishing on the boundary,
such that $g'=f^*g$ and $\alpha'=f^*\alpha+d\varphi$.
Observe that gauge equivalent magnetic systems have
the same boundary action function.

Now, we rephrase the above problem as follows:
Given a simple magnetic system, is any
other simple magnetic system on the same manifold
gauge equivalent to the former as soon as it has
the same boundary action function?
If so, we call the given magnetic system {\em magnetic boundary rigid}.

Surely, this problem can be considered under various natural
restrictions. For example, we can consider it for a fixed
Riemannian metric and try to restore a magnetic potential, or, vice versa, fix
a magnetic potential and try to restore a metric, or consider the problem
for metrics in a fixed conformal class, etc.
In particular, for the zero magnetic potentials we recover
the ordinary boundary rigidity problem for Riemannian metrics
(see recent surveys of the latter in \cite{Cr2, SU2}).

For simple magnetic systems, knowledge of the action $\mathbb{A}(x,y)$ between boundary points is equivalent to knowing the scattering relation, see section~\ref{stattering}. For non-simple systems, the problem of recovering $(g,\Omega)$ from the scattering relation is the natural problem to study. The scattering relation appears naturally in the study of the scattering operator in $\mathbb R^n$ with $g$ Euclidean outside a ball, and $\alpha$ compactly supported. Namely, for non-trapping metrics, the scattering operator associated to the semi-classical magnetic Schr\"odinger operator is a Fourier integral operator with canonical relation that determines the scattering relation on a large sphere, \cite{A1, A2, G}. (It should be noted that for magnetic Schr\"odinger operators the resolvent is also a Fourier integral operator.)

\subsection{Description of the results} In section~\ref{sec_2}, we show that for simple magnetic systems, the action determines the jets of $g$ and $\alpha$ in boundary normal coordinates. We define the scattering relation and show that for simple magnetic systems, it determines $\mathbb{A}(x,y)$ on the boundary, and vice-versa. We also show that one can recover the volume from $\mathbb{A}(x,y)$.

In section~\ref{sec_3},  we  study the linearized problem. This reduces to the magnetic ray transform $I$. We show that  potential pairs (see Definition~\ref{p-def}) belong to the  kernel of $I$. We say that $I$ is s-injective, if the kernel of $I$ consists only of potential pairs.

In section~\ref{analysis-N},  we  show that the normal operator $N=I^*I$ is a pseudo-differential operator in the interior of $M$, elliptic on solenoidal pairs that are an orthogonal complement of the potential pairs. We construct a parametrix of $N$; near the boundary, additional arguments are needed. This parametrix recovers the solenoidal projection $\mathbf{f}^s$ given $N\mathbf{f}$, up to a smoothing term. We show that s-injectivity implies a stability estimate, uniform near any $(g,\alpha)$,  in appropriate spaces, see Theorem~\ref{thm_estimate}. We consider in this section and in Appendix~B real analytic $(g,\alpha)$ and show that then $I$ is s-injective. A crucial element of the proof is that $N$ is an analytic pseudo-differential operator in the interior of $M$. This is delicate since the magnetic  exponential map is only $C^1$ smooth when $\Omega\not=0$, even in the analytic case. To handle this, the analysis is done in polar coordinates. The s-injectivity for real analytic magnetic systems and the uniform stability estimate imply s-injectivity of the magnetic ray transform for generic $(g,\alpha)$.

In section~\ref{sec_5}, we show that $I$ is s-injective for simple magnetic systems with an explicit bound on the curvature and $\Omega$. This relies on an analog of Pestov's identity for our case, see \cite{PSh, Sh1} that goes back to \cite{Mu, MR}, see also the references there.

In section~\ref{sec_6}, we consider the non-linear magnetic rigidity problem. We prove rigidity in a given conformal class and rigidity within real analytic systems.
We also show that if a simple Riemannian manifold $(M,g)$ is boundary rigid (within the class of
Riemannian metrics), then it is also magnetic boundary rigid.
In this section we also study the local problem. We show rigidity near any $(g,\alpha)$ in the generic class $\mathcal{G}$ using the analysis of the linear problem. This does not directly follow from the implicit function theorem, and the stability estimate in Theorem~\ref{thm_estimate} plays a crucial role. There is an essential difficulty compared to the Riemannian case $\alpha=0$, \cite{SU, SU1}, since we cannot decouple $g$ and $\alpha$ in the linearization argument. This difficulty is resolved by Lemma~\ref{pr_dv}.

Section \ref{sec_7} is devoted to two-dimensional systems.
Here we prove that two-dimensional simple magnetic systems are
magnetic boundary rigidity. This generalizes the boundary rigidity theorem
of \cite{PU}. Our proof essentially resembles that in \cite{PU}, establishing a connection
between the scattering relation of a magnetic system and the Dirichlet-to-Neumann map
of the Laplace-Beltrami operator of the underlying Riemannian manifold.

\section{Boundary determination, scattering relation, and volume} \label{sec_2}

\subsection{Boundary jets of the metric and magnetic potential}
Here we show that up to gauge equivalence the boundary action function
completely determines the Riemannian metric and magnetic potential
on the boundary of the manifold under study.

\begin{Lemma}\label{bnd-0}
If $(g,\alpha)$ and $(g',\alpha')$ are simple magnetic systems
on $M$ with the same boundary action function,
then
\begin{equation}\label{bnd-data-0}
i^*g=i^*g',
\quad
i^*\alpha=i^*\alpha',
\end{equation}
where $i:\partial M\to M$ is the embedding map.
\end{Lemma}

\begin{proof}
Given $x\in\partial M$ and $\xi\in T_x(\partial M)$,
let $\tau (s)$, $-\varepsilon <s<\varepsilon $, be a curve on
$\partial M$ with $\tau (0)=x$ and $\dot \tau(0)=\xi$.
It is easy to see that
$$
\lim_{s\to0}\frac{\mathbb A(x,\tau(s))}s
=|\xi|_g-\alpha(\xi).
$$

A similar equality holds for the magnetic system  $(g',\alpha')$.
Therefore,
$$
|\xi|_{g}-\alpha(\xi)=|\xi|_{g'}-\alpha'(\xi).
$$
Changing $\xi$ to $-\xi$, we get
$$
|\xi|_{g}+\alpha(\xi)=|\xi|_{g'}+\alpha'(\xi),
$$
whence we infer \eqref{bnd-data-0}
\end{proof}

Now, we prove that the boundary action function determines the full jets
of the metric and magnetic potential on the boundary.
This generalizes the corresponding results of \cite{M, LSU}.

\begin{Theorem}\label{jet}
If $(g,\alpha)$ and $(g',\alpha')$ are simple
magnetic systems on $M$ with
the same boundary action function, then $(g',\alpha')$ is gauge equivalent
to some $(\bar g,\bar\alpha)$
such that in any local coordinate system we have
$\partial^m g|_{\partial M}=\partial^m \bar g|_{\partial M}$
and $\partial^m\alpha|_{\partial M}=\partial^m \bar\alpha|_{\partial M}$
for every multi-index $m$.
\end{Theorem}

\begin{proof}
Denote by $\nu$ the inward unit normal to $\partial M$
w.r.t.\ $g$.
The ``usual'' boundary exponential map
$\exp_{\partial M}(p,t)=\exp_p(t\nu(p))$, $p\in\partial M$, $t\ge 0$,
takes a~sufficiently small
neighborhood of the set $\partial M\times\{0\}$
in $\partial M\times\mathbb R_+$ diffeomorphically
onto some neighborhood of $\partial M$ in $M$.

Let $\nu'$ and $\exp'$ denote the corresponding objects
for the metric $g'$.

The map $\exp_{\partial M}\circ\,(\exp'_{\partial M})^{-1}$
is well defined in some neighborhood of $\partial M$ in $M$
and is the identity when restricted to $\partial M$.
We extend this map from a neighborhood of $\partial M$ in $M$
to a diffeomorphism $f:M\to M$ and
put $\bar g=f^*g'$ and $\alpha''=f^*\alpha'$.
By Lemma \ref{bnd-0},
$g$ and $\bar g$ induce the same Riemannian metric on $\partial M$
and, by construction, at each point of $\partial M$
the inward unit normal w.r.t.\ $\bar g$ coincides
with the one w.r.t.\ $g$. Therefore,
$\bar g|_{\partial M}=g|_{\partial M}$.

Next, by Lemma \ref{bnd-0}, $\alpha$ and $\alpha''$ induce
the same $1$-form on $\partial M$.
Applying Lemma 2.2 of \cite{Sh3}
to the form $\omega=\alpha-\alpha''$  we find
a function $\varphi\in C^\infty_0(M)$
such that $\omega-d\varphi$ induces the zero form on
every sufficiently short geodesic
of $\bar g$ starting from $\partial M$ in the normal direction.
We now put $\bar \alpha=\alpha''+d\varphi$.

Let us prove the equality of derivatives on the boundary.
We will use the same argument as in \cite{LSU}.
We fix $x_0\in\partial M$ and introduce boundary normal
coordinates $(x',x^n)$ w.r.t.\ $g$ near $x_0$.
By construction, the same coordinates are boundary normal
coordinates w.r.t.\ $\bar g$. Thus, the line elements $ds^2$
and $d\bar s^2$ of these metrics  are given by
$$
ds^2=g_{\iota\kappa}dx'_\iota dx'_\kappa+dx^2_n,
$$
$$
d\bar s^2=\bar g_{\iota\kappa}dx'_\iota dx'_\kappa+dx^2_n,
$$
where $\iota,\kappa$ vary from $1$ to $n-1$.
Therefore, for $h=g-\bar g$
we have $h_{in}=0$ for $i=1,\dots,n$.
Also, by the construction of $\bar\alpha$,
for $\beta=\alpha-\bar \alpha$ we have
$\beta_n=0$.

It now suffices to prove that
\begin{equation}\label{der-vanish}
\partial_n^m h_{\iota\kappa}|_{x=x_0}=0,
\quad \partial_n^m \beta_\iota|_{x=x_0}=0
\quad \mbox{for } m=0,1,\dots;\ \iota,\kappa=1,\dots,n-1.
\end{equation}

The case of $m=0$ is granted.
Assume that there is a least $m\ge1$ such that \eqref{der-vanish}
is not true. The Taylor expansion of $h$ and $\beta$ then implies that
there is a unit vector $\xi_0\in T_{x_0}(\partial M)$ such that
\begin{equation}\label{h-b}
\frac 12 h_{\iota\kappa}(x)\xi^\iota \xi^\kappa-\beta_\iota(x)\xi^\iota>0
\ \mbox{(or }<0)
\end{equation}
for $x^n>0$ and $x'$ both sufficiently small and $\xi$ close to $\xi_0$.
(Here $\iota$ and $\kappa$ vary from $1$ to $n-1$ because
$h_{1n}=h_{n1}=0$ and $\beta_n=0$.)

For arbitrary $x,y\in\partial M$,
let $\gamma=\gamma_{x,y}:[0,T]\to M$
($\bar\gamma=\bar\gamma_{x,y}:[0,\bar T]\to M$)
be the unit speed magnetic
geodesic of the system $(g,\alpha)$
(system $(\bar g,\bar\alpha)$) from $x$ to $y$.

On the one hand, since $\bar\gamma$ minimizes
the time-free action w.r.t.\ $(\bar g,\bar \alpha)$,
we have
\begin{equation*}
\mathbb A(x,y)
\le\frac12\int_0^T \bar g_{ij}(\gamma(t))
\dot\gamma^i(t)\dot\gamma^j(t)\,dt+\frac12 T
-\int_0^T\bar\alpha_i(\gamma(t))\dot\gamma^i(t)\,dt.
\end{equation*}

On the other hand, since $\gamma$ minimizes the time-free action
w.r.t.\ $(g,\alpha)$,
$$
\mathbb A(x,y)
=\frac12\int_0^T \big[g_{ij}(\gamma(t))
\dot\gamma^i(t)\dot\gamma^j(t)\,dt+\frac12 T
-\int_0^T\alpha_i(\gamma(t))\dot\gamma^i(t)\big]\,dt.
$$

Therefore,
\begin{equation}\label{int-1}
\int_0^T \left[\frac 12 h_{ij}(\gamma(t))
\dot\gamma^i(t)\dot\gamma^j(t)
-\beta_i(\gamma(t))\dot\gamma^i(t)\right]\,dt\le 0.
\end{equation}

Similarly, we derive the inequality
\begin{equation}\label{int-2}
\int_0^{\bar T} \left[\frac 12 h_{ij}(\bar\gamma(t))
\dot{\bar\gamma}^i(t)\dot{\bar\gamma}^j(t)
-\beta_i(\bar\gamma(t))\dot{\bar\gamma}^i(t)\right]\,dt\ge 0.
\end{equation}

Continuing the proof of the theorem,
we now choose $x=x_0$ and choose $y=\delta(s)$, where
$\delta:(-\varepsilon,\varepsilon)\to \partial M$ is a smooth curve
with $\delta(0)=x_0$ and $\dot\delta(0)=\xi_0$.
Then for $s>0$ sufficiently small we see that, in view of
\eqref{h-b}, we cannot simultaneously have
\eqref{int-1} and \eqref{int-2}. This contradiction
concludes the proof of the theorem.
\end{proof}

\subsection{Scattering relation}\label{stattering}
Now, we define a scattering relation and
restate our problem in terms of the scattering relation.

For $(x,\xi)\in SM$, let $\gamma_{x,\xi}:[\ell^-(x,\xi),\ell(x,\xi)]\to M$
be a magnetic geodesic such that $\gamma_{x,\xi}(0)=x$,
$\dot\gamma_{x,\xi}(0)=\xi$,
and $\gamma_{x,\xi}(\ell^-(x,\xi)),\gamma_{x,\xi}(\ell(x,\xi))\in \partial M$.
Clearly, the functions $\ell^-(x,\xi)$ and $\ell(x,\xi)$ are continuous and,
on using the implicit function theorem,
they are easily seen to be smooth near a point $(x,\xi)$ such that
the magnetic geodesic $\gamma_{x,\xi}(t)$ meets $\partial M$
transversally at $t=\ell^-(x,\xi)$ and $t=\ell(x,\xi)$ respectively.
By  \eqref{strict-conv} and Lemma \ref{convexity} in Appendix A,
the last condition holds everywhere on $SM\setminus S(\partial M)$.
Thus, $\ell^-$ and $\ell$ are smooth on $SM\setminus S(\partial M)$.

Let $\partial_+SM$ and $\partial_-SM$
denote the bundles of inward and outward unit vectors
over $\partial M$:
$$
\partial_+SM=\{(x,\xi)\in SM: x\in\partial M, \langle \xi,\nu(x)\rangle \ge0\},
$$
$$
\partial_-SM=\{(x,\xi)\in SM: x\in\partial M, \langle \xi,\nu(x)\rangle \le0\},
$$
where $\nu$ is the inward unit normal to $\partial M$.
Note that $\partial(SM)=\partial_+SM\cup \partial_-SM$
and $\partial_+SM\cap \partial_-SM=S(\partial M)$.

\begin{Lemma}[{cf. \cite[Lemmas 3.2.1, 3.2.2]{Sh2}}]\label{l-l}
For a simple magnetic system, the function $\mathbb L:\partial(SM)\to\mathbb R$,
defined by
$$
\mathbb L(x,\xi):=\begin{cases}\ell(x,\xi)&\text{if}\quad (x,\xi)\in \partial_+SM,
                \\ \ell^-(x,\xi)&\text{if}\quad (x,\xi)\in \partial_-SM,
           \end{cases}
$$
is smooth. In particular, $\ell:\partial_+SM\to\mathbb R$ is smooth.
The ratio $\frac {\mathbb L(x,\xi)}{\langle \nu(x),\xi\rangle}$
is uniformly bounded on $\partial(SM)\setminus S(\partial M)$.
\end{Lemma}

\begin{proof}
Let $\rho$ be a smooth nonnegative function on $M$
such that $\partial M=\rho^{-1}(0)$
and $|\grad\rho|=1$ in some neighborhood of $\partial M$.
Put $h(x,\xi,t)=\rho(\gamma_{x,\xi}(t))$ for $(x,\xi)\in\partial(SM)$.
Then
\begin{align*}
&h(x,\xi,0)=0,
\\
&\de ht(x,\xi,0)=\langle \nu(x),\xi\rangle,
\\
&
\frac{\partial^2 h}{\partial t^2}(x,\xi,0)
=\hess_x\rho (\xi,\xi)+\langle \nu(x),Y(\xi)\rangle.
\end{align*}

Therefore, for some smooth function $R(x,\xi,t)$,
$$
h(x,\xi,t)=\langle \nu(x),\xi\rangle t
+\frac12\left(\hess_x\rho (\xi,\xi)+\langle \nu(x),Y(\xi)\rangle\right)t^2
+R(x,\xi,t)t^3.
$$
Since $h(x,\xi,\mathbb L(x,\xi))=0$, it follows that $L=\mathbb L(x,\xi)$
is a solution of the equation
\begin{equation}\label{equation-l}
F(x,\xi,L):=\langle \nu(x),\xi\rangle
+\frac12\left(\hess_x\rho (\xi,\xi)+\langle \nu(x),Y(\xi)\rangle\right)L
+R(x,\xi,t)L^2=0.
\end{equation}
By \eqref{strict-conv}, for $(x,\xi)\in S(\partial M)$
\begin{align*}
\de FL(x,\xi,0)
&=\frac12\left(\hess_x\rho (\xi,\xi)+\langle \nu(x),Y(\xi)\rangle\right)
\\
&=\frac12\left(-\Lambda(x,\xi)+\langle \nu(x),Y(\xi)\rangle\right)<0.
\end{align*}
Now, the implicit function theorem yields smoothness of $\mathbb L(x,\xi)$
in a neighborhood of $S(\partial M)$. Since $\mathbb L$ is also smooth
on $\partial(SM)\setminus S(\partial M)$, we conclude that
$\mathbb L$ is smooth on $\partial(SM)$.

Next, from \eqref{equation-l} we get for $(x,\xi)\in \partial(SM)\setminus S(\partial M)$
$$
\left[\frac12(\hess_x\rho (\xi,\xi)+\langle \nu(x),Y(\xi)\rangle)
+R(x,\xi,t)\mathbb L(x,\xi)\right]
\frac {\mathbb L(x,\xi)}{\langle \nu(x),\xi\rangle}=-1.
$$
Again by \eqref{strict-conv} this yields boundedness of the
ratio $\frac {\mathbb L(x,\xi)}{\langle \nu(x),\xi\rangle}$
for $(x,\xi)$ sufficiently close to $S(\partial M)$ (where $\mathbb L$ is sufficiently
small), and clearly implies boundedness of the ratio on the whole
set $\partial(SM)\setminus S(\partial M)$.
\end{proof}

\begin{Definition}\label{scat-def}
The {\em scattering relation}
$\mathcal S:\partial_+SM\to\partial_-SM$ of a magnetic system
$(M,g,\alpha)$
is defined as follows:
$$
\mathcal{S}(x,\xi)=(\gamma_{x,\xi}(\ell(x,\xi)),\dot\gamma_{x,\xi}(\ell(x,\xi))).
$$

The {\em restricted scattering relation} $\mathfrak s:\partial_+SM\to M$
is defined to be the postcomposition of the
scattering relation with the natural projection
of $\partial_-SM$ to $M$, i.e.,
$$
\mathfrak s(x,\xi)=\gamma_{x,\xi}(\ell(x,\xi)).
$$
\end{Definition}

By the preceding lemma, $\mathcal S$ and $\mathfrak s$ are smooth maps.

The next lemma generalizes the well-known assertion of  \cite{M}.

\begin{Lemma}\label{action-scat}
Suppose that $(g,\alpha)$ and $(g',\alpha')$ are simple
magnetic systems on $M$ such that
$g|_{\partial M}=g'|_{\partial M}$.
If the boundary action functions $\mathbb A|_{\partial M\times \partial M}$
and $\mathbb A'|_{\partial M\times \partial M}$ of both the systems coincide,
then the scattering relations $\mathcal S$
and $\mathcal S'$  of these systems coincide,
$\mathcal S=\mathcal S'$.
\end{Lemma}

In the opposite direction, we have the following:

\begin{Lemma}\label{scat-action}
Suppose that $(g,\alpha)$ and $(g',\alpha')$ are simple
magnetic systems on $M$ such that
$g|_{\partial M}=g'|_{\partial M}$ and
$i^*\alpha=i^*\alpha'$.
If the restricted scattering relations $\mathfrak s$
and $\mathfrak s'$  of the systems
coincide, $\mathfrak s=\mathfrak s'$,
then the boundary action functions of these systems coincide,
$\mathbb A|_{\partial M\times \partial M}=\mathbb A'|_{\partial M\times \partial M}$.
\end{Lemma}

To prove these lemmas, we need one fruitful result.

\begin{Lemma}\label{der}
If $(M,g,\alpha)$ is a simple magnetic system, then for $x,y\in\partial M$
$$
\de{\mathbb A(x,y)}\xi
=-\langle \dot\gamma_{x,y}(0),\xi\rangle
+\alpha(\xi)
\quad\text{for }\xi\in T_x(\partial M),
$$
$$
\de{\mathbb A(x,y)}\eta
=\langle \dot\gamma_{x,y}(T),\eta\rangle -\alpha(\eta)
\quad\text{for }\eta\in T_y(\partial M),
$$
where $\gamma_{x,y}:[0,T]\to M$ is the unit speed magnetic geodesic
from $x$ to $y$.
\end{Lemma}

\begin{proof}
Fix $x,y\in\partial M$, $x\ne y$, and $\xi\in T_x(\partial M)$.
Let $\tau(s)$, $-\varepsilon <s<\varepsilon$,
be a curve on $\partial M$ with $\tau(0)=x$ and $\dot\tau(0)=\xi$.
For every $s$, put $\gamma(t,s)=\gamma_{\tau(s),y}(t)$,
$0\le t\le T_s$, denoting $\gamma=\gamma(t,0)$
and $T=T_0$,
and consider
$$
c(t,s)=\gamma\Big(\frac{T_s}{T}t,s\Big),
\quad 0\le t\le T.
$$
Each curve $c(\cdot,s)$ is defined on the interval $[0,T]$
and its length is exactly $T_s$. We have
$$
\de{\mathbb A(x,y)}\xi
=\frac {dT_s}{ds}(0)
-\frac d{ds}\bigg\{\int_{\gamma_s} \alpha\bigg\}\bigg|_{s=0}.
$$

Using the first variation formula for length and the fact that
$c(t,0)=\gamma(t)$, we have
$$
\frac {dT_s}{ds}(0)
=-\langle \dot\gamma(0),\xi\rangle
-\int_0^T\langle \nabla_{\dot\gamma}\dot\gamma,V\rangle\,dt,
$$
where $V(t)=\de{c}{s}(t,0)$ is the variation field of $c(t,s)$.
Using the equation of magnetic geodesics \eqref{geodesic}
and the definition \eqref{lorentz} of $Y$, we obtain
$$
\frac {dT_s}{ds}(0)=-\langle \dot\gamma(0),\xi\rangle
-\int_0^T\Omega (\dot\gamma,V)\,dt.
$$

On the other hand, it is easy to see that
$$
\frac d{ds}\left\{\int_{\gamma_s} \alpha\right\}
=-\alpha(\xi)+\int_0^T\Omega(V,\dot\gamma)\,dt.
$$

This gives the first formula of the lemma.
A similar calculation gives the second
formula.
\end{proof}

\begin{proof}[Proof of Lemma \ref{action-scat}]
By  Lemma \ref{bnd-0}, $\alpha(v)=\alpha'(v)$ for all $v\in T(\partial M$).
Using Lemma \ref{der} and the fact that both metrics are
the same on the boundary, we easily conclude that
the scattering relations are the same.
\end{proof}

\begin{proof}[Proof of Lemma \ref{scat-action}]
Take $x\in\partial M$ and define
$\mathfrak s_x:\partial_+SM\setminus S(\partial M)\to \partial M\setminus\{x\}$
by $\mathfrak s_x(\xi)=\mathfrak s(x,\xi)$.
This map is a diffeomorphism. Consider its inverse
$\mathfrak s_x^{-1}: \partial M\setminus\{x\}\to \partial_+SM$.
By Lemma \ref{der}, for every $y\in\partial M$, $y\ne x$,
we have for all $\xi\in T_x(\partial M)$
\begin{equation}\label{derive}
\de{\mathbb A(x,y)}{\xi}
=-\langle \mathfrak s_x^{-1}(y),\xi\rangle_g
+\alpha(\xi).
\end{equation}
The assumptions of the lemma imply that the right-hand side
of \eqref{derive} is the same for the second magnetic system;
therefore, so does the left-hand side.
Now, the claim of the lemma is immediate.
\end{proof}

Thus, for simple magnetic systems, the boundary rigidity problem is essentially equivalent
to the problem of restoring a Riemannian metric and
a magnetic potential from the restricted scattering relation.

\subsection{Determination of volume}
Here we show that the boundary action function determines the volume
of the manifold. This generalizes the well-known assertion of
\cite[Proposition 2.13]{M}.

\begin{Theorem}\label{volume}
If $(g,\alpha)$ and $(g',\alpha')$ are simple magnetic systems
on $M$ with the same boundary action function,
then the volume $\vol_g M$ of $M$ w.r.t.\ $g$ equals the volume $\vol_{g'}M$ of $M$
w.r.t.\ $g'$.
\end{Theorem}

\begin{proof}
By Santal\'o's formula (see \eqref{santalo-f} in Appendix A) we have
\begin{equation*}
\vol_g M=\frac 1{w_{n-1}}\int_{\partial_+SM}
\mathbb \ell(x,\xi)\,d\mu(x,\xi).
\end{equation*}
Using
\begin{equation}\label{int-a}
\int_{SM} \alpha(x,\xi)\,d\Sigma^{2n-1}(x,\xi)=0
\end{equation}
and Santal\'o's formula again, we obtain
\begin{equation}\label{vol}
\vol_g M=\frac 1{w_{n-1}}\int_{\partial_+SM}
\mathbb A(\gamma_{x,\xi})\,d\mu(x,\xi).
\end{equation}

In view of Theorem \ref{jet} and Lemma \ref{action-scat},
we may assume that the right-hand side of this equality is the same for
the magnetic system $(g',\alpha')$. This yields the sought equality
of volumes.
\end{proof}

\section{Magnetic ray transform}  \label{sec_3}

\subsection{Derivatives of the action function}
Let $(g,\alpha)$ be a simple magnetic system on $M$. It is easy to see that
there is an $\varepsilon>0$ so that
every magnetic system $(g+h,\alpha+\beta)$, satisfying
\begin{equation}\label{neigh}
\|h\|_{C^2}\le\varepsilon,\quad \|\beta\|_{C^1}\le\varepsilon,
\end{equation}
is simple.

Given $h$ and $\beta$ satisfying \eqref{neigh},
consider the $1$-parameter family $(g^s,\alpha^s)$ with
$$
g^s=g+s h,
\quad \alpha^s=\alpha+s\beta, \quad s\in [0,1].
$$
Clearly, each of these systems is simple.

\begin{Lemma}\label{var}
For $x,y\in\partial M$,
\begin{equation}\label{first-der}
\frac{d\mathbb A_{g^s,\alpha^s}(x,y)}{ds}
=\frac12\int_{\gamma_s}\langle h,\dot\gamma^2_s\rangle
-\int_{\gamma_s}\beta,
\end{equation}
where $\gamma_s$ is the unit speed magnetic geodesic from $x$ to $y$
w.r.t.\ $(g^s,\alpha^s)$.

If
\begin{equation}\label{hdm-bdm}
h|_{\partial M}=0,
\quad \beta|_{\partial M}=0,
\end{equation}
then
\begin{equation}\label{second-der}
\left|\frac{d^2\mathbb A_{g^s,\alpha^s}(x,y)}{ds^2}\right|
\le C (\|h\|^2_{C^1}+\|\beta\|^2_{C^1}),
\end{equation}
with a constant $C$ independent of $h$ and $\beta$
and $C^2$ locally uniform in $(g,\alpha)$.
\end{Lemma}

\begin{proof}
Define
$$
\varphi(s,\tau):=\mathbb A_{g^\tau,\alpha^\tau}(\gamma_s)
=\frac12\int_0^{T_s}|\dot\gamma_s(t)|^2_{g^\tau}\,dt
+\frac12 T_s-\int_{\gamma_s}\alpha^\tau.
$$
Then
\begin{equation}\label{da-ds}
\frac{d\mathbb A_{g^s,\alpha^s}(x,y)}{ds}
=\de{\varphi}s(s,s)+\de{\varphi}\tau (s,s).
\end{equation}

By Lemma \ref{minimize}, unit speed magnetic geodesics minimize
the time-free action; therefore,
for a fixed $\tau$,
$\mathbb A_{g^\tau,\alpha^\tau}(\gamma_s)$
has a minimum at $s=\tau$, which yields
\begin{equation}\label{dfs}
\de{\varphi}s(\tau,\tau)=0.
\end{equation}
Next,
\begin{equation}\label{dft}
\de{\varphi}\tau
=\frac12\int_0^{T_s}\left(\de{}{\tau}|\dot\gamma_s(t)|^2_{g^\tau}\right)\,dt
-\int_{\gamma_s}\frac{d}{d\tau}\alpha^\tau
=\frac12\int_{\gamma_s}\langle h,\dot\gamma^2_s\rangle
-\int_{\gamma_s}\beta.
\end{equation}
Combining \eqref{da-ds}--\eqref{dft} gives \eqref{first-der}.

Differentiating \eqref{first-der} and using \eqref{hdm-bdm},
we obtain
\begin{multline*}
\frac{d^2\mathbb A_{g^s,\alpha^s}(x,y)}{ds^2}
=\int_0^{T_s}\left\{\frac12\left(\nabla_{\gamma'_s} h_{ij}\right)
\dot \gamma_s^i\dot \gamma_s^j
+h_{ij}\dot \gamma_s^i\left(\nabla_{\gamma'_s}\dot\gamma_s^j\right)\right.
\\
-\left.\left(\nabla_{\gamma'_s} \beta_i\right)\dot \gamma_s^i
-\beta_i\left(\nabla_{\gamma'_s}\dot\gamma_s^i\right)\right\}\,dt,
\end{multline*}
where $\gamma'_s=\partial \gamma_s/\partial s$.
Whence
\begin{equation}\label{est-d2a}
\left|\frac{d^2\mathbb A_{g^s,\alpha^s}(x,y)}{ds^2}\right|
\le C\left(\|h\|_{C^1}
+\|\beta\|_{C^1}\right)\left(\|{\gamma'_s}\|_C+\|\nabla_{\gamma'_s}\dot\gamma_s\|_C\right).
\end{equation}

By the equation of magnetic geodesics, we have
\begin{equation}\label{eq-s}
{\overset s\nabla}_{\dot\gamma_s}\dot\gamma_s=\overset sY(\dot\gamma_s),
\end{equation}
where $\overset s\nabla$ stands for the covariant derivative
related to $g^s$, and $\overset sY$ is the Lorentz force
associated with $g^s$, $\alpha_s$.
The initial conditions are given by
\begin{equation}\label{in-s}
\gamma_s(0)=x,
\quad \dot\gamma_s(0)=\big({\overset s \exp}{}^\mu_x\big)^{-1}(y),
\end{equation}
where $\overset s \exp{}^\mu_x$ stands for the magnetic exponential map
associated with $g^s$, $\alpha_s$.
From \eqref{eq-s} and \eqref{in-s} we easily infer that
\begin{equation}\label{est-s}
\|{\gamma'_s}\|_C+\|\nabla_{\gamma'_s}\dot\gamma_s\|_C
\le C \left(\|h\|_{C^1} +\|\beta\|_{C^1}\right).
\end{equation}

Combining \eqref{est-d2a} and \eqref{est-s} leads to \eqref{second-der}.
\end{proof}

\subsection{Magnetic ray transform}
Let $(M,g,\alpha)$ be a simple magnetic system and
$\phi:SM\to \mathbb R$ a smooth function on the unit sphere bundle.
We define the {\em magnetic ray transform} of $\phi$ to be
the following function on the space of unit speed magnetic
geodesics going from a boundary point to a boundary point:
$$
I\phi(\gamma)=\int_\gamma\phi:=\int_0^T\phi(\gamma(t),\dot\gamma(t))\,dt,
$$
where $\gamma:[0,T]\to M$ is any unit speed magnetic geodesic such
that $\gamma(0)\in\partial M$ and $\gamma(T)\in\partial M$.
Assuming that the magnetic geodesics are parametrized
by $\partial_+ SM$, we obtain a map
$I : C^\infty(SM) \to C(\partial_+ SM)$,
\begin{equation}\label{mag-ray}
I\phi(x,\xi)=\int_0^{\ell(x,\xi)}\phi(\psi^t(x,\xi))\,dt,
\quad (x,\xi)\in \partial_+ SM.
\end{equation}

In the space of real-valued functions on
$\partial_+SM$ define the norm
$$
\|\phi\|^2
=\int_{\partial_+SM}\phi^2\,d\mu
$$
and the corresponding  inner product.
Here $d\mu(x,\xi)=\langle \xi,\nu(x)\rangle \,d\Sigma^{2n-2}$
(see \ref{santalo-s} in Appendix A).
Denote the corresponding Hilbert space by $L^2_\mu(\partial_+SM)$.

\begin{Lemma}\label{L2}
The operator $I$ extends to a bounded operator
$$
I:L^2(SM)\to L^2_\mu(\partial_+SM).
$$
\end{Lemma}

\begin{proof}
Indeed, it is easy to see that $(I\phi)^2\le C I\phi^2$, with some constant $C$
independent of $\phi$. Therefore,
$$
\int_{\partial_+SM}(I\phi)^2\,d\mu
\le C \int_{\partial_+SM} I\phi^2\,d\mu
=C\int_{SM}\phi^2\,d\Sigma^{2n-1}
$$
by the Santal\'o formula \eqref{santalo-f}.
\end{proof}

\subsection{Solenoidal and potential pairs}\label{solenoid-potential}
For a decomposition of a symmetric tensor field into a potential and
a solenoidal part (relevant when $\Omega=0$), we refer to \cite{Sh1, Sh2}.

In view of the linearization formula \eqref{first-der},
we are  mainly interested in $I$ applied to functions
$\phi$ that are of the form
\begin{equation}  \label{S0}
\phi(x,\xi) = h_{ij}(x) \xi^i\xi^j + \beta_j(x)\xi^j
\end{equation}
(and, more generally, that are polynomials in $\xi$).
Then, given a symmetric 2-tensor $h$ and a $1$-form $\beta$,  we set
for $(x,\xi) \in\partial_+ SM$
\begin{equation}  \label{S1}
I[h,\beta] (x,\xi) = \int_0^T h_{ij}(\gamma(t)) \dot \gamma^i(t) \dot\gamma^j(t) \,dt +  \int_0^T \beta_j(\gamma(t)) \dot\gamma^j(t)\, dt,
\end{equation}
where $\gamma = \gamma_{x,\xi}$,  $T=\ell(x,\xi)$.

If $F$ is a notation for a function space
($C^k$, $L^p$, $H^k$, etc.),
then we will denote by
$\mathbf F(M)$
the corresponding space of pairs $\mathbf f=[h,\beta]$,
with $h$ a symmetric covariant $2$-tensor and $\beta$ a $1$-form,
and denote by $\mathcal F(M)$ the corresponding space of
pairs $\mathbf{w} = [v,\varphi]$, with $v$ a $1$-form
and $\varphi$ a function on $M$.
In particular, $\mathbf{L}^2(M)$ is the space of square integrable
pairs $\mathbf f=[h,\beta]$, and we endow this space with the norm
\begin{equation}\label{norm}
\|\mathbf f\|^2
=\int_M\left\{|h|_g^2+\frac{n-1}{2}|\beta|_g^2\right\}\,d\vol,
\end{equation}
with the corresponding inner product.
(The choice of the factor $(n-1)/2$ will play its role in the proof
of Theorem \ref{quadratic}.) In the space $\mathcal L^2(M)$ we will consider
the norm
\begin{equation}  \label{S7}
\|\mathbf{w}\|^2 = \int_M \left(|v|_g^2 + \varphi^2 \right)\,d\vol.
\end{equation}

Clearly, the norm of a pair $[h,\beta]$ in $\mathbf{L}^2(M)$ is equivalent to
the norm of the corresponding polynomial \eqref{S0} in $L^2(SM)$. Therefore,
Lemma \ref{L2} implies that
$$
I : \mathbf{L}^2(M) \to L^2_\mu(\partial_+ SM)
$$
and it is bounded.

Define $Y:T^*M\to T^*M$ by
$$
Y(\eta)=-(Y^j_i\eta_j),\quad \eta=(\eta_j)\in T^*M.
$$
Note that this definition agrees with the map $Y:TM\to TM$
and the isomorphism between the tangent and cotangent bundles via
the Riemannian metric.

Clearly, $I$ vanishes on functions
$\phi(x,\xi)=\mathbf{G}_\mu\psi(x,\xi)$ if $\psi$ vanishes for
$x\in \partial M$. To find a class of such functions that are
polynomials of $\xi$ of degree at most $2$,  assume that
$$
\mathbf G_\mu [v_i(x)\xi^i+\varphi(x)]=h_{ij}(x)\xi^i\xi^j +\beta_j(x)\xi^j.
$$
Since
\begin{align}\nonumber
\mathbf G_\mu (v_i(x)\xi^i+\varphi(x))
&
=\xi^j[v_{i,j}\xi^i+v_iY^i_j+\varphi_{,j} ]
\\  \label{S1aa}
&=(\d v)_{ij}\xi^i\xi^j+(\varphi_{,j}-Y(v)_j)\xi^j,
\end{align}
we get
\begin{equation}\label{p-eq}
h=\d v,\quad \beta=d\varphi-Y(v),
\end{equation}
where $\d v$ is the symmetric differential of $v$.
We used here the fact that the even part of \eqref{S1aa} w.r.t.\ $\xi\in S_xM$
determines the quadratic form $(\d v)_{ij}\xi^i\xi^j$,
while the odd part determines $\beta(\xi)$.
Next, knowing $(\d v)_{ij}\xi^i\xi^j$ for $n(n+1)/2$
generic $\xi\in S_xM$
is enough to recover $v$;
similarly knowing $\beta(\xi)$ for $n$ linear independent $\xi\in S_xM$
is enough to recover $\beta$. We have the following stronger statement.

\begin{Lemma}   \label{lemma_S_det}
Fix $x\in M$. Given any open subset $V$ of $S_xM$,
there exist $N= n(n+1)/2+n$ vectors $\xi_m\in V$, $m=1,\dots,N$,
such that $[h,\beta]$ is uniquely determined by
the values of $\psi(x,\xi) :=h_{ij}(x)\xi^i\xi^j +\beta_j(x)\xi^j$
at $\xi=\xi_m$, $m=1,\dots,N$.
\end{Lemma}

\begin{proof}
Since $x$ is fixed, we denote $\psi(\xi) =\psi(x,\xi)$.
We show first that $\psi(\xi)$, known for all $\xi\in S_xM$,
determines $[h,\beta]$. Since $\psi(\xi)$ is a linear functional of
$[h,\beta]$, and the latter belongs to a linear space
that can be identified with $\mathbb{R}^N$, it is enough to show
that $\psi(\xi) =0$ implies $h=0$, $\beta=0$.
By replacing $\xi$ by $-\xi$, we get that
$h_{ij}\xi^i\xi^j=0$, $\beta_j\xi^j=0$.
The second relation easily implies that $\beta=0$.
The first one implies $h=0$ easily as well
(see, e.g., \cite{Sh2}, where also a sharp estimate is established).

Next, assume that $\psi(\xi)=0$ in $V$. Then $\psi(\xi)=0$ on $S_xM$
by analytic continuation, thus  $[h,\beta]=0$.
Finally, since for any fixed $\xi$, $\psi(\xi)$ is a linear functional
belonging to $(\mathbb{R}^N)' \cong \mathbb{R}^N$,
and $\{\psi(\xi);\; \xi\in V\}$ is a complete set,
there exists a basis $\psi(\xi_k)$, $k=1,\dots,N$, in it.
\end{proof}

\begin{Remark}
We can say a bit more. Since the determination of $[h,\beta]$
is done by inverting a linear transform, one can choose $\xi_m$
continuously depending on $x$, if $x$ belongs to a small enough set $X$,
such that the map $\{\psi(x,\xi_m(x)),\; m=1,\dots,N\} \mapsto [h,\beta]$
is invertible with a uniform bound on the inverse.
Then this can be extended to compact sets $X$.
\end{Remark}

\begin{Definition}\label{p-def}
We call a pair $[h,\beta]\in \mathbf{L}^2(M)$ {\em potential}
if the equations \eqref{p-eq} hold with $[v,\varphi]\in \mathcal H_0^1(M)$.
\end{Definition}

This can be written as follows:
$$
\begin{pmatrix} h\\
               \beta \end{pmatrix}
=\mathbf d\begin{pmatrix} v\\
                         \varphi \end{pmatrix},
\quad \mathbf d :=\begin{pmatrix} \d&0\\
               -Y&d\end{pmatrix},
$$
where $\d$ stands for the symmetric derivative acting on covector fields
and $d$ is the usual differential acting on functions
(which coincides with the symmetric differential on functions).

Clearly, potential pairs satisfy
\begin{equation}  \label{S3a}
I\left( \mathbf{d}[v,\varphi] \right) =0.
\end{equation}
This follows from \eqref{S1aa} if $v$, $\varphi$ are smooth,
and follows by continuity for general $v$, $\varphi$, once we establish the mapping properties of
$N=I^*I$ below.

We will relate potential pairs to the non-linear problem.
Let $(g,\alpha)$ and $(g',\alpha')$ be two gauge equivalent pairs, i.e.,
$g'=f^*g$ and $\alpha'=f^*\alpha+d\varphi$ with some diffeomorphism
$f :M\to M$, fixing $\partial M$, and some function $\varphi$
vanishing on $\partial M$.
Linearize this near $f= \text{\rm Id}$ and $\varphi=0$.
In other words, let $f_\tau$ be a smooth family of
such diffeomorphisms with $f_0 = \text{\rm Id}$ and
let $\varphi^\tau$ be a smooth family of such functions with $\varphi^0=0$.
Let $g^\tau = f_\tau^*g$, $\alpha^\tau = f_\tau^*\alpha +d\varphi^\tau$,
and we will compute the derivatives at $\tau=0$.
It is well known \cite[(3.1.5)]{Sh2} and is easy to calculate that
$dg^\tau/d\tau |_{\tau=0}= 2\d v$,
where $v=df_\tau/d\tau|_{\tau=0}$.
Let $d\varphi^\tau/d\tau|_{\tau=0}=\psi$.
Since
$\alpha^\tau_i
= (\alpha_j\circ f_\tau) \partial f_\tau^j/\partial x^i
+\partial\varphi^\tau/\partial x^i$,
we get
\begin{align*}
\beta_i := \frac{d}{d\tau}\Big|_{\tau=0} \alpha^\tau_i&
=  v^j\partial_j\alpha_i + \alpha_j\partial_i v^j+\psi_{,i}\\
& = v^j\alpha_{i,j} + \alpha_j v^j_{,i} +\psi_{,i}\\
& = v^j\alpha_{i,j} + (\alpha_j v^j)_{,i}- v^j \alpha_{j,i}+\psi_{,i}\\
&= -(d\alpha)_{ij}v^j  + (\alpha_j v^j+\psi)_{,i}.
\end{align*}

Using \eqref{d-alpha}, \eqref{lorentz},
and treating $v$ as a $1$-form (by lowering the index), we obtain
$$
\beta=Y(v) +d(\langle \alpha,v\rangle+\psi).
$$
This shows that
\begin{equation}    \label{S3bb}
\frac{d}{d\tau}\bigg|_{\tau=0} \left[\frac12 g^\tau, -\alpha^\tau\right]
= \mathbf{d}[v,-\langle \alpha,v\rangle-\psi],
\quad \text{where $v=\frac {df_\tau}{d\tau}\Big|_{\tau=0}$,
$\psi = \frac{d\varphi^\tau}{d\tau}\Big|_{\tau=0} $.}
\end{equation}
Since $\psi$ can be arbitrary
(as long as it vanishes on $\partial M$),
we see that the linearization of $(g',\alpha')$ near
$f=\text{\rm Id}$, $\varphi=0$ is given by
$\mathbf{d}[v,\psi]$ with  $v$, $\psi$ that can be arbitrary.

For future references, let us mention that by \eqref{S1aa}
\begin{equation}  \label{S2a}
\mathbf{G}_\mu\langle \mathbf{w}, \xi\rangle = \langle \mathbf{dw},\xi\rangle  \quad  \Longrightarrow \quad
\frac{d}{dt} \langle  \mathbf{w}(\gamma),\dot\gamma \rangle = \langle \mathbf{d} \mathbf{w} (\gamma),\dot\gamma \rangle
\end{equation}
for any unit speed magnetic geodesic $\gamma(t)$,
where $\mathbf{w} = [v,\varphi]$,
and we used the notation
$\langle \mathbf{f}  ,\xi \rangle
= \langle h, \xi^2\rangle + \langle \beta, \xi\rangle$,
where $\mathbf{f} = [h,\beta]$, and similarly
$\langle \mathbf{w} ,\xi \rangle = \langle v,\xi\rangle +\varphi$.

\begin{Definition}  \label{def_S1}
We say that $I$ is s-injective if $I[h,\beta]=0$ with
$[h,\beta] \in\mathbf{L}^2(M)$ implies that $[h,\beta]$ is potential, i.e.,
\begin{equation}  \label{S4}
h = \d v, \quad \beta = -Y(v) +d\varphi
\end{equation}
with  $[v,\varphi]\in \mathcal H_0^1(M)$.
\end{Definition}

Let $[h,\beta]$ be orthogonal to all potential pairs. Then
$$
\int_M\left\{( h,\d v)
+\frac{n-1}{2} ( \beta,-Y(v)+d\varphi)\right\}\, d\vol=0
$$
for all $v,\varphi$ vanishing on $\partial M$. Then
$$
-\int_M\left\{\big( \delta h-\frac{n-1}{2} Y(\beta),v\big)
+\frac{n-1}{2} (\delta\beta)\varphi  \right\}\,d\vol=0,
$$
where $\delta$ is the divergence. Therefore,
\begin{equation}\label{s-eq}
\delta h-\frac{n-1}{2} Y(\beta)=0,\quad \delta\beta=0.
\end{equation}

\begin{Definition}\label{s-def}
We call a pair $[h,\beta]$
{\em solenoidal}
if the equations \eqref{s-eq} hold.
\end{Definition}

This can be written as
$$
\delt\begin{pmatrix} h\\
                               \beta\end{pmatrix}=0,
\quad \delt :=\begin{pmatrix} \delta& -\frac{n-1}{2}Y\\
                                     0&\delta\end{pmatrix}.
$$
Then $\mathbf{d} =-\delt^*$.

In terms of the operators $\delt$, $\mathbf{d}$,
the solenoidal pairs are defined as the ones in $\text{Ker}\, \delt$,
and the orthogonal complement of $\text{Ker}\, \delt$ is
$\text{Ran}\, \mathbf{d}$, consisting of potential pairs.
Next, we will describe the projections to solenoidal pairs.

We will show first that the operator $-\delt\mathbf{ d} $ is
an elliptic second order (and clearly, a formally self-adjoint) operator,
acting on pairs $\mathbf{w} = [v,\varphi]$,
where $v$ is a $1$-form and $\varphi$ is a function on $M$.
First, notice that
\begin{equation}  \label{S6}
(-\delt\mathbf{ d} \mathbf{w},\mathbf{w} ) = \|\mathbf{dw}\|^2, \quad \mathbf{w} \in \mathcal C_0^\infty(M^{\rm int});
\end{equation}
thus, $-\delt\mathbf{ d} $ is a non-negative operator.
Let $\sigma_p(P)$ stand for the principal symbol of~$P$.
Then (\ref{S6}) implies that for any such $\mathbf{w}$ and a fixed $x$
\begin{equation}  \label{S8}
(-\sigma_p(\delt\mathbf{ d} ) \mathbf{w},\mathbf{w} ) = \|\sigma_p(\mathbf{d}) \mathbf{w})\|^2,
\end{equation}
where the inner product and the norm are in
the finite dimensional space to which $\mathbf{w}$ belongs, i.e.,
they are as in (\ref{S7}) before the integration.
We will show that $-\sigma_p(\delt\mathbf{ d} )$ is in fact
a positive matrix-valued symbol for $\xi\not=0$. Since it is homogeneous
in the dual variable $\xi$, it is enough to show that
$-\sigma_p(\delt\mathbf{ d} ) \mathbf{w}=0$
implies $\mathbf{w}=0$ for $\xi\not=0$.
In view of (\ref{S8}), if  $-\sigma_p(\delt\mathbf{ d} ) \mathbf{w}=0$,
then $\sigma_p(\mathbf{d}) \mathbf{w}=0$.
Then $\sigma_p(\d) v=0$ (see \eqref{S23}), $\xi \varphi=0$,
where $[v,\varphi]=\mathbf{w}$.
It is well known \cite{Sh1,Sh2}, and easy to see directly, that
$\d$ is elliptic,
therefore $v=0$. Next, $\xi \varphi=0$ implies
$\varphi=0$ as well.

Note that the Dirichlet boundary conditions are coercive for $\delt\mathbf{d}$,
because the latter is positive elliptic. We will show that the kernel and the cokernel of this elliptic problem are trivial. If $\mathbf{w}$ belongs to the kernel of  $\delt\mathbf{d}$, and satisfies Dirichlet boundary conditions, it has to be smooth and then (\ref{S6}) holds for it as well. Then $\mathbf{dw}=0$ and $\mathbf{w}=0$ on $\partial M$. This easily implies $\mathbf{w}=0$ in $M$ by integrating (\ref{S2a}) along magnetic geodesics connecting boundary and interior points, and using Lemma~\ref{lemma_S_det}. Using standard arguments, one easily checks that the cokernel of $\delt\mathbf{d}$ equipped with Dirichlet boundary conditions, is trivial as well. Indeed, fix $\mathbf{u}$ in the cokernel. Since $\mathbf{u}$ is orthogonal to $\delt\mathbf{dw}$ for all $\mathbf{w}
\in \mathcal{C}_0^\infty ( M^{\rm int})$, we get  that $\delt\mathbf{du}=0$. This in particular implies that $\mathbf{u}$ has a trace on $\partial M$, and choosing  $\mathbf{w}$'s with a dense set of normal derivatives on $\partial M$, we show that this trace vanishes. Therefore, $\mathbf{u}=0$ by what we proved above.
Denote by $(\delt\mathbf{ d} )_D$ the Dirichlet realization
of $\delt\mathbf{ d} $ on $M$. The arguments above show that $(\delt\mathbf{ d} )_D$ is an invertible self-adjoint operator.

We define the following projections
\begin{equation}  \label{S9}
\mathcal{P} = \mathbf{d}  (\delt\mathbf{ d} )_D^{-1} \delt,
\quad \mathcal{S} = \text{\rm Id} - \mathcal{P}.
\end{equation}
Then $\mathbf f^s=\mathcal{S}\mathbf f$
is solenoidal, $\mathcal P\mathbf f=\mathbf d \mathbf w$
is potential, and we have the orthogonal decomposition
$$
\mathbf f=\mathbf f^s+\mathbf d \mathbf w
$$
into the solenoidal and potential parts.

Since $I$ vanishes on $\mathcal{P}\mathbf{L}^2(M)$,
the s-injectivity of $I$ is then equivalent to the following:
$I$ is injective on $\mathcal{S}\mathbf{L}^2(M)$.

\subsection{Adjoint of $I$}
For a fixed simple $(g,\alpha)$,
consider
$$
I : \mathbf{L}^2(M) \to L^2_\mu(\partial_+ SM)
$$
and consider its dual
$$
I^* :L^2_\mu(\partial_+ SM)\to  \mathbf{L}^2(M).
$$

We will now find an expression for $I^*$.
Let $\mathbf{f} = [h,\beta]$ and $\phi(x,\xi)$ be as in (\ref{S0}).
Let $\psi(x,\xi) \in C(\partial_+ SM)$.
Then
\begin{multline*}
(I\mathbf{f}, \psi) =  \int_{\partial_+SM}\psi(x,\xi)\,d\mu(x,\xi)
\\
\times\int_0^{\ell(x,\xi)}
\left[ h_{ij}(\gamma_{x,\xi}(t)) \dot \gamma_{x,\xi}^i(t) \dot \gamma_{x,\xi}^j(t) +\beta_i(\gamma_{x,\xi}(t) )\dot \gamma_{x,\xi}^i(t)  \right]\,dt.
\end{multline*}

By Santal\'o's formula \eqref{santalo-f} of Appendix A, we get
$$
(I\mathbf{f}, \psi) = \int_{SM}\left[ h_{ij}(x)\xi^i\xi^j +\beta_i(x)\xi^i \right] {\psi^\sharp (x,\xi)}\, d\Sigma^{2n-1}(x,\xi),
$$
where $\psi^\sharp(x,\xi)$ is defined as the function that is constant
along the orbits of the magnetic flow and that equals $\psi(x,\xi)$
on $\partial_+ SM$.
Let $d\sigma_x(\xi)$ be the measure on $S_xM$. Then
\begin{align*}
(I\mathbf{f}, \psi) = &\int_{M}  h_{ij}(x)\int_{S_xM}\xi^i\xi^j {\psi^\sharp (x,\xi)}\,d\sigma_x(\xi) \,d\vol(x)\\
&+\int_{M} \beta_i(x)\int_{S_xM}\xi^i{\psi^\sharp (x,\xi)}\,d\sigma_x(\xi) \,d\vol(x).
\end{align*}
Therefore, see also \cite{Sh3},
\begin{equation}  \label{S12}
I^*\psi = \Big[ \int_{S_xM}\xi^i\xi^j {\psi^\sharp (x,\xi)}\,d\sigma_x(\xi), \frac{2}{n-1} \int_{S_xM}\xi^i {\psi^\sharp (x,\xi)}\,d\sigma_x(\xi) \Big].
\end{equation}

\subsection{Integral representation for the normal operator}\label{represent-n}
Let $M_1\supset M$ be another manifold with boundary
(a domain in $\mathbb{R}^n$, actually)
such that $M_1^\text{int}\supset M$. Extend $g$, $\alpha$ to $M_1$.
Then $(M_1,g,\alpha)$ is still simple if $M_1$ is close enough to $M$.

Choose $x_0\in M_1\setminus M$, and consider the map
\[
(\exp^\mu_{x_0})^{-1} :M_1 \to (\exp^\mu_{x_0})^{-1} (M_1).
\]
It is $C^\infty$ away from $x=x_0$, therefore
it is a $C^\infty$ diffeomorphism to its image,
if we restrict it to any compact submanifold $M_{1/2}$ of $M_1^\text{int}$.
We can assume that $M_{1/2}$ has a smooth boundary, and that
$M_{1/2}^\text{int}\supset M$.
We denote $M_{1/2}$ by $M_1$ again, and the map above gives us
global coordinates in a neighborhood of $M_1$ that can be identified
with a subset of $\mathbb{R}^n$ with a smooth boundary. So in this section,
$M_1$ is considered as a subset of $\mathbb{R}^n$.
If $g$, $\alpha$ are analytic, and $M$ is analytic,
then we will also assume that $\partial M_1$ is analytic.

Denote by $I_1$ the magnetic ray transform on $\mathbf L^2(M_1)$,
$$
I_1:\mathbf L^2(M_1)\to L^2_\mu(\partial_+SM_1),
$$
and denote by $I^*_1$ its dual.
Set
\begin{equation}  \label{S10}
N = I^*_1I_1.
\end{equation}

Extending all tensors as 0 on $M_1$, we consider $\mathbf L^2(M)$
as a subspace of $\mathbf L^2(M_1)$.
In this way we may consider $N$ as the operator
$$
N : \mathbf{L}^2(M) \to \mathbf{L}^2(M_1).
$$
We then define s-injectivity of $N$ as in Definition~\ref{def_S1}.

\begin{Lemma}  \label{lemma_SN}
$N : \mathbf{L}^2(M) \to \mathbf{L}^2(M_1)$ is s-injective
if and only if $I$ is s-injective.
\end{Lemma}

\begin{proof}
Let $N$ be s-injective. Assume that $I\mathbf{f}=0$
with $\mathbf{f}\in \mathbf{L}^2(M)$. Then $I_1\mathbf{f}=0$
on $\partial_+SM_1$ as well,
therefore $N\mathbf{f}=0$, which implies that $\mathbf{f}$
is potential.

Now, assume that $I$ is s-injective and let $N\mathbf{f}=0$. Then
\[
0= (N\mathbf{f},\mathbf{f})_{\mathbf{L}^2(M_1)}
= \| I_1\mathbf{f}\|^2_{L^2_\mu(\partial_+ SM_1)},
\]
therefore, $I\mathbf{f}=0$ as an element of
$L^2_\mu(\partial_+ SM)$ as well, hence $\mathbf{f}$
is potential.
\end{proof}

We will find an integral representation of $N$
and will show that it is a $\Psi$DO in a neighborhood of $M$,
following \cite{SU,SU1,SSU}.

Using \eqref{S12}, and replacing $\xi$ by $v$,
we get the following.

\begin{Proposition}  \label{S_prop_N}
\begin{equation}  \label{S13}
N\mathbf{f} = \left[N_{22} h+ N_{21}\beta,  \frac{2}{n-1} (N_{12}h +N_{11}  \beta )\right],
\end{equation}
where
\begin{equation} \label{SN}
\begin{split}
(N_{11}\beta)^{i'}(x) &=  \int_{S_xM_1}v^{i'}\, d\sigma_x(v) \int \beta_{i}(\gamma_{x,v}(t)) \dot \gamma_{x,v}^i(t) \,dt,\\
(N_{12}h)^{i'} (x)&=  \int_{S_xM_1}v^{i'}\, d\sigma_x(v)  \int h_{ij}(\gamma_{x,v}(t)) \dot \gamma_{x,v}^i(t) \dot \gamma_{x,v}^j(t)\,dt,\\
(N_{21}\beta)^{i'j'}(x) &=  \int_{S_xM_1}v^{i'}v^{j'}\, d\sigma_x(v) \int \beta_{i}(\gamma_{x,v}(t)) \dot \gamma_{x,v}^i(t) \,dt,\\
(N_{22}h)^{i'j'}(x) &=  \int_{S_xM_1}v^{i'} v^{j'} \, d\sigma_x(v) \int h_{ij}(\gamma_{x,v}(t)) \dot \gamma_{x,v}^i(t) \dot \gamma_{x,v}^j(t)\,dt.
\end{split}
\end{equation}
\end{Proposition}

In each of the integrals above, split the $t$-integral into two parts: $I_+$ corresponding to $t>0$, and $I_-$, corresponding to $t<0$. In $I_+$, make the change of variables $y=\gamma_{x,v}(t) = \exp_x^\mu(tv)$.
Then $tv = (\exp_x^\mu)^{-1}(y)$, and $t,v$ are $C^\infty$
functions of $x$, $y$ away from the diagonal $x=y$.
We treat $I_-$ in a similar way by making the substitution $t'=-t$, $v'=-v$.
To this end, note that the simplicity assumption implies that the map
\begin{equation}  \label{S_exp}
\exp_x^{\mu,-}(tv)  = \pi\circ\psi^{-t}(-v), \quad t\ge0,\; v\in S_xM_1.
\end{equation}
is also a $C^1$ diffeomorphism for every $x\in M_1$.

We  analyze the Schwartz kernel of $N$ near the diagonal below,
but at this point we just want to emphasize
that it is smooth away from the diagonal and
therefore $N$ has the pseudolocal property.

\section{Analysis of the operator $N$. Generic s-injectivity}\label{analysis-N}

In this section, we analyze $N$ and prove s-injectivity for analytic $(g,\alpha)$, and for  generic ones.

Here we first assume that $[g,\alpha]\in \mathbf C^k(M)$.
As a result, various objects related to $g$, $\alpha$
will have smoothness $l(k)$ such that $l(k)\to\infty$ as $k\to\infty$.
To simplify the exposition, we will not try to estimate $l(k)$
(actually, $l(k)=k-k_0$ with $k_0$ depending on $n$ only).
We will say that a given function (or tensor field)
$f$ is smooth, or that $f\in C^l$,
if such an $l=l(k)$ exists, and $l$ may vary from line to line.
In particular, if $k=\infty$, then $l=\infty$.

At the end of the section we analyze the case in which
$(g,\alpha)$ is an analytic magnetic system.

\subsection{More about the magnetic exponential map at the origin}
In order to study the singularities of the kernel of $N$ near $x=y$,
we need more precise information about the exponential map at the origin.
One of the interesting features of the magnetic problem is that
the magnetic exponential map is not $C^2$ unless $\Omega=0$
by Lemma~\ref{exp-map}. On the other hand, in polar coordinates,
it is a smooth map. We are therefore forced to work in polar coordinates.

Consider $y=\exp_x^{\mu}(tv)= \gamma_{x,v}(t)$,
where $v\in S_xM_1$, and $t\ge0$ are  such that $y\in M_1$.
Recall that $M_1$ now is a subdomain of $\mathbb{R}^n$.
We are interested in the behavior of $y$ for small $|t|$ near a fixed $x_0$,
therefore, we can assume that we work in
$U = \{x;\; |y-x_0|\le \varepsilon\}$
with $0\le \varepsilon\ll1$ such that $U$ is strictly
convex w.r.t.\ the magnetic  geodesics as well as w.r.t.\ the Euclidean metric.

By Lemma \ref{exp-map},
the map $tv=\xi\mapsto y$ is not $C^2$ in general.
On the other hand,
the map $(t,v) \mapsto y$ is smooth, and respectively analytic,
if $g$ is analytic, and in fact extends as
a smooth/analytic map for small negative $t$ as well,
by the formula  $y=\gamma_{x,v}(t)$.
Similarly, the function
$m(t,v;x)= (\gamma_{x,v}(t) -x )/t$
has the same smoothness, therefore
\begin{equation}   \label{S16}
\gamma_{x,v}(t)-x =t m(t,v;x), \quad m(0,v;x) = v.
\end{equation}

We introduce the new variables $(r,\omega) \in \mathbb{R}\times S_xU$  by
\begin{equation}   \label{S17}
r = t|m(t,v;x)|_g, \quad \omega = m(t,v;x)/|m(t,v;x)|_g.
\end{equation}
Then $(r,\omega)$ are polar coordinates for $y-x =r\omega$
in which we allow $r$  to be ne\-gative.
Clearly,  $(r,\omega)$ are smooth/analytic at least for
$\varepsilon$ small enough, if $g$ is \nobreak{smooth}/analytic.
Consider the Jacobian of this change of variables
$$
J := \det \frac{\partial(r,\omega)}{\partial(t,v)}.
$$
It is not hard to see that $J|_{t=0}=1$,
therefore the map
$\mathbb{R}\times S_xU \ni (t,v) \mapsto (r,\omega)\in \mathbb{R}\times S_xU$
is a local $C^l$, respectively analytic, diffeomorphism from
(a neighborhood of $0)\times S_xU$ to its image.
We can decrease $\varepsilon$ if needed to ensure that it is
a (global) diffeomorphism on its domain because then it is clearly injective.
We denote the inverse functions by $t=t(r,\omega)$, $v=v(r,\omega)$.
Note that in the $(r,\omega)$ variables
\begin{equation}   \label{S18}
t=r +O(|r|), \quad v = \omega+O(|r|),\quad \dot \gamma_{x,v}(t)
= \omega +O(|r|).
\end{equation}

Another representation of the new coordinates can be given by writing
\[
\text{Exp}_x^\mu(t,v) = \gamma_{x,v}(t).
\]
Then
\[
r=\text{sign}(t)\left| \text{Exp}_x^\mu(t,v)  -x \right|_g,
\quad
\omega = \text{sign}(t)\frac{\text{Exp}_x^\mu(t,v)  -x}
{\left| \text{Exp}_x^\mu(t,v)  -x \right|_g},
\]
and
\[
(t,v) = \left( \text{Exp}_x^\mu \right)^{-1}(x+r\omega)
\]
with the additional condition that $r$ and $t$ have the same sign
(or are both zero).

\subsection{Principal symbol of $N$}
Since we showed that the Schwartz kernel is smooth away from the diagonal,
and we want to prove eventually that $N$ is a $\Psi$DO,
it is enough to study the restriction of $N$ on a small enough set $U$ as above.
We analyze $N_{22}$, for example.
Perform the change of variables $(t,v) \mapsto (r,\omega)$ in (\ref{SN}),
to get
\begin{align}   \nonumber
(N_{22} h)^{i'j'}(x) = \int_{S_xU}\int_{\mathbb{R}} & v^{i'}(r,\omega;x) v^{j'}(r,\omega;x) h_{ij}(x+r\omega) w^i(r,\omega;x) w^j(r,\omega;x) \\
&\times  J^{-1}(r,\omega;x)\,dr\,d\sigma_x(\omega), \label{S19}
\end{align}
where $w(r,\omega,\;x) = \dot\gamma_{x,v}(t)$, and $\text{supp}\, h\subset U$.

This type of integral operators is studied in Appendix \ref{singular}.
Applying Lemma~\ref{lemma_SA} and the remark after it to $N_{22}$,
see (\ref{S19}),
and proceeding similarly for the other operators $N_{kl}$, we get

\begin{Proposition}  \label{pr_Spr1}
$N_{kl}$, $k,l=1,2$, are $\Psi$DOs in $M_1^\text{int}$ with principal symbols
\begin{equation} \label{S22}
\begin{split}
\sigma_p(N_{11})^{i'i}(x,\xi) &= 2\pi\int_{S_xU} \omega^{i'} \omega^i\delta(\omega\cdot \xi)
\, d\sigma_x(\omega),\\
\sigma_p(N_{12})^{i'ij}(x,\xi) &=  2\pi\int_{S_xU} \omega^{i'} \omega^i  \omega^j  \delta(\omega\cdot \xi)
\, d\sigma_x(\omega)=0,\\
\sigma_p(N_{21})^{i'j'i}(x,\xi) &= 2\pi\int_{S_xU} \omega^{i'} \omega^{j'} \omega^i\delta(\omega\cdot \xi)
\, d\sigma_x(\omega)=0,\\
\sigma_p(N_{22})^{i'j'ij}(x,\xi) &= 2\pi\int_{S_xU} \omega^{i'} \omega^{j'} \omega^i \omega^j \delta(\omega\cdot \xi)
\, d\sigma_x(\omega).
\end{split}
\end{equation}
\end{Proposition}

Here $\omega\cdot \xi= \omega^j\xi_j$.
The reason for $\sigma_p(N_{12})=0$, $\sigma_p(N_{21})=0$
is that they are integrals of odd functions.
We therefore get that, as a $\Psi$DO of order $-1$,
$\sigma_p(N) = \text{diag}\big(\sigma_p( N_{22}),
\frac{-2}{n-1} \sigma_p(N_{11}) \big)$.
Note that the principal symbols of $N_{11}$,
$N_{22}$ are the same as in the case of ``ordinary'' geodesics,
i.e., they do not depend on $Y$.
Explicit formulas for $\sigma_p(N_{11})$, $\sigma_p(N_{22})$
can be found in \cite{SU, SSU} and they are based on the analysis of
the Euclidean case in \cite{Sh1}.
Similarly, the principal symbols of $\mathbf{d}$ and $\delt$ are
independent of $Y$ and are given by
\begin{equation}   \label{S23}
\sigma_p(\mathbf{d})[v,\phi]= \text{diag}\left(\frac12( \xi_iv_j+\xi_jv_i),  \xi_i \phi  \right),\quad
\sigma_p(\delt)[h,\beta] = \left( \xi^j h_{ij}, \xi^j\beta_j  \right).
\end{equation}
It is well known \cite{SU,SU1}, and follows immediately from
Proposition~\ref{pr_Spr1}, that $N_{22}$ is elliptic on solenoidal tensors,
i.e., $\sigma_p(N_{22})[h,\beta]=0$ and $\sigma_p(\delt)[h,\beta] =0$
imply $[h,\beta]=0$. Similarly, $N_{11}$ is elliptic on solenoidal
(divergence free) 1-forms. As a result, $N$ is elliptic on solenoidal pairs.

We want to emphasize here that pairs solenoidal in $M$ and extended as zero to $M_1\setminus M$ (that we always assume), may fail to be solenoidal in $M_1$ due to possible jumps at $\partial M$.

\subsection{Parametrix of $N$} We proceed as in \cite{SU, SU1}.
We will define the Hilbert space $\tilde{H}^2(M_1)$ as in \cite{SU,SU1}.
Let $x=(x',x^n)$ be local coordinates in a neighborhood $U$ of a point
on $\partial M$ such that $x^n=0$ defines $\partial M$. Then we  set
\[
\|f\|^2_{\tilde{H}^1(U)} = \int_U \Big(\sum_{j=1}^{n-1} |\partial_{x^j}f|^2 +|x^n\partial_{x^n}f|^2+|f|^2\Big)\, d\vol(x).
\]
This can be extended to a small enough neighborhood $V$ of $\partial M$
contained in $M_1$.
Then we set
\begin{equation}  \label{S24}
\|f\|_{\tilde{H}^2(M_1)} = \sum_{j=1}^{n}  \|\partial_{x^j}f\|_{\tilde{H}^1(V)} + \|f\|_{\tilde{H}^1(M_1)}.
\end{equation}
We also define the $\tilde{H}^2(M_1)$ space of symmetric 2-tensors
and 1-forms, and also the $\tilde{H}^2(M_1)$ space of pairs
$\mathbf{f} = [h,\beta]$ which we denote by $\tilde{\mathbf{H}}^2(M_1)$.
Clearly, the latter is a Hilbert space and $\mathbf{H}^2(M_1) \subset \tilde{\mathbf{H}}^2(M_1) \subset  \mathbf{H}^1(M_1)$.

The space $\tilde{\mathbf{H}}^2(M_1)$ has the property that
for each $\mathbf{f}\in \mathbf{H}^1(M)$ (extended as zero outside $M$),
we have $N\mathbf{f} \in \tilde{\mathbf{H}}^2(M_1)$.
This is not true if we replace $\tilde{\mathbf{H}}^2(M_1)$
by ${\mathbf{H}}^2(M_1)$.

\begin{Lemma}  \label{pr_K}
For any $t=1,2, \dots,\infty$, there exists $k>0$ and a bounded linear operator
\begin{equation}          \label{S4.0}
Q :\tilde{\mathbf{H}}^2(M_1) \to \mathcal{S} \mathbf{L}^2(M)
\end{equation}
such that
\begin{equation}          \label{S4.2}
Q N \mathbf{f} = \mathbf{f}^s + K\mathbf{f}  \quad \text{for all}\quad\mathbf{f}\in \mathbf{H}^1(M),
\end{equation}
where $K: \mathbf{H}^1(M) \to \mathcal{S}\mathbf{H}^{1+t}(M)$ extends to
$K: \mathbf{L}^2(M) \to \mathcal{S}\mathbf{H}^t(M)$. If $t=\infty$, then $k=\infty$. Moreover, $Q$ can be constructed so that $K$ depends continuously on $g$ in a small neighborhood of a fixed $g_0 \in C^k(M)$.
\end{Lemma}

\begin{proof} We will first reconstruct the solenoidal projection
$\mathbf{f}_{M_1}^s$ from $N\mathbf{f}$ modulo smooth terms.
Recall that we extend $\mathbf{f}$ as zero outside $M$,
and  $\mathbf{f}_{M_1}^s$ is the solenoidal projection of the
so-extended $\mathbf{f}$ in $M_1$.
Next, we will reconstruct $\mathbf{f}^s$ from $\mathbf{f}_{M_1}^s$.
We follow \cite{SU,SU1}.

As in  \cite{SU,SU1}, we will  work with $\Psi$DOs with symbols
of finite smoothness $k\gg1$.
All operations we are going to perform would require finitely many
derivatives of the amplitude and finitely many seminorm estimates.
In turn, this would be achieved if $g\in C^k$, $Y\in C^k$, $k\gg1$,
and the corresponding $\Psi$DOs will depends continuously on $g$, $Y$.

Since $N$ is elliptic on solenoidal pairs, we get that
\[
W := N +N_0\mathcal{P}_{M_1} : \mathbf C^\infty(M) \to \mathbf C^l(M_1)
\]
is an elliptic $\Psi$DO of order $-1$ in $M^\text{int}$,
where $N_0$ is any properly
supported parametrix of $(-\Delta)^{1/2}$
(having principal symbol $|\xi|^{-1}_g$).
Next, $W\mathbf{f}$, restricted to a small neighborhood
of $\partial M_1$ in $M_1$, is smooth because
$\text{supp}\, \mathbf{f}\subset M$.
Therefore, there exists a left parametrix $P$ to $W$ such that
$P W-\text{\rm Id} : \mathbf{L}^2(M) \to \mathbf{H}^t(M_1)$
with $t\gg1$, if $k\gg1$.
Then $\mathcal{S}_{M_1} (PW-\text{\rm Id})\mathcal{S}_{M_1}$
has the same property, therefore
$P_1 := \mathcal{S}_{M_1} P$ satisfies
$P_1 N = \mathcal{S}_{M_1} +K_2$,
where $K_2$ has the smoothing properties above, therefore,
\begin{equation}   \label{S25}
P_1N \mathbf{f} = \mathbf{f}_{M_1}^s +K_2\mathbf{f}.
\end{equation}

The next step is to compare $\mathbf{f}^s$ and $\mathbf{f}_{M_1}^s$.
We have $\mathbf{f}^s = \mathbf{f}_{M_1}^s +\mathbf{d} \mathbf{u}$,
where $\mathbf{u} = \mathbf{w}_{M_1}- \mathbf{w}$,
and the latter are the potentials related to $\mathbf{f}$ in $M_1$
and $M$, respectively.
Notice that $\mathbf d$ commutes with the extension as zero
when applied to $\mathbf{w}$ because the latter vanishes on $\partial M$.
Then $\mathbf{u}$ solves the boundary value problem
\begin{equation}  \label{S25a}
(\delt\mathbf{d}) \mathbf{u} =0 \quad \text{in $M$}, \quad \mathbf{u}|_{\partial M} =  \mathbf{w}_{M_1}|_{\partial M}.
\end{equation}
We need to express $\mathbf{w}_{M_1}|_{\partial M}$
in terms of $N\mathbf{f}$. Since $\mathbf{f}=0$ outside $M$,
relation \eqref{S25} implies that
\begin{equation}  \label{S26}
-\mathbf{d} \mathbf{w}_{M_1}= P_1N\mathbf{f}-K_2\mathbf{f} \quad \text{in $M_1\setminus M$.}
\end{equation}
Let $(x,\xi)\in S(M_1^\text{int}\setminus M^\text{int})$ be such that
the magnetic geodesic $\gamma_{x,\xi}(t)$, $t\ge0$,
can be extended in $M_1^\text{int}\setminus M^\text{int}$
for $t\in [0,\ell_1(x,\xi))$ and
$\gamma_{x,\xi}(\ell_1(x,\xi)) \in \partial M_1$;
and moreover, $\gamma_{x,\xi}$ is transversal to $\partial M_1$.
Such  $(x,\xi)$ clearly exist if $M_1$ is close enough to $M$
and if $\gamma_{x,\xi}$ is close to the outgoing
magnetic geodesic normal to $\partial M$.
Using \eqref{S2a}, integrate \eqref{S26} along such $\gamma_{x,\xi}$ to get
\begin{equation}  \label{S27}
\langle \mathbf{w}_{M_1}(x) ,\xi \rangle
=  \int_0^{\ell_1(x,\xi)} \langle ( P_1N\mathbf{f}-K_2\mathbf{f} )(\gamma_{x,\xi}) ,\dot\gamma_{x,\xi} \rangle \, dt.
\end{equation}
Here we denote $\langle \mathbf{w},\xi\rangle = v_j\xi^j+\varphi$,
where $\mathbf{w} = [v(x),\varphi(x)]$ with a $1$-form $v$
and a function $\varphi$, see also \eqref{S2a}.
Similarly we define $\langle \mathbf{f}, \xi\rangle$.
By Lemma~\ref{lemma_S_det}, for a fixed $x$, one can reconstruct $v(x)$, $\varphi$ from $\langle \mathbf{w},\xi\rangle$
known for finitely many  $\xi$'s in any neighborhood of a fixed $\xi$,
and this is done by inverting a matrix.
Moreover, one can do this near any fixed $x$,
and the norm of the solution operator is uniformly bounded in $x$,
see the remark after Lemma~\ref{lemma_S_det}.
By a compactness argument, one can construct an  operator $P_2$ such that
\begin{equation}  \label{S28}
 \left.\mathbf{w}_{M_1}(x)\right|_{\partial M}   =  P_2(P_1N-K_2)\mathbf{f}.
\end{equation}
Arguing as in \cite{SU}, we see that
\begin{equation}  \label{S29}
P_2P_1 : \tilde {\mathbf{H}}^2(M_1) \to \mathbf{H}^{1/2}(\partial M)
\end{equation}
is bounded.
Let $R :\mathbf{H}^{t-1/2}(M) \to \mathbf{H}^t(M)$ be the solution
operator $\mathbf{u} = R\mathbf{h}$ of the boundary value problem
$\delt\mathbf{d} \mathbf{u} =0$ in $M$,
$\mathbf{u}=\mathbf{h}$ on $\partial M$.
Then \eqref{S25}, \eqref{S25a}, \eqref{S28}, and \eqref{S29}
imply (see also \cite{SU})
\[
\mathbf{f}^s = (\text{\rm Id}  +dRP_2)P_1N\mathbf{f} + K\mathbf{f},
\]
where $K$ has the required smoothing properties.
We apply $\mathcal{S}$ to the identity above and set
$Q := \mathcal{S}(\text{\rm Id} +dRP_2)P_1$.

To prove the last statement of the lemma,
we note that all $\Psi$DOs we work with depend continuously on $g$, $\alpha$
if $k\gg1$. The same applies to $\mathcal{S}$, $R$ and $P_2$.
\end{proof}

\subsection{Main results for $C^k$ coefficients}
\begin{Theorem}   \label{thm_estimate}
Let $(g,\alpha)$ be a simple $C^k$ magnetic system on $M$ extended to a simple magnetic system on $M_1$. Then, for $k\gg1$,

$(a)$ $\text{\rm Ker}\,I \cap \mathcal{S}\mathbf{L^2}(M)$
is finitely dimensional and included in $\mathbf C^l(M)$,
where $l\to\infty$ as $k\to\infty$.

$(b)$ Assume that $I$ is s-injective for the pair  $(g,\alpha)$. Then
\begin{equation}   \label{S30}
\|\mathbf{f}^s\|_{\mathbf{L}^2(M)} \le C \|N\mathbf{f}\|_{\tilde{\mathbf{H}}^2(M_1) }
\end{equation}
with a constant $C>0$ that can be chosen to be uniform in $(g,\alpha)$ under a small enough $C^k$ perturbation.
\end{Theorem}

\begin{proof}
This theorem is an analog of \cite[Theorem~2]{SU} except
for the uniformity statement which is similar to an analogous result
in \cite{SU1}. Part (a) follows directly from Lemma~\ref{pr_K}.
Part (b), without the last statement, can be deduced from
Lemma~\ref{pr_K} as well, as in \cite{SU}.
Finally, the proof of the statement about the uniformity of $C$ is
identical to that of  \cite[Theorem~2]{SU1}.
\end{proof}

Without the assumption that $I$ is s-injective,
one has a  hypoelliptic a priori estimate that can
be obtained from \eqref{S30} by adding
$C_s\|f\|_{\mathbf{H}^{-s}(M_1)}$, $\forall s>0$,
to its right-hand side, see \cite[Theorem~2(a)]{SU}.

\begin{Corollary}  \label{cor_S_3}
The set of simple $C^k$ magnetic systems $(g,\alpha)$ with s-injective
magnetic ray transform $I$ is open in the $C^k$ topology, if $k\gg1$.
\end{Corollary}

The next lemma is a linear version of Theorem~\ref{jet}, see also \cite[Lemma~4]{SU1}.

\begin{Lemma}   \label{lemma_S3}
Let $(g,\alpha)$ be a simple $C^k$ magnetic system on $M$. Let $\mathbf{f}\in \mathbf{L}^2(M)$ be such that $I\mathbf{f}=0$.
Then there exists a $C^l$ pair $\mathbf{w}$, with $l\to\infty$ as $k\to\infty$, vanishing on $\partial M$, such that for $\tilde{\mathbf{f}} := \mathbf{f}-\mathbf{dw}$ we have
\begin{equation}  \label{S31}
\partial^m \tilde{\mathbf{f}}|_{\partial M} =0, \quad |m|\le l,
\end{equation}
and if\/ $\tilde{\mathbf{f}} = [\tilde{h},\tilde{\beta}]$, then in
semigeodesic  boundary normal coordinates,
\begin{equation}  \label{S32}
\tilde{h}_{in} =\tilde{\beta}_n=0 \quad \text{for all}\quad i.
\end{equation}
If $k=\infty$, then $l=\infty$ and, in particular, \eqref{S31} holds for all  $m$.
\end{Lemma}

\begin{proof}
In view of Theorem \ref{thm_estimate}(a), we may assume
without loss of generality that $\mathbf f\in \mathbf C^l(M)$.

We now show how to construct $\tilde{\mathbf{f}}$ satisfying \eqref{S32}.
We refer to \cite{E, Sh3} for similar arguments.
Let $(x',x^n)$ be semigeodesic boundary normal coordinates.
Then $g_{in}=\delta_{in}$, $\Gamma_{nn}^i= \Gamma_{in}^n=0$ for all $i$.
Let $\mathbf{f} = [h,\beta]$, $\mathbf{w} = [v,\varphi]$.
Condition \eqref{S32} is equivalent to
\begin{equation}  \label{S33}
\nabla_i v_{n}+\nabla_n v_i = 2h_{in}, \quad \partial_n \varphi -Y^\iota_n v_\iota = \beta_n.
\end{equation}
The first system can be solved by setting first $i=n$ and solving $\partial_n v_n=h_{nn}$, $v_n=0$ for $x^n=0$ by integration. Then we solve the remaining system of $n-1$ ODEs of the form
\begin{equation}  \label{S34}
\partial_n v_{\iota} - 2\Gamma^\kappa_{\iota n}v_{\kappa} = 2 h_{\iota n} - \partial_\iota v_n, \quad \iota =1,\dots, n-1,
\end{equation}
with initial conditions $v_\iota=0$ for $x^n=0$.
Finally, we solve the second equation in \eqref{S33} with
the same zero initial condition.
This defines $\mathbf{w}$ near $\partial M$.
To define $\mathbf{w}$ in the whole $M$,
we replace $\mathbf{w}$ by $\chi\mathbf{w}$,
where $\chi$ is an appropriate smooth cut-off function equal to $1$
near $\partial M$ and supported in a larger neighborhood of $\partial M$.

To prove \eqref{S31}, we argue as in Theorem~\ref{jet}, following \cite{LSU}.
Assume that \eqref{S31} is not true at some $x_0\in \partial M$.
Then by studying the Taylor expansion of
$\langle\tilde{\mathbf{f}},\xi\rangle
=\tilde{h}_{\iota\kappa}\xi^\iota \xi^\kappa+ \tilde{\beta}_\iota \xi^\iota$
near $x=x_0$, we see that there is
$\xi_0\in T(\partial M)$ such that
$\langle\tilde{\mathbf{f}},\xi\rangle$
is either strictly negative or strictly positive for $\xi$ lying in some neighborhood
of $\xi_0$.
This contradicts the fact that $I\tilde{\mathbf{f}}=0$
if we integrate over magnetic geodesics originating from $x_0$
with directions close to $\xi_0$.
\end{proof}

\subsection{Analytic magnetic systems}
Assume that $M$ is an analytic manifold
with smooth boundary $\partial M$ that does not need to be analytic
(that can always be achieved by choosing an analytic atlas;
and in case we have a simple system,
we can start with a fixed global coordinate system and
do only analytic changes of variables).
We will show that then $I$ is s-injective,
if $(g,\alpha)$ are analytic.
By ``analytic'', we mean real analytic,
and we say that $f$ is analytic in the set $X$,
not necessarily open, if $f$ is analytic in a neighborhood of $X$.
Then we write $f\in \mathcal{A}(X)$.
Notice that one can construct an analytic $M_1$ as before.

The central result of this section is the following.
\begin{Theorem}   \label{thm_S_an}
If $(g,\alpha)$ is an analytic magnetic system on $M$, then
$I$ is \mbox{s-injective}.
\end{Theorem}

We will give a proof at the end of this section.

\begin{Theorem}
Let $g$, $\alpha$ be analytic in $M_1$.
Then $N_{kl}$, $k,l=1,2$, are analytic $\Psi$DOs in
$M_1^\text{int}$ with principal symbols as in Proposition~\ref{pr_Spr1}.
\end{Theorem}

\begin{proof}
Notice first that by \eqref{SN} and the simplicity assumption,
$N_{kl}$ have Schwartz kernels that are analytic away from the diagonal.
Therefore, it is enough to prove the theorem for
$N_{1,2}$ restricted to an arbitrary small open subset of $M_1$.
This, however, follows from \eqref{S19} and Lemma~\ref{lemma_SA_an}
of Appendix \ref{singular}.
\end{proof}

\begin{Lemma}  \label{pr_S_an2}
Let $g$, $\alpha$ be analytic in $M_1$ and assume that
$I\mathbf{f}=0$, $\mathbf{f}\in\mathbf{L}^2(M)$.
Then $\mathbf{f}^s\in \mathcal{A}(M)$.
\end{Lemma}

\begin{proof}
Consider the solenoidal projection
$\mathbf{f}^s_{M_1}= \mathbf{f} -\mathbf{d}\mathbf{w}_{M_1}$
of $\mathbf{f}$ (extended as $0$ outside $M$) on $M_1$.
Since $\delt \mathbf{f}^s_{M_1}=0$
and $N\mathbf{f}^s_{M_1}=0$ in $M_1^\text{int}$, and since $\delt$ and
$N$ together form an elliptic system of analytic $\Psi$DOs
(we can apply an elliptic  $\Psi$DO of order $2$ to the left of $N$
to make the new operator and $\delt$ of the same order,
see also \cite{SU,SU1}), we get that
$\mathbf{f}^s_{M_1} \in  \mathcal{A}( M_1^\text{int})$.
On the other hand, $\mathbf{w}_{M_1}$
solves $\delt\mathbf{d} \mathbf{w}_{M_1} =0$
in $M_1^\text{int}\setminus M$, $\mathbf{w}_{M_1} =0$ on $\partial M_1$,
and by elliptic boundary regularity, see \cite[Lemma~3]{SU1}
and references therein, we have that $\mathbf{w}_{M_1}$
is analytic up to $\partial M_1$, therefore
$\mathbf{f}^s_{M_1} \in  \mathcal{A}( M_1)$.

Next, we have $\mathbf{f}^s = \mathbf{f}^s_{M_1} +\mathbf{du}$,
where $\mathbf{u} =\mathbf{w}_{M_1}-\mathbf{w}$,
see also the proof of Lemma~\ref{pr_K}.
Then $\mathbf{u}$ solves the boundary value problem \eqref{S25a},
and all we need is to prove
that $\mathbf{w}_{M_1}|_{\partial M}$ is analytic.
Note that in general, $\mathbf{w}_{M_1}$ is not analytic
across $\partial M$, because $\mathbf{dw}_{M_1}$ may have a jump there
but it belongs to $H^1$ because $\mathbf{d}$ is elliptic;
therefore the trace on $\partial M$ is well defined.
Since $\mathbf{d}\mathbf{w}_{M_1}=\mathbf{f}^s_{M_1}$
in $M_1^\text{int}\setminus M$, and it is analytic in its closure,
in the notation of \eqref{S27}, we have
\[
\langle \mathbf{w}_{M_1}(x),\xi\rangle
= \int_0^{\ell_1(x,\xi)} \langle\mathbf{f}^s_{M_1} (\gamma_{x,\xi}), \dot \gamma_{x,\xi} \rangle\, dt
\]
for $x\in\partial M$ and $\xi$ in a small neighborhood of
the unit exterior normal to $\partial M$.
This, combined with Lemma~\ref{lemma_S_det} and the remark after it,
shows that $\mathbf{w}_{M_1}|_{\partial M}$ is analytic.
\end{proof}

\begin{Lemma} \label{lemma_S4}
Under the assumptions of Lemma~\ref{lemma_S3},
assume that $g$ and $\alpha$ are analytic in a neighborhood of $\partial M$.
Then $\tilde{\mathbf{f}}=0$ in a neighborhood of $\partial M$.
\end{Lemma}

\begin{proof}
It is enough to notice that near $\partial M$, $\mathbf{w}$
in the proof of Lemma~\ref{lemma_S3} is obtained by solving ODEs
with analytic coefficients.
Therefore, $\tilde{\mathbf{f}}$ is analytic in a  neighborhood of $\partial M$,
and it must vanish there by \eqref{S31}.
\end{proof}

\begin{proof}[Proof of Theorem~\ref{thm_S_an}]
We first find $\mathbf{w}$ as in Lemma~\ref{lemma_S3}
and Lemma~\ref{lemma_S4} such that
$\tilde{\mathbf{f}} = \mathbf{f}^s-\mathbf{dw}$ vanishes near $\partial M$,
i.e.,
\begin{equation}  \label{S45a}
\mathbf{f}^s = \mathbf{dw} \quad \text{near $\partial M$.}
\end{equation}
Our next goal is to show that one can extend $\mathbf{w} = [v,\varphi]$
to the whole $M$ analytically; then \eqref{S45a} would be preserved
by analytic continuation (recall that $\mathbf{f}^s$ is analytic by
Lemma~\ref{pr_S_an2}), and this would imply $\mathbf{f}^s=0$.

Let $u_\pm(x,\xi)$ be the solution of the kinetic equation
$\mathbf{G}_\mu u= \langle \mathbf{f}^s,\xi\rangle $
with $u=0$ on $\partial_\pm SM$. Integrate to get
\[
u_\pm(x,\xi) = \mp\int_0^{\ell^\pm(x,\xi)}\langle \mathbf{f}^s(\gamma_{x,\xi}),\dot\gamma_{x,\xi}\rangle \, dt,
\]
where $\ell^+$ stands for $\ell$.
Since $I\mathbf{f}=0$, we have $u_++u_-=0$.
For $(x,\xi)$ such that $x$ is close enough to $\partial M$
and $\xi$ close enough to $-e_n$ (in boundary normal coordinates),
we have $u_+(x,\xi) = \langle \mathbf{dw},\xi\rangle = \langle v(x),\xi\rangle +\varphi(x)$
by \eqref{S45a} and \eqref{S2a}.
This also implies $-u_-(x,\xi) =  \langle v(x),\xi\rangle +\varphi(x)$.
If we replace $\xi$ by $-\xi$, since $-\xi$ is close to $e_n$,
we get by  \eqref{S45a},  $-u_-(x,-\xi) =  -\langle v(x),\xi\rangle +\varphi(x)$;
therefore, $u_+(x,-\xi) =  -\langle v(x),\xi\rangle +\varphi(x)$.

In particular, for $x$ close enough to $\partial M$ and $\xi\in S_xM$
close enough to either $e_n$ or $-e_n$,
\begin{equation}   \label{S45b}
\varphi(x) = u_\text{+,\,even}(x,\xi),
\quad v_j(x) = \partial_{\eta^j}|\eta|_gu_\text{+,\,odd}(x,\eta/|\eta|_g)\big|_{|\eta|_g=1},
\end{equation}
where $u_\text{+,\,odd/even}$ stands for the odd/even part of $u_+$ w.r.t.\ $\xi$.
The derivative on the right-hand side of \eqref{S45b} can be written as
$Pu_\text{+\,,odd}$, where $P$ is a first order differential operator on
$S_xM$ with coefficients analytically depending on $x$ and $\xi$.
A direct differentiation shows that
$(Pu)_j = \partial_{\xi^j}u -\xi_j\xi^i\partial_{\xi^i}u+\xi_ju $,
and the first order part is clearly tangent to $S_xM$.

We define $\varphi(x,\xi)$ and $v(x,\xi)$ by \eqref{S45b} on the whole
$SM$. Since  $u_+(x,\xi)$ is analytic on $SM$, so are
$\varphi(x,\xi)$ and $v(x,\xi)$.
By \eqref{S45b}, in some open set they are independent of $\xi$,
i.e., $\text{grad}_\xi \phi(x,\xi)=0$, $\text{grad}_\xi v(x,\xi)=0$
there, where $\text{grad}_\xi$ stands for the gradient on $S_xM$.
Those equalities extend to the whole $SM$ by analytic continuation;
therefore, \eqref{S45b} defines $\xi$-independent $\varphi$ and $v$
in the whole $M$, and they are analytic functions of $x$.
\end{proof}

We are ready now to state the generic s-injectivity result.

\begin{Definition}\label{generic-set}
For a fixed manifold $M$, we define
$\mathcal{G}^k$ to be the set  of simple $C^k$ pairs $(g,\alpha)$
with s-injective magnetic ray transform $I=I_{g,\alpha}$.
\end{Definition}

\begin{Theorem} \label{S_thm_lin}
There exists $k_0>0$, such that for $k\ge k_0$,
the set $\mathcal{G}^k$ is open and dense in the set of all simple
$C^k$ pairs $(g,\alpha)$ and contains all real analytic simple pairs.
\end{Theorem}

\begin{proof}
By Corollary~\ref{cor_S_3}, $\mathcal{G}^k$ is open. By Theorem~\ref{thm_S_an},
it is dense.
\end{proof}

\section{Energy estimates method}  \label{sec_5}

It is easy to see that if $(M,g,\alpha)$ is a simple smooth magnetic system
and the magnetic ray transform of a smooth function $\phi:SM\to\mathbb R$
vanishes, then $\phi$ is the flow derivative of
a unique function $u$ which is continuous on $SM$, smooth on
$SM\setminus S(\partial M)$, and
vanishes on $\partial(SM)$:
\begin{equation}\label{kinetic}
\mathbf G_\mu u=\phi,
\quad u|_{\partial(SM)}=0.
\end{equation}
Indeed, $u$ is defined by the formula
\begin{equation}\label{def-u}
u(x,\xi)=-\int_0^{\ell(x,\xi)}\phi(\psi^t(x,\xi))\,dt,
\quad (x,\xi)\in SM.
\end{equation}

In this section we will analyze the linear problem \eqref{kinetic} for
$\phi(x,\xi)$ of degree at most $1$ in $\xi$
(the linear problem for $1$-tensors)
and  for $\phi(x,\xi)$ of degree at most $2$ in $\xi$
(the linear problem for $2$-tensors).

\subsection{Semibasic tensor fields}
We recall the notion of semibasic tensor field
and its derivatives (we prefer to adhere
to the notations of \cite{DP}).

Let $\pi:TM\setminus\{0\}\to M$ be the natural projection,
and let $\beta^r_s M:=\pi^*\tau^r_s M$ denote the bundle of semibasic tensors
of degree $(r,s)$, where $\tau^r_s M$ is the bundle of tensors
of degree $(r,s)$ over $M$.
Sections of the bundles $\beta^r_sM$ are called semibasic tensor fields,
and we denote the space of smooth sections by $C^\infty(\beta^r_sM)$
(in particular, $C^\infty(\beta^0_0M)=C^\infty(TM\setminus\{0\})$).
For such a field $T$, the coordinate representation
$$
T=(T^{i_1\dots i_r}_{j_1\dots j_s})(x,\xi)
$$
holds in the domain of a standard local coordinate system $(x^i,\xi^i)$
on $TM\setminus\{0\}$
associated with a local coordinate system $(x^i)$ in $M$.
Under a change of a local coordinate system, the components of
a semibasic tensor field are transformed by  the same formula
as those of an ordinary tensor field on $M$.

Every ``ordinary'' tensor field on $M$ defines
a semibasic tensor field by the rule $T\mapsto T\circ\pi$, so that
the space of tensor fields on $M$ can be treated as embedded in the space
of semibasic tensor fields.

For a semibasic tensor field $(T^{i_1\dots i_r}_{j_1\dots j_s})(x,\xi)$,
the horizontal derivative is defined by
\begin{multline*}
T^{i_1\dots i_r}_{j_1\dots j_s|k}
=\de{}{x^k}T^{i_1\dots i_r}_{j_1\dots j_s}
-\Gamma^p_{kq}\xi^q\de{}{\xi^p}T^{i_1\dots i_r}_{j_1\dots j_s}
\\
+\sum_{m=1}^r\Gamma^{i_m}_{kp}
T^{i_1\dots i_{m-1}pi_{m+1}\dots i_r}_{j_1\dots j_s}
-\sum_{m=1}^r\Gamma^{p}_{kj_m}
T^{i_1\dots i_r}_{j_1\dots j_{m-1}pj_{m+1}\dots j_s},
\end{multline*}
the vertical derivative by
$$
T^{i_1\dots i_r}_{j_1\dots j_s\cdot k}
=\de{}{\xi^k}T^{i_1\dots i_r}_{j_1\dots j_s},
$$
and the modified horizontal derivative by
$$
T^{i_1\dots i_r}_{j_1\dots j_s: k}
=T^{i_1\dots i_r}_{j_1\dots j_s|k}+|\xi|Y^j_k T^{i_1\dots i_r}_{j_1\dots j_s\cdot j}.
$$

The operators
$$
\nabla_{|}:C^\infty(\beta^r_sM)\to C^\infty(\beta^r_{s+1}M),
\quad
\nabla_{\cdot}:C^\infty(\beta^r_sM)\to C^\infty(\beta^r_{s+1}M),
$$
and
$$
\nabla_{:}:C^\infty(\beta^r_sM)\to C^\infty(\beta^r_{s+1}M),
$$
are defined as
$$
(\nabla_{|} T)^{i_1\dots i_r}_{j_1\dots j_s k}
=\nabla_{|k} T^{i_1\dots i_r}_{j_1\dots j_{s}}
:=T^{i_1\dots i_r}_{j_1\dots j_s|k},
$$
$$
(\nabla_{\cdot} T)^{i_1\dots i_r}_{j_1\dots j_{s}k}
=\nabla_{\cdot k} T^{i_1\dots i_r}_{j_1\dots j_{s}}
=T^{i_1\dots i_r}_{j_1\dots j_s\cdot k}
$$
and
$$
(\nabla_{:} T)^{i_1\dots i_r}_{j_1\dots j_s k}
=\nabla_{:k} T^{i_1\dots i_r}_{j_1\dots j_{s}}
:=T^{i_1\dots i_r}_{j_1\dots j_s:k}.
$$
For convenience, we also define $\nabla^{|}$, $\nabla^{\cdot}$,
and $\nabla^{:}$ as
$$
\nabla^{|i}=g^{ij} \nabla_{|j},\quad
\nabla^{\cdot i}=g^{ij}\nabla_{\cdot j},
\quad \nabla^{:i}=g^{ij} \nabla_{:j}.
$$
In \cite{PSh, Sh1}, the operators $\nabla_{|}$ and $\nabla_{\cdot}$
were denoted by $\overset h \nabla$
and $\overset v \nabla$ respectively.

Given $u\in C^\infty(TM\setminus\{0\})$, we define
$$
\mathbf X u(x,\xi)=\xi^i u_{:i}=\xi^i(u_{|i}+|\xi|Y^j_iu_{\cdot j}).
$$
Note that $\mathbf X$ restricted to $SM$ coincides with $\mathbf G_\mu$.

For $V=(V^i)\in C^\infty(\beta^1_0M)$, we set
$$
\divh V:=V^i_{|i},
\quad
\divm V:=V^i_{:i},
\quad
\divv V:=V^i_{\cdot i},
\quad
(\mathbf X V)^i
=\xi^kV^i_{:k}.
$$
Note if $\gamma$ is a magnetic geodesic, then
$$
(\mathbf X V)(\gamma(t),\dot\gamma(t))
=\nabla_{\dot\gamma}(V(\gamma(t),\dot\gamma(t)).
$$

Given a function $u:SM\to\mathbb R$, we will also denote by $u$
its extension to a positively homogeneous function of degree $0$
on $TM\setminus\{0\}$
(hoping that this will not yield any confusion).
For a smooth $\phi$ the smoothness properties of $u$
defined by \eqref{def-u} are
determined by those of $\ell(x,\xi)$. As mentioned
in the beginning of Subsection \ref{stattering}, the latter function
is smooth on $SM\setminus S(\partial M)$. All points of $S(\partial M)$
are singular for $\ell$ as some derivatives of $\ell$ are unbounded
in a neighborhood of such a point. Nonetheless, some derivatives are bounded.
Let $\rho$ be a smooth function on $M$ such that $\partial M=\rho^{-1}(0)$
and $|\grad\rho|=1$ in some neighborhood of $\partial M$.
Define
\begin{equation*}
\nabla^{:}_\rho u(x,\xi)
=\nabla^:u(x,\xi)-\langle \nabla^:u(x,\xi),\grad\rho(x)\rangle\grad \rho(x).
\end{equation*}
Note that $\nabla^{:}_\rho u(x,\xi)$ is completely determined by the restriction of $u$
to the level set of $\rho$ that contains $(x,\xi)$.

\begin{Lemma}[{cf. \cite[Lemma 3.2.3]{Sh2}}]\label{properties-l}
The semibasic vector fields $\nabla^{:}_\rho\ell$ and $\nabla^{\cdot}\ell$
are bounded on $SM\setminus S(\partial M)$.
\end{Lemma}
\begin{proof}
Clearly, $\nabla^:_\rho\rho=0$.
Let $h(x,\xi,t)=\rho(\gamma_{x,\xi}(t))$.
Since $h(x,\xi,0)=\rho(x)$, it follows that
$\nabla^:_\rho h(x,\xi,0)=0$
and therefore
$$
\nabla^:_\rho h(x,\xi,t)=a(x,\xi,t)t
$$
for some smooth field $a$.
Differentiating the equality $h(x,\xi,\ell(x,\xi))=0$ for
$(x,\xi)\in SM\setminus S(\partial M)$, we obtain
\begin{multline*}
\nabla^:_\rho h(x,\xi,\ell(x,\xi))
+\de{\rho}{t}(x,\xi,\ell(x,\xi))\nabla^:_\rho\ell(x,\xi)
\\
=\nabla^:_\rho h(x,\xi,\ell(x,\xi))
+\langle\grad\rho(y),\eta\rangle \nabla^:_\rho\ell(x,\xi)=0,
\end{multline*}
where $(y,\eta)=\psi^{\ell(x,\xi)}(x,\xi)\in\partial_-SM$.
Hence,
$$
\nabla^:_\rho\ell(x,\xi)
=-\frac{\ell(x,\xi)}{\langle\grad\rho(y),\eta\rangle}\ a(x,\xi,\ell(x,\xi)).
$$
Since $\ell(x,\xi)\le |\mathbb L(y,\eta)|$, boundedness of
$\nabla^:_\rho\ell(x,\xi)$ now follows from Lemma \ref{l-l}.

Boundedness of $\nabla^{\cdot}\ell$ is established similarly.
\end{proof}

\subsection{Identities}
The next identities are particular cases of those in
\cite[Lemma 4.6 and Lemma 4.7]{DP}. For $(x,\xi)\in SM$:
\begin{align}
2\langle\nabla^{:} u,\nabla^{\cdot}(\mathbf X u)\rangle
&=|\nabla^{:} u|^2
+\mathbf X (\langle\nabla^{:} u,\nabla^{\cdot} u\rangle)
-\divm((\mathbf X  u)\nabla^{\cdot} u)
\nonumber
\\
&\quad+\divv((\mathbf X u)\nabla^{:} u)
-\langle\tilde{\mathbf R}_{\xi}(\nabla^{\cdot} u),\nabla^{\cdot} u\rangle
+\langle Y(\nabla^{\cdot}u),\nabla^{:}u\rangle,
\label{pestov-dif}
\\
\label{xnabla}
2\langle\nabla^{:} u,\nabla^{\cdot}(\mathbf X u)\rangle
&=|\nabla^{:}u|^2+|\nabla^{\cdot} (\mathbf X  u)|^2
+\langle Y(\xi),\nabla^{\cdot}u\rangle^2
-|\mathbf X(\nabla^{\cdot}u)|^2,
\end{align}
where
$\tilde {\mathbf R}_\xi(Z)=\left(\tilde R^i_{jkl}\xi^j\xi^l Z^k\right)$,
$$
\tilde R^i_{jkl}=R^i_{jkl}+\xi_j(Y^i_{k|l}-Y^i_{l|k})
+g_{js}(Y^s_l Y^i_k-Y^s_k Y^i_l),
$$
and $(R^i_{jkl})$ is the Riemann curvature tensor.

For a unit speed magnetic geodesic~$\gamma$, the operator $\mathcal C$
on smooth vector fields along $\gamma$ is defined by (see \cite[(45)]{DP})
\begin{equation*}
{\mathcal C}(Z)=R(\dot{\gamma},Z)\dot{\gamma}- Y(\dot{Z})
-(\nabla_{Z}Y)(\dot{\gamma}),
\end{equation*}
where $\dot Z$ stands for the covariant derivative along $\gamma$,
$\dot Z=\nabla_{\dot\gamma}Z$.

The operator $\tilde{\mathcal C}$
on smooth semibasic vector fields on $TM\setminus\{0\}$
is defined by (see \cite[Subsection 4.6]{DP})
$$
\tilde{\mathcal C}(Z)(x,\xi)=\mathbf R_\xi(Z)-Y(\mathbf X Z)-(\nabla_{|Z} Y)(\xi),
$$
where $\mathbf R_\xi(Z)=\left(R^i_{jkl}\xi^j\xi^l Z^k\right)$
is the curvature operator.

Notice that
$$
\tilde{\mathcal C}(Z)(\gamma(t),\dot\gamma(t))
=\mathcal C(Z(\gamma(t),\dot\gamma(t))).
$$
In view of this identity, we henceforth omit the tilde in the notation of
$\tilde{\mathcal C}$, hoping that no confusion will arise.

\begin{Lemma}\label{dif-ids}
The following hold for $(x,\xi)\in SM$:
\begin{multline}\label{dif-1}
|\mathbf X(\nabla^{\cdot}u)|^2
-\langle {\mathcal C}(\nabla^{\cdot}u),\nabla^{\cdot}u\rangle
-\langle Y(\xi),\nabla^{\cdot}u\rangle^2
=
|\nabla^{\cdot} (\mathbf X  u)|^2
-\divv((\mathbf X u)\nabla^{:} u)
\\
-\divh(\langle\nabla^{:} u,\nabla^{\cdot} u\rangle\xi
-(\mathbf X  u)\nabla^{\cdot} u)
+\divv(\langle\nabla^{:} u,\nabla^{\cdot} u\rangle Y(\xi)),
\end{multline}
\begin{multline}\label{dif-2}
|\mathbf X(\nabla^{\cdot}u)|^2
-2\langle {\mathcal C}(\nabla^{\cdot}u),\nabla^{\cdot}u\rangle
-\langle Y(\xi),\nabla^{\cdot}u\rangle^2
\\
=|\nabla^{\cdot} (\mathbf X  u)|^2
-2\divv((\mathbf X u)\nabla^{:} u)
-|\nabla^{:} u|^2
+2\langle\nabla^{:} u,\nabla^{\cdot}(\mathbf X u)\rangle
\\
-2\divh(\langle\nabla^{:} u,\nabla^{\cdot} u\rangle\xi
-(\mathbf X  u)\nabla^{\cdot} u)
+\divv(\langle\nabla^{:} u,\nabla^{\cdot} u\rangle Y(\xi)).
\end{multline}
\end{Lemma}

\begin{proof}
Note that
\begin{multline}\label{x-trans}
\mathbf X(\langle\nabla^{:} u,\nabla^{\cdot} u\rangle)(x,\xi)
=\xi^j(\langle\nabla^{:} u,\nabla^{\cdot} u\rangle)_{:j}
=(\langle\nabla^{:} u,\nabla^{\cdot} u\rangle\xi^j)_{:j}
-\langle\nabla^{:} u,\nabla^{\cdot} u\rangle\xi^j_{:j}
\\
=\divm(\langle\nabla^{:} u,\nabla^{\cdot} u\rangle\xi)
-\langle\nabla^{:} u,\nabla^{\cdot} u\rangle Y^j_j
=\divm(\langle\nabla^{:} u,\nabla^{\cdot} u\rangle\xi),
\end{multline}
because $Y^j_j=0$ by the skew-symmetry of $Y$.

Using \eqref{x-trans}, we first change \eqref{pestov-dif} to
\begin{multline}\label{pestov-r}
2\langle\nabla^{:} u,\nabla^{\cdot}(\mathbf X u)\rangle
=|\nabla^{:} u|^2
+\divm V +\divv((\mathbf X u)\nabla^{:} u)
\\
-\langle \tilde{\mathbf R}_{\xi}(\nabla^{\cdot} u),\nabla^{\cdot} u\rangle
+\langle Y(\nabla^{\cdot}u),\nabla^{:}u\rangle,
\end{multline}
with
$$
V=\langle\nabla^{:} u,\nabla^{\cdot} u\rangle\xi
-(\mathbf X  u)\nabla^{\cdot} u.
$$
Next,
\begin{align*}
\divm V&=V^i_{:i}=V^i_{|i}+Y^j_i V^i_{\cdot j}
=\divh V+\divv(Y(V))
\\
&=\divh V+\divv(\langle\nabla^{:} u,\nabla^{\cdot} u\rangle Y(\xi))
-\divv((\mathbf X  u)Y(\nabla^{\cdot} u))
\end{align*}
and
\begin{align*}
\divv((\mathbf X  u)Y(\nabla^{\cdot} u))
&=\big((\mathbf Xu)Y^j_i g^{ik}u_{\cdot k}\big)_{\cdot j}
=(\mathbf X u)_{\cdot j}Y^j_i g^{ik}u_{\cdot k}
+(\mathbf X u)g^{ki} Y^j_i u_{\cdot k\cdot j}
\\
&=\langle\nabla^{\cdot } (\mathbf X u),Y(\nabla^{\cdot} u)\rangle,
\end{align*}
because $g^{ki} Y^j_i u_{\cdot k\cdot j}=0$ by the skew-symmetry of $Y$
and symmetry of mixed derivatives of $u$.

Then \eqref{pestov-r} takes the form
\begin{multline}\label{pestov-changed}
2\langle\nabla^{:} u,\nabla^{\cdot}(\mathbf X u)\rangle
=|\nabla^{:} u|^2
+\divh V+\divv(\langle\nabla^{:} u,\nabla^{\cdot} u\rangle Y(\xi))
-\langle\nabla^{\cdot } (\mathbf X u),Y(\nabla^{\cdot} u)\rangle
\\
+\divv((\mathbf X u)\nabla^{:} u)
-\langle \tilde{\mathbf R}_{\xi}(\nabla^{\cdot} u),\nabla^{\cdot} u\rangle
+\langle Y(\nabla^{\cdot}u),\nabla^{:}u\rangle.
\end{multline}

We have (see \cite[Subsection 4.6]{DP}):
\begin{multline*}
\langle {\mathcal C}(\nabla^{\cdot}u),\nabla^{\cdot}u\rangle
=\langle {\mathbf R}_{\xi}(\nabla^{\cdot}u),\nabla^{\cdot} u\rangle
+\langle \nabla^{\cdot}(\mathbf X  u),
Y(\nabla^{\cdot}u)\rangle
-\langle\nabla^{:}u,Y(\nabla^{\cdot}u)\rangle
\\
+\langle Y(\xi),\nabla^{\cdot}u\rangle^2
-\langle(\nabla_{|(\nabla^{\cdot} u)}Y)(\xi),\nabla^{\cdot}u\rangle.
\end{multline*}

Also,
\begin{align*}
\langle \tilde{\mathbf R}_{\xi}(\nabla^{\cdot} u),\nabla^{\cdot} u\rangle
&=(R^i_{jkl}+\xi_j(Y^i_{k|l}-Y^i_{l|k})
+g_{js}(Y^s_l Y^i_k-Y^s_k Y^i_l))u^{\cdot k} v^l v^j u_{\cdot i}\nonumber
\\
&=\langle R(\nabla^{\cdot} u,\xi)\xi,\nabla^{\cdot} u\rangle
+\langle \nabla_\xi Y(\nabla^{\cdot} u),\nabla^{\cdot}u\rangle
-\langle\nabla_{(\nabla^{\cdot}u)}Y(\xi),\nabla^{\cdot}u\rangle \nonumber
\\
&\quad+\langle Y(\xi),\xi\rangle\langle Y(\nabla^{\cdot} u),\nabla^{\cdot}u\rangle
-\langle Y(\nabla^{\cdot}u),\xi\rangle\langle Y(\xi),\nabla^{\cdot}u\rangle \nonumber
\\
&=\langle \mathbf R_\xi(\nabla^{\cdot} u),\nabla^{\cdot} u\rangle
-\langle(\nabla_{(\nabla^{\cdot}u)}Y)(\xi),\nabla^{\cdot}u\rangle
+\langle Y(\xi),\nabla^{\cdot}u\rangle^2.      
\end{align*}
Therefore,
\begin{equation}\label{tilder}
\langle \tilde{\mathbf R}_{\xi}(\nabla^{\cdot} u),\nabla^{\cdot} u\rangle
=\langle {\mathcal C}(\nabla^{\cdot}u),\nabla^{\cdot}u\rangle
-\langle \nabla^{\cdot}(\mathbf X  u),
Y(\nabla^{\cdot}u)\rangle
+\langle\nabla^{:}u,Y(\nabla^{\cdot}u)\rangle.
\end{equation}

Using \eqref{tilder} in \eqref{pestov-changed}, we get
\begin{multline}\label{pestov-dif-final}
2\langle\nabla^{:} u,\nabla^{\cdot}(\mathbf X u)\rangle
=|\nabla^{:} u|^2
+\divh V+\divv(\langle\nabla^{:} u,\nabla^{\cdot} u\rangle Y(\xi))
+\divv((\mathbf X u)\nabla^{:} u)
\\
-\langle {\mathcal C}(\nabla^{\cdot}u),\nabla^{\cdot}u\rangle.
\end{multline}

Subtracting \eqref{pestov-dif-final} from \eqref{xnabla},
we arrive at \eqref{dif-1}.
Multiplying \eqref{pestov-dif-final} by $2$ and subtracting
the result from \eqref{xnabla},
we arrive at \eqref{dif-2}.
\end{proof}

\subsection{Linear problem for $1$-tensors}

\begin{Theorem}\label{linear}
Let $(M,g,\alpha)$ be a simple magnetic system,  $v$
a square integrable $1$-form, and $\varphi$ a square integrable
function on $M$.
If the magnetic ray transform of the function
$\phi(x,\xi)=v_i(x)\xi^i+\varphi(x)$\ vanishes, then $\varphi=0$
and $v=dh$ for some function $h\in H^1_0(M)$.
\end{Theorem}

\begin{proof}
Recall that the space $\mathcal L^2(M)$
consists of the pairs
$[v,\varphi]$, where $v$ is a square integrable $1$-form
and $\varphi$ is a square integrable function on $M$, furnished with
the norm \eqref{S7} (see Section \ref{solenoid-potential}).
In $\mathcal L^2(M)$, we consider the subspace
$\mathcal S\mathcal L^2(M)$ of solenoidal pairs
$[v,\varphi]$, defined by the condition $\delta v=0$,
and the subspace $\mathcal P\mathcal L^2(M)$ of potential pairs
$[dh,0]$, where $h\in H^1_0(M)$. Then $\mathcal L^2(M)$
decomposes as the orthogonal direct sum of $\mathcal S\mathcal L^2(M)$
and $\mathcal P\mathcal L^2(M)$.

Associating each pair $[v,\varphi]\in\mathcal L^2(M)$
with the function $\phi(x,\xi)=v_i(x)\xi^i+\varphi(x)$, define
$$
\mathcal I[v,\varphi]=I\phi.
$$

Clearly, $\mathcal I$ vanishes on $\mathcal P\mathcal L^2(M)$ and,
to prove the theorem,
it suffices to show that $\ker \mathcal I\cap \mathcal S\mathcal L^2(M)=0$.

By considerations similar to those in Section 4, one can prove, as
in part (a) of Theorem \ref{thm_estimate}, that
$\text{\rm Ker}\,\mathcal I \cap \mathcal{S}\mathcal L^2(M)
\subset \mathcal C^l(M)$.

Thus, we may assume under assumptions of the theorem that
$[v,\varphi]\in \mathcal C^l(M)$.

Define $u:SM\to \mathbb R$ by means of \eqref{def-u}. Then $u$ satisfies
the boundary value problem \eqref{kinetic}.  As before, we preserve
the notation $u$ for the extension of $u$ to a positively homogeneous function
of degree $0$ on $TM\setminus\{0\}$.
Then for $(x,\xi)\in TM\setminus\{0\}$ we have
$$
\mathbf X u(x,\xi)=v_i(x)\xi^i+|\xi|\varphi(x).
$$

Now, we wish to integrate identity \eqref{dif-1} over $SM$.
However, $u$ has singularities on $T(\partial M)$, and we need
some precautions against them. We proceed in the same way as in the proof
of Theorem 3.4.3 in \cite{Sh2}.
Let $\rho:M\to\mathbb R$ be a~nonnegative smooth function such that
$\partial M=\rho^{-1}(0)$ and $|\grad\rho|=1$ near $\partial M$.
For $\varepsilon>0$, let $M_\varepsilon=\{x\in M:\rho(x)\ge \varepsilon\}$.
Then $u$ is smooth on $SM_\varepsilon$.
Integrating \eqref{dif-1} over $SM_\varepsilon$,
we get
\begin{multline}\label{integral-1}
\int_{SM_\varepsilon}\left\{|\mathbf X(\nabla^{\cdot}u)|^2
-\langle {\mathcal C}(\nabla^{\cdot}u),\nabla^{\cdot}u\rangle
-\langle Y(\xi),\nabla^{\cdot}u\rangle^2\right\}\,d\Sigma^{2n-1}(x,\xi)
\\
=
\int_{SM_\varepsilon}|\nabla^{\cdot} (\mathbf X  u)|^2\,d\Sigma^{2n-1}
-\int_{SM_\varepsilon}\divv((\mathbf X u)\nabla^{:} u)\,d\Sigma^{2n-1}
\\
-\int_{SM_\varepsilon}\divh(\langle\nabla^{:} u,\nabla^{\cdot} u\rangle\xi
-(\mathbf X  u)\nabla^{\cdot} u) \,d\Sigma^{2n-1}
\\
+\int_{SM_\varepsilon}\divv(\langle\nabla^{:} u,\nabla^{\cdot} u\rangle Y(\xi))\,d\Sigma^{2n-1}(x,\xi).
\end{multline}

Now, we transform the last three integrals on the right-hand side by
using the Gauss--Ostrogradski\u\i{} formulas of \cite[Theorem 2.7.1]{Sh2}.
They give:
\begin{align*}
\int_{SM_\varepsilon}\divv((\mathbf X u)\nabla^{:} u)\,d\Sigma^{2n-1}
&=n \int_{SM_\varepsilon}\langle(\mathbf X u)\nabla^{:} u,\xi\rangle\,d\Sigma^{2n-1}(x,\xi)
\\
&=n \int_{SM_\varepsilon}(\mathbf X u)^2\,d\Sigma^{2n-1}
\end{align*}
because $\langle \nabla^{:} u,\xi\rangle=\mathbf X u$,
\begin{multline*}
\int_{SM_\varepsilon}\divh(\langle\nabla^{:} u,\nabla^{\cdot} u\rangle\xi
-(\mathbf X  u)\nabla^{\cdot} u) \,d\Sigma^{2n-1}
\\
=(-1)^{n-1}\int_{\partial(SM_\varepsilon)}\langle \langle\nabla^{:} u,\nabla^{\cdot} u\rangle\xi
-(\mathbf X  u)\nabla^{\cdot} u,\grad\rho(x)\rangle\,d\Sigma^{2n-2}(x,\xi)
\\
=(-1)^{n-1}\int_{\partial(SM_\varepsilon)}
\left(\langle \nabla^{:}_\rho u,\nabla^{\cdot} u\rangle \langle\xi,\grad\rho\rangle
-\langle \nabla^{:}_\rho u,\xi\rangle\langle \nabla^{\cdot}u,\grad \rho\rangle\right)
\,d\Sigma^{2n-2}
\end{multline*}
because a straightforward calculation gives
$$\langle \langle\nabla^{:} u,\nabla^{\cdot} u\rangle\xi
-(\mathbf X  u)\nabla^{\cdot} u,\grad\rho(x)\rangle
=\langle \nabla^{:}_\rho u,\nabla^{\cdot} u\rangle \langle\xi,\grad\rho\rangle
-\langle \nabla^{:}_\rho u,\xi\rangle\langle \nabla^{\cdot}u,\grad \rho\rangle,
$$
and
\begin{multline*}
\int_{SM_\varepsilon}\divv(\langle\nabla^{:} u,\nabla^{\cdot} u\rangle Y(\xi))\,d\Sigma^{2n-1}(x,\xi)
\\
=(n-1)\int_{SM_\varepsilon}\langle(\langle\nabla^{:} u,\nabla^{\cdot} u\rangle Y(\xi),\xi\rangle\,d\Sigma^{2n-1}(x,\xi)=0
\end{multline*}
because $\langle Y(\xi),\xi\rangle=0$.

Hence, \eqref{integral-1} takes the form
\begin{multline}\label{integral-2}
\int_{SM_\varepsilon}\left\{|\mathbf X(\nabla^{\cdot}u)|^2
-\langle {\mathcal C}(\nabla^{\cdot}u),\nabla^{\cdot}u\rangle
-\langle Y(\xi),\nabla^{\cdot}u\rangle^2\right\}\,d\Sigma^{2n-1}(x,\xi)
\\
=
\int_{SM_\varepsilon}|\nabla^{\cdot} (\mathbf X  u)|^2\,d\Sigma^{2n-1}
-n \int_{SM_\varepsilon}(\mathbf X u)^2\,d\Sigma^{2n-1}
\\
+(-1)^n\int_{\partial(SM_\varepsilon)}
\left(\langle \nabla^{:}_\rho u,\nabla^{\cdot} u\rangle \langle\xi,\grad\rho\rangle
-\langle \nabla^{:}_\rho u,\xi\rangle\langle \nabla^{\cdot}u,\grad \rho\rangle\right)
\,d\Sigma^{2n-2}
\end{multline}

Observe that by Lemma \ref{properties-l} and \eqref{def-u},
the derivatives $\nabla^:_\rho u$ and $\nabla^{\cdot}u$ are bounded
in $SM\setminus S(\partial M)$. This fact allows us to pass to the limit
in \eqref{integral-2}, arriving at the identity
\begin{multline*}
\int_{SM}\left\{|\mathbf X(\nabla^{\cdot}u)|^2
-\langle {\mathcal C}(\nabla^{\cdot}u),\nabla^{\cdot}u\rangle
-\langle Y(\xi),\nabla^{\cdot}u\rangle^2\right\}\,d\Sigma^{2n-1}(x,\xi)
\\
=
\int_{SM}|\nabla^{\cdot} (\mathbf X  u)|^2\,d\Sigma^{2n-1}
-n \int_{SM}(\mathbf X u)^2\,d\Sigma^{2n-1}
\\
+(-1)^n\int_{\partial(SM)}
\left(\langle \nabla^{:}_\rho u,\nabla^{\cdot} u\rangle \langle\xi,\grad\rho\rangle
-\langle \nabla^{:}_\rho u,\xi\rangle\langle \nabla^{\cdot}u,\grad \rho\rangle\right)
\,d\Sigma^{2n-2}.
\end{multline*}

Note that the last integral on the right-hand side vanishes,
because $u|_{\partial(SM)}=0$ and hence $\nabla^{\cdot}u|_{\partial(SM)}=0$.

For $(x,\xi)\in SM$,
$\mathbf X u(x,\xi)=v_i(x)\xi^i+\varphi(x)$.
As follows from Lemma 4.5.3 of \cite{Sh1}
(and is easy to check directly),
$$
\int_{SM}|\nabla^{\cdot} (v_i(x)\xi^i+\varphi(x))|^2\,d\Sigma^{2n-1}
\le n \int_{SM}(v_i(x)\xi^i+\varphi(x))^2\,d\Sigma^{2n-1}.
$$

Setting $Z=\nabla^{\cdot}u$, we thus have
$$
\int_{SM}\left\{|\mathbf XZ|^2
-\langle {\mathcal C}(Z),Z\rangle
-\langle Y(\xi),Z\rangle^2\right\}\,d\Sigma^{2n-1}\le 0.
$$
By Santal\'o's formula \eqref{santalo-f} and
the Index Lemma \ref{index-lemma} of Appendix A,
this is possible if and only if $Z=0$. So $u(x,\xi)$ is independent of $\xi$,
which clearly implies the sought conclusion.
\end{proof}

\subsection{Linear problem for $2$-tensors}

For a magnetic system $(M,g,\alpha)$
and $(x,\xi)\in SM$, put
$$
k_\mu(x,\xi)=\sup_{\eta}
\big\{2K(x,\sigma_{\xi,\eta})
+\langle Y(\eta),\xi\rangle^2
+(n+3)|Y(\eta)|^2
-2\langle(\nabla_\eta Y)(\xi),\eta\rangle\big\},
$$
where the supremum is taken over all unit vectors $\eta\in T_xM$ orthogonal to $\xi$,
and
$K(x,\sigma_{\xi,\eta})$ is the sectional curvature of
the $2$-plane $\sigma_{\xi,\eta}$
spanned by $\xi$ and $\eta$.

Define
$$
k_\mu^+(x,\xi)=\max\{0,k(x,\xi)\}
$$
and
$$
k(M,g,\alpha)=\sup_\gamma T_\gamma\int_0^{T_\gamma} k_\mu^+(\gamma(t),\dot\gamma(t))\,dt,
$$
where the supremum is taken over all unit speed magnetic geodesics
$\gamma:[0,T_\gamma]\to M$ running between boundary points.

\begin{Theorem}\label{quadratic}
If $(M,g,\alpha)$ be a simple magnetic system with
$k(M,g,\alpha)\le 4$, then $I$ is s-injective.
\end{Theorem}

\begin{Remark}
Note that the condition $k(M,g,\alpha)\le 4$ holds if $(M,g)$
is negatively curved and the $C^1$-norm of $Y$ is small enough.
Also, it is easy to see that this conditions is valid
for every sufficiently small simple piece of any magnetic system.
\end{Remark}

\begin{proof}
Assume that $[h,\beta]\in \ker I\cap\mathcal S\mathbf L^2(M)$.
Then, by Theorem \ref{thm_estimate}(a), the pair $[h,\beta]$ is smooth.
Define $u:SM\to \mathbb R$ by means of \eqref{def-u}. Then $u$ satisfies
the boundary value problem \eqref{kinetic}.  As before, we preserve
the notation $u$ for the extension of $u$ to a positively homogeneous function
of degree $0$ on $TM\setminus\{0\}$. Then for $(x,\xi)\in TM\setminus\{0\}$ we have
\begin{equation}\label{kinetic-u}
\mathbf Xu=|\xi|^{-1}h_{ij}(x)\xi^i\xi^j+\beta_j(x)\xi^j.
\end{equation}

Now, we wish to integrate identity \eqref{dif-2} over $SM$.
Again $u$ has singularities on $T(\partial M)$. We overcome this obstacle
in the same way as in the proof of the previous theorem. Thus, we deduce
\begin{multline*}
\int_{SM}\left\{|\mathbf X(\nabla^{\cdot}u)|^2
-2\langle {\mathcal C}(\nabla^{\cdot}u),\nabla^{\cdot}u\rangle
-\langle Y(\xi),\nabla^{\cdot}u\rangle^2\right\}\,d\Sigma^{2n-1}(x,\xi)
\\
=\int_{SM}|\nabla^{\cdot} (\mathbf X  u)|^2\,d\Sigma^{2n-1}
-2\int_{SM}\divv((\mathbf X u)\nabla^{:} u)\,d\Sigma^{2n-1}
-\int_{SM}|\nabla^{:} u|^2\,d\Sigma^{2n-1}
\\
+2\int_{SM}\langle\nabla^{:} u,\nabla^{\cdot}(\mathbf X u)\rangle\,d\Sigma^{2n-1}
-2\int_{SM}\divh(\langle\nabla^{:} u,\nabla^{\cdot} u\rangle\xi
-(\mathbf X  u)\nabla^{\cdot} u)\,d\Sigma^{2n-1}
\\
+\int_{SM}\divv(\langle\nabla^{:} u,\nabla^{\cdot} u\rangle Y(\xi))\,d\Sigma^{2n-1}.
\end{multline*}

Again we use the Gauss-Ostrogradski\u\i{} formulas to transform
integrals of divergent form. Denoting the leftmost side of the above formula
by $A$, we obtain
\begin{multline}\label{intt-2}
A=\int_{SM}|\nabla^{\cdot} (\mathbf X  u)|^2\,d\Sigma^{2n-1}
-\int_{SM}\left\{2n(\mathbf X u)^2+|\nabla^{:} u|^2\right\}\,d\Sigma^{2n-1}
\\
+2\int_{SM}\langle\nabla^{:} u,\nabla^{\cdot}(\mathbf X u)\rangle\,d\Sigma^{2n-1}
\\
\le \int_{SM}|\nabla^{\cdot} (\mathbf X  u)|^2\,d\Sigma^{2n-1}
-(2n+1)\int_{SM}(\mathbf X u)^2\,d\Sigma^{2n-1}
\\
+2\int_{SM}\langle\nabla^{:} u,\nabla^{\cdot}(\mathbf X u)\rangle\,d\Sigma^{2n-1},
\end{multline}
where we have used the inequality
$|\mathbf Xu|=|\langle \nabla^{:}u,\xi\rangle|\le |\nabla^{:}u|$.

From \eqref{kinetic-u}
we have for $(x,\xi)\in SM$
$$
(Xu)_{\cdot l}=2h_{lj}\xi^j-\xi_lh_{ij}\xi^i\xi^j+\beta_l.
$$
Therefore,
\begin{multline}\label{intt-3}
\int_{SM}|\nabla^{\cdot} (\mathbf X  u)|^2\,d\Sigma^{2n-1}
=\int_{SM}|2h(\xi)-\langle h,\xi^2\rangle\xi+\beta|^2\,d\Sigma^{2n-1}
\\
=\int_{SM}\big\{4|h(\xi)|^2-3\langle h,\xi^2\rangle^2+|\beta|^2
+4\langle h(\xi),\beta\rangle -\langle h,\xi^2\rangle\langle\beta,\xi\rangle\big\}\,d\Sigma^{2n-1}
\\
=\int_{SM}\big\{4|h(\xi)|^2-3\langle h,\xi^2\rangle^2+|\beta|^2\big\}\,d\Sigma^{2n-1},
\end{multline}
\begin{equation}\label{intt-4}
\int_{SM}(\mathbf X u)^2\,d\Sigma^{2n-1}
=\int_{SM}\big\{\langle h,\xi^2\rangle^2
+\langle \beta,\xi\rangle^2\big\}\,d\Sigma^{2n-1}.
\end{equation}

Now, we transform the last integral in \eqref{intt-2}.
Once
$$
\nabla_{:}u=\nabla_{|}u-Y(\nabla_{\cdot}u),
$$
it follows that
\begin{equation}\label{nunxu}
\langle\nabla^{:} u,\nabla^{\cdot}(\mathbf X u)\rangle
=\langle \nabla_{|}u,\nabla_{\cdot}(\mathbf X u)\rangle
-\langle Y(\nabla_{\cdot}u),\nabla_{\cdot}(\mathbf X u)\rangle.
\end{equation}

Observe that
\begin{multline}\label{nunxu-1}
\langle\nabla_{|} u,\nabla_{\cdot}(\mathbf X u)\rangle
=g^{kl}u_{|k}(2h_{lj}\xi^j-\xi_lh_{ij}\xi^i\xi^j+\beta_l)
\\
=2(ug^{kl}h_{lj}\xi^j)_{|k}
-2ug^{kl}h_{lj,k}\xi^j
-(u\xi^kh_{ij}\xi^i\xi^j)_{|k}
+u\xi^k(h_{ij}\xi^i\xi^j)_{|k}
+(u\beta^k)_{|k}-u\beta^k_{,k}
\\
=\divh(W)
-2u(\delta h)_j\xi^j
+u\mathbf G(h_{ij}\xi^i\xi^j)-u\delta(\beta)
\\
=\divh(W)
-(n-1)\langle uY(\beta),\xi\rangle
+u\mathbf G(\langle h,\xi^2\rangle),
\end{multline}
with
$$
W^k=2ug^{kl}h_{lj}\xi^j-u\xi^kh_{ij}\xi^i\xi^j+u\beta^k.
$$
Here we have used the equations \eqref{s-eq} for solenoidal pairs.

Integrating \eqref{nunxu-1} and transforming integrals by means of Gauss--Ostrogradski\u\i{}
formulas, we get
\begin{multline}
\int_{SM}\langle\nabla_{|} u,\nabla_{\cdot}(\mathbf X u)\rangle\,d\Sigma^{2n-1}
\\
=-(n-1)\int_{SM}\langle uY(\beta),\xi\rangle\,d\Sigma^{2n-1}
+\int_{SM}u\mathbf G(\langle h,\xi^2\rangle)\,d\Sigma^{2n-1}
\\
=-\int_{SM}\divv(uY(\beta))\,d\Sigma^{2n-1}
-\int_{SM}\mathbf G(u)\langle h,\xi^2\rangle\,d\Sigma^{2n-1}
\\
=-\int_{SM}\langle Y(\beta),\nabla^{\cdot}u\rangle\,d\Sigma^{2n-1}
-\int_{SM}\mathbf G_{\mu}(u)\langle h,\xi^2\rangle\,d\Sigma^{2n-1}
\\
+\int_{SM}(\xi^iY^j_iu_{\cdot j})\langle h,\xi^2\rangle\,d\Sigma^{2n-1}
\\
=\int_{SM}\langle Y(\nabla^{\cdot}u),\beta\rangle\,d\Sigma^{2n-1}
-\int_{SM}\mathbf \langle h,\xi^2\rangle^2\,d\Sigma^{2n-1}
\\
+\int_{SM}\langle \nabla^{\cdot}u,Y(\xi)\rangle \langle h,\xi^2\rangle\,d\Sigma^{2n-1}.
\label{nu-nxu}
\end{multline}

Next,
\begin{multline}\label{ynu-nxu}
\int_{SM}\langle Y(\nabla_{\cdot}u),\nabla_{\cdot}(\mathbf X u)\rangle\,d\Sigma^{2n-1}
\\
=\int_{SM}\langle Y(\nabla_{\cdot}u),(2h_{lj}\xi^j-\xi_lh_{ij}\xi^i\xi^j+\beta_l)\rangle\,d\Sigma^{2n-1}
\\
=2\int_{SM}\langle Y(\nabla^{\cdot}u),h(\xi)\rangle\,d\Sigma^{2n-1}
+\int_{SM}\langle \nabla^{\cdot}u,Y(\xi)\rangle\langle h,\xi^2\rangle\,d\Sigma^{2n-1}
\\
+\int_{SM}\langle Y(\nabla^{\cdot}u),\beta\rangle\,d\Sigma^{2n-1},
\end{multline}
where $h(\xi)=(g^{ki}h_{ij}\xi^j)$.

From \eqref{nunxu}, \eqref{nu-nxu}, and \eqref{ynu-nxu}
we get
\begin{multline}\label{intt-5}
\int_{SM}\langle\nabla^{:} u,\nabla^{\cdot}(\mathbf X u)\rangle\,d\Sigma^{2n-1}
\\
=-\int_{SM}\mathbf \langle h,\xi^2\rangle^2\,d\Sigma^{2n-1}
-2\int_{SM}\langle Y(\nabla^{\cdot}u),h(\xi)\rangle\,d\Sigma^{2n-1}
\\
\le -\int_{SM}\mathbf \langle h,\xi^2\rangle^2\,d\Sigma^{2n-1}
+\frac{2}{n+2}\int_{SM}|h(\xi)|^2\,d\Sigma^{2n-1}
\\
+\frac {n+2}{2}\int_{SM}|Y(\nabla^{\cdot}u)|^2\,d\Sigma^{2n-1}.
\end{multline}

Using \eqref{intt-3}, \eqref{intt-4}, and \eqref{intt-5} in
\eqref{intt-2}, we find that
\begin{multline*}
A\le \int_{SM}\big\{4|h(\xi)|^2-3\langle h,\xi^2\rangle^2+|\beta|^2\big\}\,d\Sigma^{2n-1}
\\
-(2n+1)\int_{SM}\big\{\langle h,\xi^2\rangle^2+\langle \beta,\xi\rangle^2\big\}\,d\Sigma^{2n-1}
-2\int_{SM}\mathbf \langle h,\xi^2\rangle^2\,d\Sigma^{2n-1}
\\
+\frac 4{n+2}\int_{SM}|h(\xi)|^2\,d\Sigma^{2n-1}
+(n+2)\int_{SM}|Y(\nabla^{\cdot}u)|^2\,d\Sigma^{2n-1}
\\
=\frac{4(n+3)}{n+2}\int_{SM}|h(\xi)|^2\,d\Sigma^{2n-1}
-2(n+3)\int_{SM}\langle h,\xi^2\rangle^2\,d\Sigma^{2n-1}
\\
+\int_{SM}|\beta|^2\,d\Sigma^{2n-1}
-(2n+1)\int_{SM}\langle \beta,\xi\rangle^2\,d\Sigma^{2n-1}
\\
+(n+2)\int_{SM}|Y(\nabla^{\cdot}u)|^2\,d\Sigma^{2n-1}.
\end{multline*}

As follows from Lemma 4.5.3 of \cite{Sh1}
(and is easy to check directly),
$$
\int_{SM}|h(\xi)|^2\,d\Sigma^{2n-1}\le \frac{n+2}2\int_{SM}\langle h,\xi^2\rangle^2\,d\Sigma^{2n-1},
$$
$$
\int_{SM}|\beta|^2\,d\Sigma^{2n-1}\le n\int_{SM}\langle \beta,\xi\rangle^2\,d\Sigma^{2n-1},
$$
Therefore,
\begin{equation*}
A\le (n+2)\int_{SM}|Y(\nabla^{\cdot}u)|^2\,d\Sigma^{2n-1}.
\end{equation*}

Setting $Z=\nabla^{\cdot}u$, we rewrite this as follows:
\begin{equation}\label{ineq-sm}
\int_{SM}\left\{|\mathbf X(Z)|^2
-2\langle {\mathcal C}(Z),Z\rangle
-\langle Y(Z),\xi\rangle^2
-(n+2)|Y(Z)|^2\right\}\,d\Sigma^{2n-1}\le 0.
\end{equation}

We need the following lemma whose proof is given below.

\begin{Lemma}\label{index-2}
Let  $\gamma:[0,T]\to M$ be a unit speed magnetic geodesic.
If
$$
\int_0^T k_\mu^+(\gamma(t),\dot\gamma(t))\,dt\le 4/T,
$$
then for every smooth vector field  $Z$ along $\gamma$
vanishing at the endpoints of~$\gamma$ and orthogonal to~$\dot\gamma$ we have
\begin{equation}\label{index-in-2}
\int_0^T\left\{|\dot Z|^2
-2\langle {\mathcal C}(Z),Z\rangle
-\langle Y(Z),\dot\gamma\rangle^2
-(n+2)|Y(Z)|^2\right\}\,dt\ge 0
\end{equation}
with equality if and only if $Z=0$.
\end{Lemma}

Continuing the proof of the theorem, observe that, by
Santal\'o's formula \eqref{santalo-f} and
Lemma \ref{index-2}, inequality \eqref{ineq-sm} may hold
if and only if $Z=0$. So $\nabla^{\cdot}u=0$, which means that
$u(x,\xi)$ is independent of $\xi$. Since the pair $[h,\beta]$ is solenoidal,
this readily yields the sought conclusion.
\end{proof}

\begin{proof}[Proof of Lemma \ref{index-2}]
Denote the left-hand side of \eqref{index-in-2} by $B$.
Since
\begin{equation*}
{\mathcal C}(Z)=\mathbf R_{\dot{\gamma}}(Z)- Y(\dot{Z})-(\nabla_{Z}Y)(\dot{\gamma}),
\end{equation*}
we have
\begin{multline*}
B=\int_0^T\left\{|\dot Z|^2
-2\langle\mathbf R_{\dot{\gamma}}(Z),Z\rangle
+2\langle Y(\dot{Z}),Z\rangle
+2\langle(\nabla_{Z}Y)(\dot{\gamma}),Z\rangle \right.
\\
-\left.\langle Y(Z),\dot\gamma\rangle^2
-(n+2)|Y(Z)|^2\right\}\,dt.
\end{multline*}

Let $e_i$, $i=1,\dots,n$ be an orthonormal frame at $\gamma(0)$,
with $e_n=\dot \gamma(0)$.
Let $E_i(t)$, $0\le t\le T$, be a vector field along $\gamma$
satisfying the equation
$$
\dot E_i=Y(E_i)
$$
with the initial condition $E_i(0)=e_i$. Surely, $E_n=\dot\gamma$.
Since
$$
\frac d{dt}\langle E_i,E_j\rangle
=\langle \dot E_i,E_j\rangle+\langle E_i,\dot E_j\rangle
=\langle Y(E_i), E_j\rangle+\langle E_i,Y(E_j)\rangle=0,
$$
we see that $E_1(t),\dots,E_n(t)$ is an orthonormal frame
for each $t$.

Consider the expansion of $Z$ in this frame:
$$
Z(t)=z^\kappa(t) E_\kappa(t),\quad \kappa\in\{1,\dots,n-1\}.
$$
Then
$$
\dot Z=\dot z^\kappa E_\kappa+z^\kappa\dot E_\kappa
=\dot z^\kappa E_\kappa+z^\kappa Y(E_\kappa)
=\dot z^\kappa E_\kappa+Y(Z),
$$
$$
|\dot Z|^2=\sum_\kappa\dot z^\kappa \dot z^\kappa
+2Y_{\kappa\lambda}z^\kappa\dot z^\lambda +|Y(Z)|^2,
$$
$$
\langle Y(\dot Z),Z\rangle
=\langle\dot z^\kappa Y(E_\kappa)+Y^2(Z),
z^\lambda E_\lambda\rangle
=-Y_{\kappa\lambda} z^\kappa \dot z^\lambda-|Y(Z)|^2
$$
with
$$
Y_{\kappa \lambda}=\langle Y(E_\kappa),E_\lambda\rangle.
$$

Therefore,
\begin{multline*}
B=\int_0^T\Big\{\sum_\kappa\dot z^\kappa \dot z^\kappa
-\big[2K(\sigma_{\dot\gamma,Z})|Z|^2
+\langle Y(Z),\dot\gamma\rangle^2
\\
+(n+3)|Y(Z)|^2
-2\langle(\nabla_{Z}Y)(\dot{\gamma}),Z\rangle\big]\Big\}\,dt
\\
\ge
\int_0^T\Big\{\sum_\kappa\dot z^\kappa \dot z^\kappa
-k_\mu|Z|^2\Big\}\,dt
=\int_0^T\sum_\kappa\Big\{\dot z^\kappa \dot z^\kappa
-k_\mu z^\kappa z^\kappa\Big\}\,dt.
\end{multline*}

Now, the claim follows from Lyapunov's inequality \cite[p.~346]{Ha}.
\end{proof}

\section{General rigidity theorems}\label{sec_6}

\subsection{Rigidity in a given conformal class}
We first state a rigidity theorem, fixing the conformal class of
a metric. The theorem below generalizes the corresponding
well-known result for the ordinary boundary rigidity problem,
see \cite{Cr1, Mu, MR}.

\begin{Theorem}\label{conf}
Let $(g,\alpha)$ and $(g',\alpha')$ be simple
magnetic systems on $M$ whose boundary action functions
$\mathbb A|_{\partial M\times \partial M}$
and $\mathbb A'|_{\partial M\times \partial M}$ coincide.
If $g'$ is conformal to $g$, i.e., $g'=\omega^2(x)g$
for a~smooth positive function $\omega$ on $M$,
then $\omega\equiv 1$ and $\alpha'=\alpha+dh$ for some
smooth function $h$ on $M$ vanishing on $\partial M$,
hence $(g',\alpha')$ is gauge equivalent to $(g,\alpha)$.
\end{Theorem}

\begin{proof}
In view of Lemma \ref{bnd-0}, $\omega=1$ on the boundary of $M$.
Next, using Lemma \ref{action-scat}, we see that the scattering
relations $\mathcal S$ and $\mathcal S'$
of both magnetic systems coincide.

Let us show that $\omega=1$ on the whole of $M$.
Note that for $(x,\xi)\in\partial_+SM$
\begin{equation}\label{a'}
\mathbb A'(\gamma'_{x,\xi})\le\mathbb A'(\gamma_{x,\xi})
=\frac12\int_0^{\ell(x,\xi)} \omega^2(\gamma_{x,\xi}(t))\,dt +\frac12 \ell(x,\xi)
-\int_{\gamma_{x,\xi}}\alpha'.
\end{equation}

Using \eqref{vol} and \eqref{int-a}, whence we obtain
\begin{align*}
\vol_{g'}(M)
&=\frac{1}{\omega_{n-1}}\int_{\partial_+SM}
\mathbb A'(\gamma'_{x,\xi})\,d\mu(x,\xi)
\\
&\le\frac{1}{2\omega_{n-1}}\int_{\partial_+SM}
\bigg\{\int_0^{\ell(x,\xi)} \omega^2(\gamma_{x,\xi}(t))\,dt +\ell(x,\xi)\bigg\}
\,d\mu(x,\xi)
\\
&=\frac 12 \int_M \omega^2\,d\vol+\frac 12 \vol_g(M).
\end{align*}

On the one hand, by H\"older's inequality
\begin{equation}\label{hoelder}
\int_M \omega^2\,d\vol\leq \left\{\int_M \omega^n\,d\vol\right\}^{\frac 2 n}
\left\{\int_M\,d\vol\right\}^{\frac {n-2} n}
=\vol_{g'}(M)^{\frac 2n}\vol_g(M)^{\frac {n-2} n},
\end{equation}
with equality if and only if $\omega\equiv 1$.

It follows that
\begin{equation}\label{v-v}
\vol_{g'}(M)\le\frac 12\vol_{g'}(M)^{\frac 2n}\vol_g(M)^{\frac {n-2} n}
+\frac 12 \vol_g(M).
\end{equation}
However, by Theorem \ref{volume}, $\vol_{g'}(M)=\vol_g(M)$,
which implies that \eqref{v-v} holds with the equality sign.
This means that \eqref{hoelder} holds with the equality sign.
Thus, $\omega\equiv 1$.

Now, \eqref{a'} and the equality
$\mathbb A'(\gamma'_{x,\xi})=\mathbb A(\gamma_{x,\xi})$ yield
$$
\int_{\gamma_{x,\xi}}\alpha
\ge \int_{\gamma_{x,\xi}}\alpha'.
$$
In view of \eqref{santalo-f} and \eqref{int-a},
we conclude now that
$$
\int_{\gamma_{x,\xi}}(\alpha'-\alpha)=0
$$
for all $(x,\xi)\in\partial_+SM$.
By Theorem \ref{linear}, we see that $\alpha'-\alpha=dh$
for some function $h$ on $M$ vanishing on
$\partial M$. This completes the proof of the theorem.
\end{proof}

\subsection{Rigidity of analytic systems}
\begin{Theorem}\label{analytic}
If $M$ is a real-analytic compact manifold with boundary,
and $(g,\alpha)$ and $(g',\alpha')$ are simple
real-analytic magnetic systems on $M$ with the same boundary action function,
then these systems are gauge equivalent.
\end{Theorem}

\begin{proof}
Consider the diffeomorphism $f:M\to M$ constructed
in the proof of Theorem \ref{jet}.
Recall that in a neighborhood of $\partial M$,
$f$ coincides with  the map
$\exp_{\partial M}\circ\,(\exp'_{\partial M})^{-1}$,
which is obviously analytic for analytic systems.

As shown in that proof, $g$ and $f^*g'$ have the same jets on $\partial M$.
Hence, by analyticity, $g$ and $f^*g'$ coincide in some
connected neighborhood of $\partial M$, and
$f$ is thus an analytic isometry of a neighborhood of $\partial M$
of the the analytic Riemannian manifold $(M,g)$
onto a neighborhood of $\partial M$
of the analytic Riemannian manifold $(M,g')$, which is moreover
the identity when restricted to $\partial M$.
Now, the same arguments as in the proof of \cite[Theorem C(a)]{LeU}
show that $f$ extends from a neighborhood of $\partial M$
to an analytic isometry $\tilde f$ of $(M,g)$ onto $(M,g')$.
Now, the sought conclusion follows from Theorem \ref{conf} applied
to the magnetic systems $(M,g,\alpha)$ and $(M,\tilde f^*g',\tilde f^*\alpha')$.
\end{proof}

\subsection{Rigidity of reversible systems}
A magnetic system $(M,g,\alpha)$ is said to be reversible if
the flip $(x,\xi)\mapsto (x,-\xi)$ conjugates $\psi^t$ with $\psi^{-t}$.
It is easy to see that a system is reversible if and only if $d\alpha=0$, i.e.,
if and only if the Lorentz force $Y$ of the system vanishes.
In this case magnetic geodesics are nothing but the ordinary geodesics
of the Riemannian manifold $(M,g)$; moreover, simplicity of
a reversible magnetic system
$(M,g,\alpha)$ is equivalent to simplicity of the Riemannian manifold $(M,g)$.

It is interesting to know whether reversibility of a magnetic system
can be established by boundary measurements. Observe that
if a system is reversible then
\begin{equation}\label{reverse}
\mathcal S(-\mathcal S(x,\xi))=(x,-\xi)
\end{equation}
for all $(x,\xi)\in\partial_+SM$. Let us call the systems satisfying
\eqref{reverse} {\em boundary reversible}.

\begin{Theorem}\label{thm_reverse}
A simple magnetic system is boundary reversible if and only if it is reversible.
\end{Theorem}

\begin{proof}\label{}
Let $(M,g,\alpha)$ be a boundary reversible simple magnetic system.
Consider the magnetic system $(g,-\alpha)$ on $M$.
Note that if $\gamma:[0,T]\to M$ is a unit speed magnetic geodesic
of the system $(g,\alpha)$, then the curve $\bar{\gamma}:[0,T]\to M$,
defined as $\bar\gamma(t)=\gamma(T-t)$, is a unit speed magnetic geodesic
of the system $(g,-\alpha)$.

For $x,y\in M$, let
$\gamma_{x,y}:[0,T_{x,y}]\to M$ denote the unit speed magnetic geodesic
of the system $(g,\alpha)$ from $x$ to $y$.
By Lemma \ref{minimize},
\begin{equation*}
\mathbb A(\gamma_{x,y})\le \mathbb A(\bar\gamma_{y,x}),
\end{equation*}
i.e.,
\begin{equation}\label{reverse-1}
T_{x,y}-\int_{\gamma_{x,y}}\alpha\le T_{y,x}+\int_{\gamma_{y,x}}\alpha.
\end{equation}
Interchanging $x$ and $y$, we receive
\begin{equation}\label{reverse-2}
T_{y,x}-\int_{\gamma_{y,x}}\alpha\le T_{x,y}+\int_{\gamma_{x,y}}\alpha.
\end{equation}

From \eqref{reverse-1} and \eqref{reverse-2} we get
\begin{equation}\label{reverse-3}
\int_{\gamma_{x,y}}\alpha+\int_{\gamma_{y,x}}\alpha\ge0.
\end{equation}

For $(x,\xi)\in\partial_+SM$ and $y=\gamma_{x,\xi}(\ell(x,\xi))$,
we have in view of \eqref{reverse}
$$
\gamma_{x,y}=\gamma_{x,\xi},
\quad
\bar\gamma_{y,x}=\tilde\gamma_{x,\xi},
$$
where $\tilde\gamma_{x,\xi}:[0,\tilde\ell(x,\xi)]\to M$ is the unit speed magnetic
geodesic of the magnetic system $(M,g,-\alpha)$ with
initial conditions $\tilde\gamma_{x,\xi}(0)=x$,
$\dot{\tilde\gamma}_{x,\xi}(0)=\xi$,
and with $\tilde\gamma_{x,\xi}(\tilde\ell(x,\xi))=y$.
Therefore, \eqref{reverse-3} yields
\begin{equation}\label{reverse-ineq}
\int_{\gamma_{x,\xi}}\alpha-\int_{\tilde\gamma_{x,\xi}}\alpha\ge0.
\end{equation}

By Santal\'o's formula \eqref{santalo-f}
and \eqref{int-a}, we have
\begin{equation}\label{reverse-santalo}
\int_{\partial_+SM}\left\{\int_{\gamma_{x,\xi}}\alpha
-\int_{\tilde\gamma_{x,\xi}}\alpha\right\}\,d\mu(x,\xi)=0.
\end{equation}
Combining \eqref{reverse-santalo} and \eqref{reverse-ineq} yields
\begin{equation}\label{reverse-eq}
\int_{\gamma_{x,\xi}}\alpha=\int_{\tilde\gamma_{x,\xi}}\alpha
\end{equation}
for all $(x,\xi)\in\partial_+SM$.
From \eqref{reverse-1} we then obtain
$$
\ell(x,\xi)\le \tilde\ell(x,\xi),
$$
and from \eqref{reverse-2},

$$
\tilde\ell(x,\xi)\le\ell(x,\xi).
$$
Thus,
\begin{equation}\label{ell-eq}
\tilde\ell(x,\xi)=\ell(x,\xi).
\end{equation}

Now, by \eqref{reverse-eq} and \eqref{ell-eq}
$$
\mathbb A(\tilde\gamma_{x,\xi})=\mathbb A(\gamma_{x,\xi})= \mathbb A(x,y).
$$
This implies that $\tilde\gamma_{x,\xi}$ is a unit speed magnetic geodesic
of the system $(g,\alpha)$ which joins $x$ to $y$.
By simplicity, we then have
$$
\tilde\gamma_{x,\xi}=\gamma_{x,\xi}.
$$
Valid for all $(x,\xi)\in\partial_+SM$, this equality means clearly
that the magnetic system $(M,g,\alpha)$ is reversible.
\end{proof}

\begin{Theorem}\label{thm_reverse_rig}
A simple reversible magnetic system $(M,g,\alpha)$
is magnetic boundary rigid
if and only if the simple Riemannian manifold $(M,g)$ is boundary rigid.
\end{Theorem}

\begin{proof}
Suppose that a simple reversible magnetic system $(g,\alpha)$ on $M$
is magnetic boundary rigid. Let $g'$ be a simple Riemannian metric on
$M$ whose boundary distance function equals the boundary distance function
of $g$. We must prove that $g'$ is isometric to $g$ by an isometry
which is the identity on the boundary.

Without loss of generality we may assume (see, e.g., \cite[Theorem 2.1]{LSU}) that
\begin{equation}\label{bnd}
g'|_{\partial M}=g|_{\partial M}.
\end{equation}
Then equality of the boundary distance functions of $g$ and $g'$
and \eqref{bnd} are well known to imply
equality of the scattering relations of the Riemannian manifolds $(M,g)$
and $(M,g')$ (see, e.g., \cite{M}).

Consider the magnetic system $(g',\alpha)$ on $M$.
Since $d\alpha=0$, this system is reversible as well.
It is trivial that for a reversible system, the scattering relation
of the magnetic system coincides with
the scattering relation of the underlying Riemannian metric.
Therefore,  the magnetic systems $(g,\alpha)$
and $(g',\alpha)$ on $M$ have the same scattering relations.
In view of \eqref{bnd} and Lemma \ref{scat-action},
the boundary actions functions of these systems coincide.
From the magnetic boundary rigidity of $(M,g,\alpha)$ we infer that
the metric $g'$ is isometric to $g$ by an isometry which is the identity
on the boundary, as required.

Now, let $(M,g,\alpha)$ be a simple reversible
magnetic system such that the simple Riemannian manifold
$(M,g)$ is boundary rigid.
Let $(g',\alpha')$ be a simple magnetic system on $M$ whose
boundary action function $\mathbb A'|_{\partial M\times\partial M}$
equals the boundary action function $\mathbb A|_{\partial M\times\partial M}$
of the magnetic system $(g,\alpha)$. By Theorem \ref{jet},
we may assume that
$g'|_{\partial M}=g|_{\partial M}$ and
$\alpha'|_{\partial M}=\alpha|_{\partial M}$.
Now, Lemma \ref{action-scat} implies that the scattering relations
of the magnetic systems $(M,g,\alpha)$ and $(M,g',\alpha')$
are the same. Since $(M,g,\alpha)$ is reversible
and hence boundary reversible, $(M,g',\alpha')$ too
is boundary reversible and hence, by
Theorem \ref{thm_reverse}, it is also reversible.
This implies readily that the scattering relations
of the Riemannian manifolds $(M,g)$ and $(M,g')$ are the same.
Applying Lemma \ref{scat-action} to the magnetic systems
$(M,g,0)$ and $(M,g',0)$, we see that the boundary distance functions
of the metrics $g$ and $g'$ are the same. Hence, by the boundary rigidity
of the Riemannian manifold $(M,g)$, there is a diffeomorphism
$f:M\to M$, $f|_{\partial M}= \text{identity}$, such that
$g'=f^*g$. Then Theorem \ref{conf}, applied to the magnetic systems
$(M,g,\alpha)$ and $(M,f^*g',f^*\alpha')$,
gives the gauge equivalence of the magnetic systems
$(g,\alpha)$ and $(g',\alpha')$, as required.
\end{proof}

\subsection{Generic local boundary rigidity}
We will prove that near each $(g_0,\alpha_0)$ in
the generic set $\mathcal{G}^k$
of Definition \ref{generic-set},
the action  on the boundary determines $(g,\alpha)$.
Note that,  by Theorem \ref{quadratic}, this generic set $\mathcal G^k$
contains the magnetic systems $(g,\alpha)$ with $k(M,g,\alpha)\le 4$,
in particular magnetic systems in which the underlying
Riemannian metric is negatively curved and
the magnetic field is sufficiently small.

\begin{Theorem}  \label{S_thm_gen_rigidity}
Let $k_0$ be as in Theorem~\ref{S_thm_lin}.
There exists $k\ge k_0$ such that for every
$(g_0,\alpha_0) \in \mathcal{G}^k$, there is $\varepsilon>0$
such that for any two magnetic systems $(g,\alpha)$, $(g',\alpha')$ with
\[
\|g-g_0\|_{C^k(M)}+ \|\alpha-\alpha_0\|_{C^k(M)}\le \varepsilon,\quad
\|g'-g_0\|_{C^k(M)}+ \|\alpha'-\alpha_0\|_{C^k(M)}\le \varepsilon
\]
we have the following:
\begin{equation}  \label{S_pro}
\mbox{$\mathbb{A}_{g,\alpha}=\mathbb{A}_{g',\alpha'}$}  \quad \text{on $\partial M\times \partial M$}
\end{equation}
implies that $(g,\alpha)$ and $(g',\alpha')$ are gauge equivalent,
i.e., $g'=\psi^* g$, $\alpha'=\psi^* \alpha+d\phi$
with some $C^{k+1}(M)$-diffeomorphism $\psi:M \to M$,
fixing the boundary, and some $C^{k+1}$ function $\phi$
vanishing on $\partial M$.
\end{Theorem}

Observe that  $g$ is solenoidal (with respect to itself),
and if we replace $\alpha$ by its solenoidal projection
$\alpha^s := \alpha-d\phi$, where
$\Delta_g\phi =\delta\alpha$, $\phi|_{\partial M}=0$,
then $\delta \alpha^s=0$ as well.
Therefore, $[\frac12 g,-\alpha^s]$ is not necessarily
a solenoidal pair in the sense of Definition~\ref{s-def}, instead
\begin{equation}  \label{nl_1}
\delt \left[\frac12 g,-\alpha^s\right] = \left[\frac{n-1}{2} Y(\alpha^s),0\right].
\end{equation}

We will prove an analogue of \cite[Theorem~2.1]{CDS}.

\begin{Lemma}  \label{thm_shift}
Let $(g,\alpha) \in C^{k,\mu}$, $0<\mu<1$, $k\ge2$,
be a simple magnetic system on $M$ with $\delta \alpha=0$.
Then for any other $(g',\alpha')$ close enough to $(g,\alpha)$,
there exists a magnetic system $(\tilde g',\tilde \alpha')$
gauge equivalent to $(g',\alpha')$ and satisfying
\eqref{nl_1}, i.e.,
\begin{equation}  \label{nl_2}
\delt \left[\frac12 \tilde g',-\tilde \alpha'\right] = \left[\frac{n-1}{2} Y(\alpha),0\right],
\end{equation}
where $\delt$ and $Y$ are related to $(g,\alpha)$.

Moreover, if\/ $\|(g',\alpha') - (g,\alpha)\|_{C^{k,\mu}} \le \varepsilon$
with $\varepsilon\ll1$, then
$\|(\tilde g', \tilde\alpha') - (g,\alpha)\|_{C^{k,\mu}} \le \varepsilon_1$,
with $\varepsilon_1=\varepsilon_1(\varepsilon)\to0$ as $\varepsilon\to0$.
\end{Lemma}

\begin{proof}
Our argument is much the same as that in the proof of \cite[Theorem~2.1]{CDS}.
So we merely sketch it.

We have to show that there exist a diffeomorphism $f$ of $M$,
fixing $\partial M$, and a function $\varphi$,
vanishing on $\partial M$, such that $\tilde g'=f^*g'$,
$\tilde \alpha'=f^*\alpha'+d\varphi$, and \eqref{nl_2} holds.
If $f$ is close enough to the identity, we can identify $f$
with a certain vector field $v$ as follows.
If $v$ is a vector field vanishing on $\partial M$
and $|\nabla v|_g<1$, the map
\[
e_v(x) = \exp_x(v(x))
\]
is well defined on $M$, with image in $M$ again,
and it is a diffeomorphism for $\|v\|_{C^{k,\mu}}$ small enough.
Then $v\mapsto e_v$ has an inverse  defined as follows.
If $f$ is a diffeomorphism close enough to the identity, set
$v_f(x) = \exp^{-1}_x(f(x))$, i.e,
$v_f(x) = \dot\gamma(0)$, where $\gamma: [0,1]\to M$ is the  geodesic
such that $\gamma(0)=x$, $\gamma(1)=f(x)$.
The existence of such a geodesic follows from
the strict convexity of the boundary, following in turn
from the simplicity assumption.
Clearly, the differential of the map $v\mapsto e_v$ at $v=0$
is the identity transformation.

Let $h=\frac12(g'-g)$, $\beta=-(\alpha'-\alpha)$.
Condition \eqref{nl_2}  then takes the form
\begin{align}   \label{nl_3}
\delta\left( e_v^*(g/2+h) \right)
+\frac{n-1}{2} Y\left((e_v^*(\alpha+\beta)+d\varphi \right)&
=   \frac{n-1}{2} Y(\alpha),
\\
-\delta\left( e_v^*(\alpha+\beta) +d\varphi\right)&=0. \label{nl_4}
\end{align}

We define the map (recall that $(g,\alpha)$ is fixed)
\begin{multline*}\label{def-f}
F([v,\varphi],[h,\beta])
\\
=\left[\delta\left( e_v^*(g/2+h) \right)
+\frac{n-1}{2} Y\left((e_v^*(\alpha+\beta)+d\varphi \right),
-\delta\left( e_v^*(\alpha+\beta) +d\varphi\right)\right].
\end{multline*}
Then \eqref{nl_3}, \eqref{nl_4} can be rewritten as
\begin{equation}\label{def-f}
F([v,\varphi],[h,\beta]) = \left[\frac{n-1}{2} Y(\alpha),0\right].
\end{equation}
We want to solve this equation for $[v,\varphi]$
if $[h,\beta]$ are small enough
and $g$, $\alpha$ are fixed. If $[h,\beta]=[0,0]$, then $[v,\varphi]=[0,0]$
is a solution by \eqref{nl_1}.
To show solvability for $[h,\beta]$ small enough,
we apply an implicit function theorem in Banach spaces as in \cite{CDS}.
For this purpose, we need to compute the  derivative
$F'_{[v,\varphi]}([0,0],[0,0])$.

Set $[h,\beta]=0$ in \eqref{def-f} and use
\eqref{S3bb} to deduce
\begin{align*}
F'_{[v,\varphi]}&([0,0],[0,0])[v,\varphi]\\
 &= \Big[ \delta \d v +\frac{n-1}{2}Y( Y(v)+d\langle v,\alpha\rangle +d\varphi ),
 -\delta(Y(v)+ d \langle v, \alpha\rangle +d\varphi )  \Big]\\
  &= \delt\mathbf{d}[v,  -\langle v,\alpha\rangle -\varphi],
\end{align*}
where $v$ is considered as a $1$-form by lowering the index.
Thus, the derivative above is the superposition of the map
$[v,\varphi] \mapsto [v, - \langle v,\alpha\rangle -\varphi]$,
which is an isomorphism (and this map and its inverse preserve
the zero boundary conditions), and $\delt\mathbf{d}$
with Dirichlet boundary conditions, which is also invertible in
appropriate spaces. We refer to \cite{CDS} for more technical details.
\end{proof}

The next lemma states, loosely speaking,  that gauge equivalent pairs differ by a potential one, modulo a quadratic term.

\begin{Lemma}   \label{pr_dv}
Let $(g',\alpha')$ and $(g,\alpha)$ be in $C^k$, $k\ge2$,
and gauge equivalent, i.e.,
\[
g'=\psi^*g, \quad \alpha'=\psi^*\alpha+d\varphi
\]
for some diffeomorphism $\psi$ fixing $\partial M$ and some function
$\varphi$ be  vanishing on $\partial M$.
Set $\mathbf{f} = [\frac12(g'-g),-(\alpha'-\alpha)]$.
Then there exists $\mathbf{w}$, vanishing on $\partial M$, such that
\[
\mathbf{f} = \mathbf{dw} +\mathbf{f}_2,
\]
and for  $(g,\alpha)$ belonging to any bounded set $U$ in $C^k$,
there exists $C(U)>0$ such that
\[
\|\mathbf{f}_2\|_{C^{k-2}}\le C(U)\|\psi-\text{\rm Id}\|^2_{C^{k-1}},
\quad
\|\mathbf{w}\|_{C^{k-1}} \le C(U)\|\psi-\text{\rm Id}\|_{C^{k-1}}.
\]
\end{Lemma}

\begin{proof}
As in the proof above,
set $v(x)= \exp^{-1}_x(\psi(x))$, which is a well defined vector field
if $\psi$ is close enough to the identity in $C^2$
(it is enough to prove the claim in this case only)
and $v=0$ on $\partial $M.
Set $\psi_\tau(x)  = \exp_x(\tau v(x))$, $0\le \tau\le1$.
Let $g^\tau = \psi_\tau^* g$. Then the Taylor formula implies
\[
g' = g+ \frac{d}{d\tau}\Big|_{\tau=0} g^\tau +h = g+ 2\d v+h,
\]
where
\[
|h|\le \frac12 \max_{\tau\in[0,1]} \Big|\frac{d^2 g^\tau}{d\tau^2}\Big|,
\]
and $2\d v$ is the linearization of $g^\tau$ at $\tau=0$.
To estimate $h$, write
\[
g^\tau_{ij} = g_{kl}\circ \psi_\tau \frac{\partial \psi^k_\tau}{\partial x^i} \frac{\partial \psi^l_\tau}{\partial x^j},
\]
and differentiate twice w.r.t.\ $\tau$. Notice that
\[
\Big|\frac{\partial^2 \psi_\tau}{\partial\tau^2}\Big| \le C\|v\|^2_{L^\infty}, \quad
\Big|\frac{\partial^2 \nabla \psi_\tau}{\partial\tau^2}\Big| \le C\|v\|^2_{C^1}.
\]
This yields the stated estimate  for the first component of $\mathbf{f}_2$
for $k=2$. The estimates for $k>2$ go along similar lines
by expressing the remainder $h$ in its Lagrange form
and estimating the  derivatives of $h$.

The analysis of the second component is similar, using \eqref{S3bb}.
In particular, we get that $\mathbf{w}= [v,-\alpha(v)-\phi]$
which corresponds well with the linearization formula \eqref{S3bb}.
\end{proof}

\begin{proof}[Proof of Theorem~\ref{S_thm_gen_rigidity}]
We can assume that $(g',\alpha')$ is replaced by its gauge equivalent pair
that satisfies the assumptions of Lemma~\ref{thm_shift}
(where $\delt$ and $Y$ are related to $(g,\alpha)$).
Set $\mathbf{f} = [(g'-g)/2,-(\alpha'-\alpha)]$.
Then $\delt\mathbf{f}=0$. We assume that $(g,\alpha)$ belongs to
a small enough neighborhood of $(g_0,\alpha_0)$,
so that $(g,\alpha)\in \mathcal{G}^k$, and the constant $C$ in
Theorem~\ref{thm_estimate} is uniform in that neighborhood.
By the second statement of Lemma~\ref{thm_shift}
and the assumptions of Theorem~\ref{S_thm_gen_rigidity},
\begin{equation}  \label{nl_5}
\|\mathbf{f}\|_{C^k} \le \varepsilon_1(\varepsilon),
\end{equation}
where $\varepsilon_1\to0$ as $\varepsilon\to0$, and $k\gg1$ is fixed
(this requires the original $(g',\alpha')$
to be in $\mathcal{G}^{k+1}$ if we want to avoid the $C^{k,\mu}$ spaces).
Moreover, $\|g\|_{C^k} +\|\alpha\|_{C^k} \le A$,
where $A>0$ depends on $(g_0,\alpha_0)$ and on an upper bound
$\varepsilon_0$ of $\varepsilon$.
All constants $C$ below will depend only on $A$
and will be uniform in $\varepsilon\le \varepsilon_0$.

Let $\psi$ be a diffeomorphism in $M$ that maps boundary normal coordinates
w.r.t.\ $g$ into boundary normal coordinates w.r.t.\ $g'$
near $\partial M$, then extended to the whole $M$. In other words,
if $\psi_1$ that maps a neighborhood of $\partial M$
into $\partial M\times [0,\delta]$
defines the semigeodesic normal coordinates related to $g$ and if
$\psi_2$ is defined in the same way corresponding to $g'$,
then $\psi = \psi_2\circ \psi_1^{-1}$.
Clearly,
\begin{equation}  \label{nl_5aa}
\|\psi-\text{\rm Id}\|_{C^{k-1}} \le C\|g'-g\|_{C^k} \le C'\|\mathbf{f}\|_{C^k}.
\end{equation}
Set $\tilde g' = \psi^*g'$, $\tilde \alpha' = \psi^*\alpha'+d\varphi$,
where $\varphi$ is such that in boundary normal coordinates
$\tilde \alpha'_n=0$, see Lemma~\ref{lemma_S3}.
Set $\tilde{\mathbf{f}}= [\tilde g'/2,-\tilde\alpha'] - [g/2,-\alpha]$.
By Lemma~\ref{pr_dv} and \eqref{nl_5aa},
\begin{equation}   \label{nl_5a}
\tilde{\mathbf{f}} =\mathbf{f} + \mathbf{dw} +\mathbf{h}, \quad
\|\mathbf{h}\|_{C^{k-2}} \le C\|\mathbf{f}\|_{C^{k}}^2,
\end{equation}
with $\mathbf{w}|_{\partial M}=0$; so, roughly speaking,
$\tilde{\mathbf{f}}$ and $\mathbf{f}$ differ,
up to a quadratic term, only by a potential term.
By Theorem~\ref{jet}, $\tilde{\mathbf{f}}$ vanishes on $\partial M$
together with its derivatives up to any fixed order $m$, if $k\gg1$.
We will choose $m$ below.
We are going to prove that ${\mathbf{f}} =0$, and this would prove the theorem.

Since
$\mathbb{A}_{g,\alpha}=  \mathbb{A}_{ g',  \alpha'}
= \mathbb{A}_{\tilde g',\tilde \alpha'}$ on $\partial M\times \partial M$,
by the linearization Lemma \ref{var}
\[
\|I\tilde{\mathbf{f}}\|_{L^\infty} \le C\|\tilde {\mathbf{f}}\|^2_{C^1}.
\]
By \eqref{S12},
$\|N\tilde{\mathbf{f}}\|_{L^\infty(M_1)} \le C\|I\tilde{\mathbf{f}}\|_{L^\infty} $;
therefore,
\begin{equation}  \label{nl_6}
\|N\tilde{\mathbf{f}}\|_{L^\infty(M_1)} \le C\|\tilde{\mathbf{f}}\|^2_{C^1}.
\end{equation}
We apply Theorem~\ref{thm_estimate} to get
\begin{equation}   \label{nl_7a}
\|\tilde{\mathbf{f}}^s\|_{L^2(M)}\le  C\|N\tilde{\mathbf{f}} \|_{ H^2(M_1)},
\end{equation}
where we estimated the $\tilde H^2$ norm by the $H^2$ norm,
which is finite for $\tilde{\mathbf{f}}$. Note that for the solenoidal projection $\tilde{\mathbf{f}}^s$ of $\tilde{\mathbf{f}}$ we have, by \eqref{nl_5a},
\begin{equation}   \label{nl_7b}
\tilde{\mathbf{f}}^s = \mathbf{f} +\mathbf{h}^s;
\end{equation}
in particular,
\begin{equation}   \label{nl_7c}
\|\tilde{\mathbf{f}}^s\|_{L^2(M)}\ge \|{\mathbf{f}}\|_{L^2(M)} - \|\mathbf{h}^s\|_{L^2(M)}\ge \|{\mathbf{f}}\|_{L^2(M)} -C\|\mathbf{f}\|_{C^{2}}^2.
\end{equation}

Using interpolation estimates in $H^s(M_1)$,
see \cite[Theorem 4.3.1/1]{Tr}, and \eqref{nl_6}, we get
\begin{equation}  \label{nl_8}
\|N\tilde{\mathbf{f}} \|_{ H^2(M_1)} \le
C\|N\tilde{\mathbf{f}} \|_{ H^s(M_1)}^{2/s} \|N \tilde{\mathbf{f}} \|_{ L^2(M_1)}^{1-2/s}  \le
C'\|\tilde{\mathbf{f}}\|^{2-4/s}_{C^1}, \quad s>2.
\end{equation}
To estimate the first factor in the middle term,
we used  the fact that $\partial^\alpha \tilde{\mathbf{f}}|_{\partial M}=0$
for $|\alpha|\le m$ and therefore, if $m\ge s-2$,
\[
\|N\tilde{\mathbf{f}} \|_{ H^s(M_1)} \le
C\| \tilde{\mathbf{f}}  \|_{ H^{s-1}(M)}\le
C'\|\mathbf{f}\|_{C^{s+1}(M)}\le C''
\]
if $k\ge s+1$. The last inequality follows by comparing $\mathbf{f}$, $\tilde{\mathbf{f}}$  using \eqref{nl_5a}
and the estimate on $\mathbf{w}$ in Lemma~\ref{pr_dv}, combined with \eqref{nl_5aa}.
Combine \eqref{nl_7a}, \eqref{nl_7c}, and \eqref{nl_8} to get
\[
\|\mathbf{f}\|_{L^2(M)}\le C\left( \|\mathbf{f}\|^{2-4/s}_{C^1} +\|\mathbf{f}\|_{C^{2}}^2\right).
\]
Using again the fact that $\|\mathbf{f}\|_{C^k(M)}\le C$ ($k\gg1$), Sobolev embedding estimates, and interpolation estimates in $H^s(M)$, one gets
\[
\|\mathbf{f}\|_{L^2(M)}\le C\|\mathbf{f}\|_{L^2(M)}^{(2-4/s)\mu}
\]
with $0<\mu<1$ that can be chosen as close to $1$ as needed,
provided that $k\gg1$. It is enough now to choose $s>2$ and
$\mu$ so that $(2-4/s)\mu>1$ to deduce  that $\mathbf{f}=0$
if $\varepsilon\ll1$, see \eqref{nl_5}.
\end{proof}

\section{Rigidity of two-dimensional systems}\label{sec_7}
The main result of this section is the following rigidity theorem for two-dimensional systems:

\begin{Theorem}\label{2-dim}
If\/ $\dim M=2$ and $(g, \alpha)$ and $(g',\alpha')$
are simple magnetic systems on $M$
with the same boundary action function,
then these systems are gauge equivalent.
\end{Theorem}

This theorem generalizes the boundary rigidity theorem
for simple Riemannian surfaces which was established in \cite{PU}.
Our proof of Theorem \ref{2-dim} mimics that of the mentioned theorem in \cite{PU}.

First of all, by Theorem \ref{jet} and Lemma \ref{action-scat}
we may assume
that $g'$ coincides with $g$ on $\partial M$,
\begin{equation*}
g'|_{\partial M}=g|_{\partial M},
\end{equation*}
and that the scattering relations of the
magnetic systems $(M,g,\alpha)$ and $(M,g'\alpha')$ coincide:
\begin{equation}\label{scat}
\mathcal S'=\mathcal S.
\end{equation}

The crucial step then consists in establishing that the scattering relation
of a two-dimensional magnetic system $(M,g,\alpha)$ determines
the Dirichlet-to-Neumann (DN) map associated to the Laplace-Beltrami
operator of the Riemannian manifold $(M,g)$.

It is proved in \cite{LeU, LaU} for two-dimensional manifolds that the DN map determines the
conformal class of the Riemannian metric up to an isometry that is the identity
on the boundary.

Afterwards,  the proof of Theorem \ref{2-dim} is finished
by applying Theorem \ref{conf}
which claims magnetic boundary rigidity within a given conformal class.

Derivation of the connection between the scattering relation
and the DN map is based on
the properties of the magnetic ray transform
and the commutation formula between the magnetic flow and the fiberwise
Hilbert transform.

We proceed with describing the needed properties of the magnetic ray transform.

\subsection{More about the magnetic ray transform}
We recall that, given a notation $F$ for a function space
($C^k$, $L^p$, $H^k$, etc.),
we denote by
$\mathcal F(M)$
the corresponding space of pairs $[v,\varphi]$, with $v$ a $1$-form
and $\varphi$ a function on $M$.
Also, recall that in the space  $\mathcal L^2(M)$
we consider the norm \eqref{S7} defined as
\begin{equation*}
\|\mathbf{w}\|^2 = \int_M \left(|v|_g^2 + \varphi^2 \right)\,d\vol.
\end{equation*}

Associating each pair $[v,\varphi]\in \mathcal L^2(M)$ with the function
$\phi(x,\xi)=v_i(x)\xi^i+\varphi(x)$,
we may consider $\mathcal L^2(M)$
as a subspace of $L^2(SM)$.

Consider the restriction $\mathcal I$ of the magnetic ray transform  $I$
to $\mathcal L^2(M)$:
\begin{align*}
\mathcal I[v,\varphi](x,\xi)&=\int_0^{\ell(x,\xi)}
v_i(\gamma_{x,\xi}(t))\dot\gamma_{x,\xi}^i(t)\, dt
+\int_0^{\ell(x,\xi)}\varphi(\gamma_{x,\xi}(t))\, dt
\\
&=\mathcal I_1 v(x,\xi)+\mathcal I_0\varphi(x,\xi), \quad (x,\xi)\in\partial_+SM.
\end{align*}
By Lemma \ref{L2},  we receive an operator
\begin{equation*}
\mathcal I: \mathcal L^2(M)\to L^2_\mu(\partial_+SM).
\end{equation*}
Let
\begin{equation*}
\mathcal I^*: L^2_\mu(\partial_+SM)\to \mathcal L^2(M)
\end{equation*}
be its dual.
The same calculations as in Subsection 3.4 show that
\begin{equation*}
\mathcal I^*w=[\mathcal I^*_1w,\mathcal I^*_0w]
\end{equation*}
with
\begin{equation*}
\mathcal I^*_1w=\left(\int_{S_xM}\xi^i w^\sharp(x,\xi)\,d\sigma_x(\xi)\right),
\qquad
\mathcal I^*_0 w=\int_{S_xM}w^\sharp(x,\xi)\,d\sigma_x(\xi).
\end{equation*}

The following holds:
\begin{equation}\label{ker-i}
\ker \mathcal I=\mathcal P\mathcal L^2(M),
\end{equation}
\begin{equation}\label{im-i*}
\Image \mathcal I^*\subset \mathcal S\mathcal L^2(M).
\end{equation}
Equality \eqref{ker-i} is just the claim of Theorem \ref{linear}.
Equality \eqref{im-i*} follows, as soon as $\Image \mathcal I^*$ is in the orthogonal
complement to $\ker \mathcal I$ and  the orthogonal complement to
$\mathcal P\mathcal L^2(M)$ is $\mathcal S\mathcal L^2(M)$.

Consider the operator
\begin{equation*}
\mathcal N=\mathcal I^*\mathcal I: \mathcal L^2(M)\to \mathcal L^2(M),
\end{equation*}
which can be written down as
$$
\mathcal N[v,\varphi]=[\mathcal N_{11} v+\mathcal N_{10}\varphi,
\mathcal N_{01}v+\mathcal N_{00}\varphi].
$$

Considerations similar to those in the proof of \ref{pr_Spr1}
show the following:

\begin{Proposition}\label{pdo}
$\mathcal N$ is a $\Psi$DO in $M^\text{\rm int}$ of order
$-1$ with principal symbol
$$
\sigma_p(\mathcal N) = \diag\big(\sigma_p( \mathcal N_{11}),
\sigma_p(\mathcal N_{00}) \big),
$$
\begin{equation*}
\begin{split}
\sigma_p(\mathcal N_{11})^{j}_i(x,\xi) &= c_n|\xi|^{-1}(\delta^{j}_i-\xi^j\xi_i/|\xi|^2),\\
\sigma_p(\mathcal N_{00})(x,\xi) &=c_n|\xi|^{-1},
\end{split}
\end{equation*}
where $|\xi|^2=g^{ij}(x)\xi_i\xi_j$  and $\xi^j=g^{ji}(x)\xi_i$.
\end{Proposition}

Define
$$
C^\infty_\alpha(\partial_+SM)=\{w\in C^\infty(\partial_+SM)
: w^\sharp\in C^\infty(SM)\},
$$
where $w^\sharp$ is the function that is constant
along the orbits of the magnetic flow and equals $w$
on $\partial_+ SM$.  This space can be described in terms of the scattering
relation $\mathcal S$ alone; by the forthcoming Lemma \ref{c-alpha}
\begin{equation}\label{c-alpha-eq}
C^\infty_\alpha(\partial_+SM)=\{w\in C^\infty(\partial_+SM)
: Aw\in C^\infty(\partial(SM))\},
\end{equation}
where
$$
Aw(x,\xi)=\begin{cases} w(x,\xi),&(x,\xi)\in\partial_+SM,\\
                       w\circ\mathcal S^{-1}(x,\xi), &(x,\xi)\in\partial_-SM.\end{cases}
$$

The following is an analog of
\cite[Theorems 3.3.3, 3.3.4]{P}
(see also \cite[Theorem 1.4]{PU}, \cite[Theorems 4.1, 4.2]{PU2}).

\begin{Theorem}\label{surj}
Let\/ $(M,g,\alpha)$ be a simple magnetic system. Then,  for
every pair $[v,\varphi]\in \mathcal C^\infty(M)$, there exist
$w\in C^\infty_\alpha(\partial_+ SM)$
and $f\in C^\infty(M)$ such that
$$
[v+\nabla f,\varphi]=\mathcal I^*w.
$$
\end{Theorem}

(Note that if $v$ is solenoidal then,
in view of \eqref{im-i*}, $f$ is harmonic in~$M$.)

\begin{proof}
Our argument is the same as in \cite{P, PU, PU2}.
Embed $M$ into a closed manifold $\tilde M$ and extend $g$ and $\alpha$
smoothly to a Riemannian metric $g$ and a $1$-form $\alpha$ on $\tilde M$.

If $U\subset \tilde M$ is an open neighborhood of $M$
with smooth boundary,  then  $(\bar U,g,\alpha)$
is also a simple magnetic system
if $\partial U$ is close enough to $\partial M$.
Henceforth such an $U$ is assumed to be fixed.

Denote the magnetic ray transform for $(\bar U,g,\alpha)$ by
$\mathcal I_U$,  and denote
by $r_M$ the operator of restriction to $M$.
We have the following analog of
\cite[Theorems 3.3.1, 3.3.2]{P} (see also
\cite[Theorem 3.1]{PU}, \cite[Theorem 4.3]{PU2}):

\begin{Lemma}\label{inverse-n}
For every pair $[v,\varphi]\in \mathcal H^{s+1}(M)$, $s\ge0$,
there exists a~pair $[u,\psi]\in \mathcal H^s(U)$ and
$f\in H^{s+2}(U)$ such that
\[
[v+\nabla f,\varphi]=r_M \mathcal I_U^*\mathcal I_U[u,\psi].
\]
\end{Lemma}

Again, the proof is similar to the one in \cite{P, PU, PU2};
we will merely sketch it.

Cover $\tilde M$ by finitely many open sets $U_k$
such that $U=U_1$, $U_k\cap M=\emptyset$ for $k\ge2$,  and
$(\bar U_k,g, \alpha)$ is a simple magnetic system for every $k$.
Let $\{h_k\}$ be a~subordinate partition of unity such that $h_1|_M=1$.
Consider the operators $\mathcal I_{(k)}$,
$\mathcal I^*_{(k)}$ for $(\bar U_k,g, \alpha)$,
and define the following operator on the bundle
$T\tilde M\oplus(\tilde M\times\mathbb R)$:
$$
P[u,\psi]=\sum_k h_k \mathcal I^*_{(k)}\mathcal I_{(k)}([u,\psi]|_{U_k}).
$$
In view of Proposition \ref{pdo}, $P$ is a $\Psi$DO of order $-1$
with principal symbol
$$
c_n\diag\left(\frac{\delta^{j}_i}{|\xi|}-\frac{\xi^j\xi_i}{|\xi|^3},
\frac{1}{|\xi|}\right).
$$

The operator $\Lambda: C^\infty(T\tilde M)\to C^\infty(T\tilde M)$,
$$
\Lambda=-c_n\nabla(-\Delta)^{-3/2}\delta,
$$
is a $\Psi$DO of order $-1$ with principal symbol
$$
c_n\frac{\xi^j\xi_i}{|\xi|^3}.
$$
Therefore, the principal symbol of the operator
$$
P+\diag(\Lambda,0)
$$
equals
$$
c_n\diag\left(\delta^{j}_i|\xi|^{-1},|\xi|^{-1}\right),
$$
which means that $P+\diag(\Lambda,0)$ is an elliptic $\Psi$DO of order $-1$.

Now, the same arguments as in \cite{P, PU, PU2} show that
the operator
$$
r_M(P+\diag(\Lambda,0)): \mathcal H^s(U)\to\mathcal H^{s+1}(M)
$$
has closed range and finite codimension.

Since on $\mathcal H^s(U)$
\begin{equation}\label{rmp}
r_MP=r_M\mathcal I^*_U\mathcal I_U,
\end{equation}
\eqref{im-i*} and the argument of \cite{P, PU, PU2}
show that the equation
$$
r_M(P+\diag(\Lambda,0))[u,\psi]=[v,\varphi]
$$
has a solution $[u,\psi]\in \mathcal H^s(U)$
for every pair $[v,\varphi]\in \mathcal H^{s+1}(M)$.
Then by \eqref{rmp}
$$
[v,\varphi]+[\nabla f,0]=r_M\mathcal I^*_U\mathcal I_U[u,\psi]
$$
with $f=c_n(-\Delta)^{-3/2}\delta u$,  which proves the lemma.

Continuing the proof of Theorem \ref{surj},  observe that
by Lemma \ref{inverse-n}, for every pair
$[v,\varphi]\in\mathcal C^\infty(M)$
there exist $[u,\psi]\in C^{\infty}(U)$
and $f\in C^\infty(M)$
such that
\begin{equation}\label{rmi1}
[v+\nabla f,\varphi]=r_{M}\mathcal I^{\ast }_U \mathcal I_U[u,\psi].
\end{equation}

For $(x,\xi)\in SU$, define
$\ell^\pm(x,\xi)$ by
$$
\pi(\psi^{\ell^\pm}(x,\xi))\in\partial U,\quad \ell^-<0, \quad \ell^+>0.
$$

Define
$$
w^+(x,\xi)=\int_0^{\ell^+(x,\xi)}(u_i(\gamma_{x,\xi}(t)\dot\gamma_{x,\xi}^i(t)+\psi(\gamma_{x,\xi}(t)))\,dt,
$$
$$
w^-(x,\xi)=\int^0_{\ell^-(x,\xi)}(u_i(\gamma_{x,\xi}(t)\dot\gamma_{x,\xi}^i(t)+\psi(\gamma_{x,\xi}(t)))\,dt.
$$
These functions belong to $C^\infty(SU)$.
Let $w=(w^++w^-)|_{\partial _{+}SM}$. It is easy to see that
$(w^++w^-)$ is constant on the orbits of the magnetic flow;
therefore, $w^\sharp=(w^++w^-)|_{SM}$ and $w^\sharp\in C^\infty(SM)$.
Since $\mathcal I_U[u,\psi]=(w^++w^-)|_{\partial_+SU}$, we see from \eqref{rmi1}
that $\mathcal I^{\ast}w=[v+\nabla f,\varphi]$,
which completes the proof of the theorem.
\end{proof}

\subsection{Scattering relation and folds}
The main aim of this section is to prove the characterization
\eqref{c-alpha-eq} for the space $C^\infty_\alpha(\partial_+SM)$
(Lemma \ref{c-alpha}).
We will proceed in the same way as in \cite{PU}.

Preserving the notations of the previous section, define the map
$$
\Phi:\partial(SM)\rightarrow\partial_{-} SU
$$
by
\begin{equation}\label{phi}
\Phi(x,\xi)=\psi^{\ell^+(x,\xi)}(x,\xi),\quad
(x,\xi)\in\partial(SM).
\end{equation}
Since $\ell^+$ is smooth in $SU$,
$\Phi$ is smooth as well.

\begin{Lemma}[{cf. \cite[Theorem 4.1]{PU}}]\label{fold}
$\Phi$  is a fold map with fold
$S(\partial M)$.
\end{Lemma}

We recall that a smooth map $f:\mathcal M\to \mathcal N$ between two smooth manifolds
$\mathcal M$ and $\mathcal N$ of the same dimension is said to be
a {\em Whitney fold} (with fold L) at a point $m\in L$
if $f$ drops rank by one simply at $m$,
so that $\{x: df(x) \mbox { is singular}\}$ is
a smooth hypersurface near $m$ and
$\ker (df(m))$ is transverse to $T_{m}L$.

\smallskip
\noindent{\bf Example.}
Let $\Sigma$ be a manifold, $\mathcal R$ a smooth function on $\Sigma$
having $0$ as a regular value, and $\mathcal M=\{\mathcal R^{-1}(0)\}$.
Let $X$ be a nonzero vector field on $\Sigma$
such that $X\mathcal R(m)=0$, and $XX\mathcal R(m)\ne0$
for a point $m\in \mathcal M$.
Let $\mathcal N$ be a hypersurface in $\Sigma$ transversal to $X$
such that the projection $f:\mathcal M\to\mathcal N$
along the integral curves of $X$ is well defined.
Then $f$ is a Whitney fold at $m$ with fold
$\mathcal M\cap\{(X\mathcal R)^{-1}(0)\}$.

Indeed, we may assume  without loss of generality that $\Sigma$ is
a domain in $\mathbb R^n$, $X=\partial/\partial x^n$, and $\mathcal N$
is the hyperplane $\{x_n=0\}$. In this case,
we easily verify the claim by straightforward calculations.

\begin{proof}[Proof of Lemma \ref{fold}]
Let $\rho$ be a defining function of $M$ in $U$
with $\grad\rho(x)=\nu(x)$ for $x\in\partial M$. Then
$$
\mathbf G_\mu (\rho\circ\pi)(x,\xi)=\langle \xi,\grad\rho(x)\rangle,
$$
$$
\mathbf G_\mu^2(\rho\circ\pi)(x,\xi)=\langle \xi,\nabla_\xi \grad\rho(x)\rangle
+\langle Y_x(v),\grad\rho(x)\rangle.
$$
Therefore, for $(x,\xi)\in S(\partial M)$ we have
$\mathbf G_\mu (\rho\circ\pi)(x,\xi)=0$, while
$\mathbf G_\mu^2(\rho\circ\pi)(x,\xi)
=-\Lambda(x,\xi)+\langle Y_x(\xi),\nu(x)\rangle\ne0$
by strict magnetic convexity of $\partial M$.

We arrive directly at the above Example if we take
$\Sigma=S\tilde M$, $\mathcal R=\rho\circ\pi$, $\mathcal M=\partial(SM)$,
$\mathcal N=\partial(S\bar U)$, and $X=\mathbf G_\mu$.
This completes the proof of the lemma.
\end{proof}

Define the extension operator
$$
A: C(\partial_+SM)\to C(SM)
$$
by
$$
Aw(x,\xi)=\begin{cases} w(x,\xi),&(x,\xi)\in\partial_+SM,\\
                       w\circ\mathcal S^{-1}(x,\xi), &(x,\xi)\in\partial_-SM,\end{cases}
$$
where $\mathcal S$ is the scattering relation.

\begin{Lemma}[cf. {\cite[Lemma 1.1]{PU}}]\label{c-alpha}
If $(M,g,\alpha)$ is a simple magnetic system, then
$$
C^\infty_\alpha(\partial_+SM)=\{w\in C^\infty(\partial_+SM)
: Aw\in C^\infty(\partial(SM))\}.
$$
\end{Lemma}

\begin{proof}
If $w^\sharp\in C^\infty(SM)$, then $Aw=w^\sharp|_{\partial(SM)}$
is smooth as well.
Let us prove the converse. If $Aw\in C^\infty(\partial(SM))$,
then from Lemma \ref{fold} and \cite[Theorem C.4.4]{Hor}
we deduce the existence of a smooth function $v$ on a neighborhood of the range
$\Phi(\partial(SM))$ such that $w=v\circ\Phi$.

Let $\Psi:SM\to \partial_-SM$ denote the map
$$
\Psi(x,\xi)=\psi^{\ell(x,\xi)}(x,\xi), \quad (x,\xi)\in SM,
$$
and $\Psi_U:S\bar U\to \partial_-S\bar U$ the map
$$
\Psi_U(x,\xi)=\psi^{\ell^+(x,\xi)}(x,\xi),\quad (x,\xi)\in S\bar U.
$$
Note that $w^\sharp=\mathcal S^{-1}\circ \Psi$. Therefore,
$w^\sharp=v\circ\Phi\circ\mathcal S^{-1}\circ \Psi$.
It is easy to see that
$\Phi\circ\mathcal S^{-1}\circ \Psi=\Psi_U|_{SM}$. Since
$\Psi_U$ is smooth on $SM$, we conclude that $w^\sharp\in C^\infty(SM)$,
i.e., $w\in C^\infty_\alpha(\partial_+SM)$.
\end{proof}

\subsection{Hilbert transform}
From now on, we let $(M,g)$ be an oriented Riemannian surface.
Given $v\in T_xM$,
we will denote by $v_\perp$ the vector
obtained by rotating $v$ by $\pi/2$
according to the orientation of $M$.
In coordinates
$(v _{\perp })_{i}=\varepsilon _{ij}v ^{j},$\
\vspace{1pt}where
\[
\varepsilon =\sqrt{\det g}\begin{pmatrix}
  {0} & {1} \\
  {-1} & {0}
\end{pmatrix}.
\]

Consider the fiberwise Hilbert transform \cite[(1.4)]{PU}
\begin{equation*}
Hu(x,\xi)=\frac{1}{2\pi}\int_{S_{x}M}\frac{1+(\xi,\eta)}{(\xi_{\perp},\eta)}u(x,\eta)\,d\sigma_{x}(\eta),\quad
\xi\in S_{x}M.
\end{equation*}
If we fix $x\in M$ and a reference point $a\in S_xM$, any function on
the fiber $S_x$
can be treated as a function of an angular variable.
Then
\[
Hu(x,\theta)=\frac1{2\pi}\int_0^{2\pi}\cot\left(\frac{\varphi-\theta}{2}\right)u(x,\varphi)\,d\varphi.
\]

Define
$$
(\mathbf{G}_{\perp}u)(x,\xi)=\langle \xi_{\perp},\nabla^{|}u\rangle
=-\langle\xi,\nabla_{\perp}u\rangle,
$$
where $\nabla_{\perp}u=\varepsilon\nabla^{|} u$.
The following
commutation formula holds \cite[Theorem 1.5]{PU}:
\begin{equation}\label{hgu}
[H,\mathbf G]u=\mathbf{G}_{\perp}(u_{0})+(\mathbf{G}_{\perp}u)_{0},
\end{equation}
where
\[
u_{0}(x)=\frac{1}{2\pi}\int_{S_{x}M}u(x,\xi)\,d\sigma_{x}(\xi)
\]
is the average value on a fiber.

Let $V$ be the infinitesimal generator of the action of $S^1$
on the fibers of the canonical projection $\pi:SM\to M$. Then
the generator of the magnetic flow $\psi^t$ is given by
$$
\mathbf G_\mu=\mathbf G+\lambda V,
$$
where $\lambda$ is a function on $M$
such that $\Omega=\lambda\Omega_a$, with $\Omega_a$ the area form of $g$.

Since $H$ commutes with $\lambda V$, we get from \eqref{hgu}
\begin{equation}\label{hgmu}
[H,\mathbf G_\mu]u=\mathbf{G}_{\perp}(u_{0})+(\mathbf{G}_{\perp}u)_{0}.
\end{equation}

Substitute $u=w^\sharp$, $w\in C_\alpha^\infty(\partial_+SM)$,
into \eqref{hgmu}. This yields
\begin{equation}\label{ghw}
\mathbf G_\mu Hw^\sharp=-\mathbf G_\perp (w^\sharp)_0
-(\mathbf{G}_{\perp}w^\sharp)_{0}.
\end{equation}

Define the operator
$$
B: C(SM)\to C(\partial_+SM)
$$
by
$$
B u(x,\xi)=u(x,\xi)-u\circ \mathcal S(x,\xi), \quad (x,\xi)\in\partial_+SM.
$$

Clearly,
\begin{equation}\label{fundamental}
I\mathbf G_\mu u =-B u.
\end{equation}

Using \eqref{fundamental},
we deduce from \eqref{ghw} the identity
\begin{equation}\label{crux}
BHAw=I(\mathbf G_\perp (w^\sharp)_0)
+I((\mathbf{G}_{\perp}w^\sharp)_{0}),
\end{equation}
since $w^\sharp|_{\partial (SM)}=Aw$.

Note that
\begin{equation*}
\mathbf G_\perp(w^\sharp)_0
=-\frac 1{2\pi} \langle \xi, \nabla_{\perp}\mathcal I^*_0 w\rangle
\end{equation*}
and
\begin{equation*}
(\mathbf{G}_{\perp}w^\sharp)_{0}=-\frac 1{2\pi} \delta_{\perp}\mathcal I^*_1 w,
\end{equation*}
where
$$
\delta_{\perp}v=-\delta v_{\perp}.
$$

Hence,  we can rewrite \eqref{crux} as
\begin{equation}\label{main}
B HAw=-\frac 1{2\pi}\mathcal I[\nabla_{\perp}\mathcal I^*_0 w,\delta_{\perp}\mathcal I^*_1 w].
\end{equation}

\subsection{Dirichlet-to-Neumann map}
Given a  compact Riemannian manifold  $(M,g)$ with boundary, denote the
Laplace-Beltrami operator associated with $g$ by $\Delta_ g$.
Consider the Dirichlet
problem
$$ \Delta_ g u  = 0\hbox{ on }M,\quad u|_{\partial M}
=  f.$$
The DN map is defined by
$$
\Lambda_ g(f)=\langle\nu,\nabla u|_{\partial M}\rangle.
$$

In the two dimensional case, the DN map can also be described as follows.
Let $(h,h_{\ast})$ be a pair of conjugate harmonic functions on
$M$,
\[
\nabla h=\nabla_{\perp}h_{\ast},\quad \nabla
h_{\ast}=-\nabla_{\perp}h.
\]
Let $h^0$ and $h^0_*$ denote their traces on $\partial M$.
Then
\begin{equation}\label{dn-con}
\Lambda(h^0_*)=\langle\nu,\nabla h^0_*|_{\partial M}\rangle
=-\langle\nu,(\nabla_\perp h)|_{\partial M}\rangle
=\langle\nu_\perp,(\nabla h)|_{\partial M}\rangle
=\langle\nu_\perp,\nabla_{\partial M} h^0\rangle,
\end{equation}
where $\nabla_{\partial M}$ is the gradient
w.r.t.\ the induced metric on $\partial M$.

The following theorem is an analog of \cite[Theorem 1.3]{PU} and
states that the scattering relation
of a simple magnetic surface $(M,g,\alpha)$
completely determines the DN map of the metric $g$.

\begin{Theorem}\label{dn}
Let $(g,\alpha)$ and $(g',\alpha')$ be simple magnetic systems on a compact
surface $M$ with boundary such that
$g|_{\partial M}=g'|_{\partial M}$.
Assume that the scattering relations $\mathcal S$
and  $\mathcal S'$ of these systems coincide.
Then $\Lambda_{g}=\Lambda_{g'}$.
\end{Theorem}

\begin{proof}
Assume that  $h,h_{\ast}$ is a pair of smooth conjugate harmonic functions on
$M$. Then $\mathbf G_\perp h=\mathbf G h_*=\mathbf G_\mu h_*$.
By Theorem \ref{surj}, there are $w\in C_\alpha^\infty(\partial_+SM)$
and $f\in C^\infty(M)$
satisfying $\mathcal I^*_0 w=h$, $\mathcal I^*_1w=\nabla f$.
From \eqref{main} we then obtain
\begin{equation}\label{final}
BHAw=-\frac 1{2\pi} Bh^0_* ,
\end{equation}
since $\delta_{\perp}\nabla f=0$.
Hence, the following holds:

\begin{Lemma}[{cf. \cite[Theorem 1.6]{PU}}]\label{direct-dn}
If\/ $h,h_{\ast}$ is a pair of smooth conjugate harmonic functions on
$M$, then there is $w\in C^\infty_{\alpha}(\partial_{+}SM)$
such that $h=\mathcal I^*_0 w$ and equation \eqref{final} holds
with $h^0_{\ast}$ the trace of
$h_{\ast}$ on $\partial M$.
\end{Lemma}

In the opposite direction we have:

\begin{Lemma}[{cf. \cite[Theorem 1.6]{PU}}]\label{reverse-dn}
Suppose $h^0_{\ast}\in C^\infty(\partial M)$ and
$w\in C^\infty_{\alpha}(\partial_{+}SM)$ satisfy
equation \eqref{final}. Define $h:=\mathcal I^*_0 w$ and let
$h_{\ast}$ be the harmonic continuation of $h^0_{\ast}$ to $M$.
Then $h$ and $h_*$
are conjugate harmonic functions.
\end{Lemma}

\begin{proof}
Let $q$ be an arbitrary smooth extension of $h_{\ast}^{0}$
to  $M$. Note that
$$
\mathbf G_\mu q=\mathbf G q=\langle \xi, \nabla q\rangle.
$$
Using \eqref{main},
we can therefore rewrite \eqref{final} as
\[
-\mathcal I[\nabla_{\perp}h, \varphi]=\mathcal I[\nabla q,0],
\]
with $\varphi=\delta_{\perp}\mathcal I^*_1 w$.

Thus,
$$
\mathcal I[\nabla_{\perp}h+\nabla q, \varphi]=0.
$$

By \eqref{ker-i}
we then have $\varphi=0$ and $\nabla q+\nabla _{\perp }h=\nabla p$
for some smooth function $p$ on $M$
with $p|_{\partial M}=0$. Therefore,
$h$ and $q-p$ are conjugate harmonic functions.
Since $(q-p)|_{\partial M}=h^{0}_{\ast}$, we have
$q-p=h_*$, which implies that $h$ and $h_*$ are conjugate
harmonic functions.
\end{proof}

Continuing the proof of Theorem \ref{dn},
we have the following procedure to obtain the DN map
from the scattering relation. For an arbitrary given smooth
function $h_{\ast}^{0}$ on $\partial M$ we find a function
$w\in C_\alpha^{\infty}(\partial_{+}SM)$
that solves equation \eqref{final}.
Then the functions $h^{0}=2\pi (Aw)_{0}$ (notice that $2\pi
(Aw)_{0}=\mathcal I^{\ast}_0w|_{\partial M} $) and $h_{\ast}^{0}$ are
traces of conjugate harmonic functions. This gives the DN map
by means of \eqref{dn-con}.
\end{proof}

\subsection{Proof of Theorem \ref{2-dim}}
We end up with the proof of this theorem in the same way as in \cite{PU}.
If two simple magnetic systems $(g,\alpha)$ and $(g',\alpha')$
on a compact surface with boundary have the same boundary action functions,
then by Theorem \ref{jet} we may assume that
$g|_{\partial M}=g'|_{\partial M}$ and, by Lemma \ref{action-scat},
their scattering relations coincide.
Theorem \ref{dn} then tells us that the DN maps of the metrics
$g$ and $g'$  coincide.
Now, the result of \cite{LeU, LaU} implies the existence of a diffeomorphism
$f:M\to M$, which is the identity on $\partial M$, and of a function
$\omega$ such that $g'=\omega^2 f^*g$. Next, Theorem \ref{conf} yields
$\omega=1$ and $\alpha'=f^*\alpha+d\varphi$
for a smooth $\varphi$ vanishing
on $\partial M$. This concludes the proof of the theorem.

\appendix
\section{Geometry of magnetic systems}

\subsection{Ma\~n\'e's critical value and simplicity}
Here we adapt a certain part of the theory of convex superlinear Lagrangians
to the case of manifolds with boundary.

Let $M$ be a compact Riemannian
manifold with boundary and let $L:TM\to\mathbb R$
be a $C^\infty$ Lagrangian satisfying the following hypotheses:

\begin{itemize}
\item {\em Convexity}: For all $x\in M$ the restriction of $L$ to $T_xM$
has everywhere positive definite Hessian.
\item {\em Superlinear growth}:
$$
\lim_{|v|\to\infty}\frac{L(x,v)}{|v|}=+\infty
$$
uniformly on $x\in M$.
\end{itemize}

The action of $L$ on an absolutely continuous curve
$\gamma:[a,b]\to M$ is
$$
\mathbb A_L(\gamma)=\int_a^b L(\gamma(t),\dot\gamma(t))\,dt.
$$

For each $k\in \mathbb R$, the {\em Ma\~n\'e action potential}
$\mathbb A_k:M\times M\to\mathbb R\cup \{-\infty\}$ is defined by
$$
\mathbb A_k(x,y)=\inf_{\gamma\in \mathcal C(x,y)} \mathbb A_{L+k}(\gamma),
$$
where
$\mathcal C(x,y)=\{\gamma:[0,T]\to M: \gamma(0)=x,\ \gamma(T)=y,\ \gamma
\text{ is absolutely continuous}\}$.

The {\em critical level} $c=c(L)$ is defined as
\begin{align*}
c(L)&=\sup\{k\in\mathbb R: \mathbb A_{L+k}(\gamma)<0\text{ for some closed curve }\gamma\}
\\
&=\inf\{k\in\mathbb R: \mathbb A_{L+k}(\gamma)\ge 0\text{ for every closed
curve }\gamma\}.
\end{align*}

\begin{Proposition}
For $k<c(L)$, $\mathbb A_k(x,y)=-\infty$ for all $x,y\in M$.
For $k\ge c(L)$, $\mathbb A_k(x,y)\in\mathbb R$ for all $x,y\in M$.
\end{Proposition}

\begin{proof}
The same as in \cite[Proposition 2-1.1]{CI2}.
\end{proof}

One more characterization of $c(L)$ is useful.
Recall that the Hamiltonian $H:T^*M\to\mathbb R$ associated with $L$
is defined by the Fenchel transform
$$
H(x,p)=\sup\{p(v)-L(x,v):v\in T_xM\},
$$
and the supremum is achieved at $v$ such that $p=\de Lv(x,v)$.

\begin{Proposition}\label{kcl}
If there exists a $C^1$ function $f:M\to\mathbb R$ such that $H(df)<k$,
then $k\ge c(L)$.
\end{Proposition}

\begin{proof}
The same as in \cite[Lemma 5]{CIPP1}.
\end{proof}

Recall that the {\em energy function} $E:TM\to \mathbb R$ for $L$
is defined by
$$
E(x,v)=\de Lv(x,v)\cdot v-L(x,v),
$$
and that the energy function is constant on every solution $x(t)$ of the
Euler--Lagrange equation
\begin{equation}\label{e-l}
\frac d{dt}\de Lv(x(t),\dot x(t))=\de Lx(x(t),\dot x(t)).
\end{equation}

Let $\psi^t:TM\to TM$ be the Euler--Lagrange flow, defined
by $\psi^t(x,v)=(\gamma(t),\dot\gamma(t))$, where $\gamma$
is the solution of \eqref{e-l} with $\gamma(0)=x$ and $\dot\gamma(0)=v$.
For $x\in M$ and $k\in\mathbb R$, the {\em exponential map} at $x$
of energy $k$
is defined to be the partial map
$\exp_x:T_xM\to M$ given by
$$
\exp^k_x(tv)=\pi\circ\psi^t(v),\quad t\ge 0,\ v\in T_xM,\ E(x,v)=k.
$$
Then $\exp^k_x$  is a $C^1$-smooth partial map on $T_xM$
which is $C^\infty$-smooth on $T_xM\setminus\{0\}$.

The next proposition is similar to \cite[Theorem D]{CIPP2}
and has a similar proof.

\begin{Proposition}\label{kcl-simple}
If $\exp^k_x:(\exp^k_x)^{-1}(M)\to M$
is a diffeomorphism for every $x\in M$, then $k\ge c(L)$.
\end{Proposition}

\begin{proof}
Fix $q\in M$. Given $x\in M$, let $\gamma_{q,x}:[0,T_{q,x}]\to M$ be a solution
of the Euler--Lagrange equation with energy $k$, joining
$q$ to $x$.
Consider the function $f(x)=\mathbb A_k(\gamma_{q,x})$.
It is easy to see that the assumption of the lemma implies
this function is smooth in $M\setminus\{q\}$.
It follows from the first variation
formula of \cite[Lemma 4]{CIPP2} that for $x\in M\setminus\{q\}$
$$
d_xf(w)=\de Lv(x,\dot\gamma_{q,x}(T_{q,x}))\cdot w,
$$
which implies
$$
H(x,d_xf)=E(x,\dot\gamma_{q,x}(T_{q,x}))=k.
$$
This last equation also implies that $|d_xf|$ is uniformly bounded for all $x\in M\setminus\{q\}$.
Thus $f$ is Lipschitz in $M$ and a smoothing argument as in \cite{CIPP1, FM} shows that
for any $\varepsilon>0$ there exists $\tilde{f}\in C^{\infty}(M)$ for which
$H(x,d_x\tilde{f})< k+\varepsilon$ for all $x\in M$. Thus by Proposition \ref{kcl},
$k\geq c(L)$.
\end{proof}

The next proposition is an analog of \cite[Proposition 3-5.1]{CI2}
and has the same proof.

\begin{Proposition}\label{minimize-L}
If $k>c(L)$ and $x,y\in M$ $x\ne y$, then
there is $\gamma\in \mathcal C(x,y)$ such that
$$
\mathbb A_k(x,y)=\mathbb A_{L+k}(\gamma).
$$
Moreover, the energy of $\gamma$ is $E(\gamma,\dot\gamma)\equiv k$.
\end{Proposition}

Now, we apply the above to the case of magnetic systems.
For a simple magnetic system $(M,g,\alpha)$,
the magnetic flow can also obtained as the Euler--Lagrange flow with
the corresponding Lagrangian defined by
$$
L(x,v)=\frac 12|v|^2_g-\alpha_x(v).
$$

\begin{Lemma}\label{minimize}
Let $(g,\alpha)$ be a simple magnetic system on $M$.
For $x,y\in M$, $x\ne y$,
\begin{equation*}
\mathbb A_{1/2}(x,y)=\mathbb A_{L+1/2}(\gamma_{x,y})
=T_{x,y}-\int_{\gamma_{x,y}}\alpha,
\end{equation*}
where $\gamma_{x,y}:[0,T_{x,y}]\to M$ is the unit speed magnetic geodesic from $x$
to $y$.
\end{Lemma}

\begin{proof}
It is easy to see that the simplicity assumption implies that
for this Lagrangian the assumptions of Proposition \ref{kcl-simple}
hold for all $k$ sufficiently close to $1/2$. Therefore,
the proposition gives $1/2>c(L)$. Then Proposition \ref{minimize-L}
shows that, given $x\ne y$ in $M$, there is $\gamma\in\mathcal C(x,y)$
with energy $1/2$ (i.e., $\gamma$ is parametrized by arc length)
such that $\mathbb A(x,y)=\mathbb A(\gamma)$.
Using simplicity, one can then prove that $\gamma$ is a unit speed
magnetic geodesic, i.e., $\gamma=\gamma_{x,y}$.
\end{proof}

\subsection{Magnetic convexity}
Let $M$ be a compact manifold with boundary,
endowed with a Riemannian metric~$g$ and a closed $2$-form $\Omega$.
Consider a manifold $M_1$ such that $M_1^\text{int}\supset M$.
Extend $g$ and $\Omega$ to $M_1$ smoothly, preserving the former notation
for extensions.
We say that $M$ is
{\em magnetic convex} at $x\in \partial M$
if there is a neighborhood $U$ of $x$ in $M_1$ such that
all  unit speed magnetic geodesics in $U$, passing through $x$
and tangent to $\partial M$ at $x$,
lie in $M_1\setminus M^{\text{int}}$.
If, in addition, these geodesics
do not intersect $M$ except for $x$,
we say that $M$ is {\em strictly magnetic convex} at $x$.
It is not hard to show that these definitions depend neither on the
choice of $M_1$ nor on the way we extend $g$ and $\Omega$ to $M_1$.

As before, we let $\Lambda$ denote the second fundamental form of $\partial M$
and $\nu(x)$ the inward unit normal to $\partial M$ at $x$.

\begin{Lemma}\label{convexity}
If $M$ is magnetic convex at $x\in\partial M$, then
\begin{equation}\label{nonstrict-conv}
\Lambda(x,v)\ge\langle Y_x(v),\nu(x)\rangle
\quad\text{for all } v\in S_x(\partial M).
\end{equation}
If the inequality is strict, then $M$ is strictly magnetic
convex at $x$.
\end{Lemma}

\begin{proof}
Suppose $M$ is convex at $x$.
Choosing a smaller $U$ if necessary, we may assume
that there is a smooth function $\rho$ on $U$ such that $|\grad\rho|=1$
and $\partial M\cap U=\rho^{-1}(0)$. Further we may assume that all the above
geodesics lie in $U^-=\{x:\rho(x)\le 0\}$.

Let $v\in S_x (\partial M)$ and $\gamma(t)$ be the
magnetic geodesic with $\gamma(0)=x$, $\dot \gamma(0)=v$. By assumption,
$\rho\circ \gamma(t)\le 0$ for all small $t$. Therefore,
$$
\frac{d^2}{dt^2}\big[\rho\circ\gamma(t)\big]\Big|_{t=0}\le 0.
$$
Since
\begin{multline*}
\frac{d^2}{dt^2}\big[\rho\circ\gamma(t)\big]
=\frac{d}{dt}\langle \grad\rho(\gamma(t)),\dot\gamma(t)\rangle
\\
=\langle \nabla_{\dot\gamma(t)}\grad\rho(\gamma(t)),\dot\gamma(t)\rangle
+\langle \grad\rho(\gamma(t)),\ddot\gamma(t)\rangle
\\
=\hess_{\gamma(t)}\rho (\dot\gamma(t),\dot\gamma(t))+\langle \grad\rho(\gamma(t)),Y(\dot\gamma(t))\rangle
\end{multline*}
and since $\Lambda(x,\xi)=-\hess_x\rho(v,v)$ and $\grad\rho(x)=\nu(x)$
when $(x,v)\in\partial(SM)$, we arrive at \eqref{nonstrict-conv}.

Now, assume that \eqref{nonstrict-conv} is strict.
Then there is $\delta>0$ such that for every magnetic geodesic $\gamma$
in $N$ with $\gamma(0)=x$ and $\dot\gamma(0)=v\in S_x(\partial M$),
$$
\frac{d^2}{dt^2}\big[\rho\circ\gamma(t)\big]\Big|_{t=0}\le -\delta.
$$
Thus, there is a small $\varepsilon>0$ such that
$$
\rho\circ\gamma(t)\le -\frac 14\delta t^2
\quad \text{for all }t\in (-\varepsilon,\varepsilon).
$$
This implies the second statement.
\end{proof}

\subsection{Exponential map}

\begin{Lemma} \label{exp-map}
The map $\exp_{x}^{\mu}:T_{x}M\to M$ is\/ $C^1$
and $C^{\infty}$ on $T_{x}M\setminus\{0\}$.
This map is $C^2$ if and only if $\Omega=0$.
\end{Lemma}

\begin{proof}
Recall that $\exp_{x}^\mu(tv)=\pi\circ \psi^{t}(x,v):=\gamma(t,x,v)$
where $v$ has norm one and $t\geq 0$. Clearly this implies that
$\exp_{x}^{\mu}$ is $C^{\infty}$ on $T_{x}M\setminus\{0\}$.
Since
\[\frac{d}{dt}\bigg|_{t=0}\pi\circ\psi^{t}(x,v)=v\]
we see that the directional derivative of $\exp_{x}^\mu$
at $0$ in the direction of $v$ exists and equals $v$.
In a coordinate system around $x$ write
\begin{equation}\label{taut}
\exp_{x}^{\mu,i}(tv)=\gamma^{i}(t,x,v)
\end{equation}
where $\gamma^i$ depends smoothly in $(t,x,v)$. Differentiating with respect
to $v$ we obtain
\[\frac{\partial \exp_{x}^{\mu,i}}{\partial v^j}(tv)\,t=
\frac{\partial \gamma^{i}}{\partial v^j}(t,x,v).\]
Since $\gamma^{i}(0,x,v)=x^i$,
$\frac{\partial\gamma^{i}}{\partial v^{j}}(0,x,v)=0$.
Therefore,
\begin{align*}
\lim_{t\to 0^{+}}\frac{\partial \exp_{x}^{\mu,i}}{\partial v^j}(tv)
&=\lim_{t\to 0^{+}}\frac 1t\frac{\partial \gamma^{i}}{\partial v^j}(t,x,v)
\\
&=\frac{\partial^{2}\gamma^{i}}{\partial t\partial v^{j}}(0,x,v)
=\frac{\partial^{2}\gamma^{i}}{\partial v^j\partial t}(0,x,v)
=\delta_{j}^i
\end{align*}
since $\dot{\gamma}^i(0,x,v)=v^i$. Since $S_{x}M$ is compact, the above limit
is uniform in $v\in S_{x}M$ and thus $\exp_{x}^{\mu}$ has continuous partial
derivatives at $0$, i.e., $\exp_{x}^{\mu}$ is $C^1$.

Suppose now that $\exp_{x}^\mu$ is $C^2$. Differentiate (\ref{taut}) twice
with respect to $t$ ($v\in S_{x}M$)
\[v^{j}v^k \frac{\partial^2 \exp_{x}^{\mu,i}}{\partial v^j v^k}(tv)
=\ddot{\gamma}^i(t,x,v).\]
Using the equations of a magnetic geodesic
\[\ddot{\gamma}^i+\dot{\gamma}^j\dot{\gamma}^k\Gamma_{jk}^i
=Y_{k}^i\dot{\gamma}^k\]
we obtain
\[v^{j}v^k \frac{\partial^2 \exp_{x}^{\mu,i}}{\partial v^j v^k}(tv)
=-\dot{\gamma}^j\dot{\gamma}^k\Gamma_{jk}^i+Y_{k}^i\dot{\gamma}^k.\]
Let $t\to 0^{+}$. Then
\[v^{j}v^k\left( \frac{\partial^2 \exp_{x}^{\mu,i}}{\partial v^j v^k}(0)
+\Gamma_{jk}^i(x)\right)=Y_{k}^i(x) v^k.\]
Since this holds for all $v\in S_{x}M$ we must have $Y_{k}^i(x)=0$,
i.e., $\Omega=0$.
\end{proof}

\subsection{Santal\'o's formula}\label{santalo-s}
If $(M,g)$ is a compact Riemannian manifold with boundary,
we endow its unit sphere bundle $SM$
with its usual Liouville (local product) measure
$d\Sigma^{2n-1}$, and
endow the bundle $\partial_+SM$
with its standard measure $d\Sigma^{2n-2}$ (again a local product measure
where the measure on the fiber is the measure of a hemisphere).
Denote by $d\mu$ the measure on $\partial_+SM$ given by
$$
d\mu(x,\xi)=\langle \xi,\nu(x)\rangle\,d\Sigma^{2n-2}(x,\xi),
$$
where $\nu(x)$ is the inward unit normal of $\partial M$
at a point $x$.

The following version of Santal\'o's formula holds for magnetic flows.

\begin{Lemma}
Suppose that $(M,g,\alpha)$ is simple.
Then for every continuous function
$\varphi:SM\to\mathbb R$ we have
\begin{equation}\label{santalo-f}
\int_{SM}\varphi\,d\Sigma^{2n-1}
=\int_{\partial_+SM}d\mu(x,\xi)
\int_0^{\ell(x,\xi)}\!\!\varphi(\gamma_{x,\xi}(t),\dot\gamma_{x,\xi}(t))\,dt.
\end{equation}
\end{Lemma}

\begin{proof}
The argument we use is the same as in \cite{Sh2}.
We give it for the sake of completeness.

First, we recall the well-known fact that the Liouville measure
is invariant under the magnetic flow (for example, because
$\omega_0{}^{\wedge n}=(\omega_0+\pi^*\Omega)^{\wedge n}$
while $\omega_0+\pi^*\Omega$ is flow invariant). Now, let
$D=\{(x,\xi;t)\in\partial_+SM:0\le t\le \ell(x,\xi)\}$, and define
$\Psi:D\to SM$ by $\Psi (x,\xi;t)=\psi^t(x,\xi)$,
where $\psi^t$ is the magnetic flow. Then
\begin{equation}\label{pullback}
\int_{SM}\varphi\,d\Sigma^{2n-1}
=\int_D(\varphi\circ\Psi)\,\Psi^*(d\Sigma^{2n-1}).
\end{equation}

By construction, $\Psi$ conjugates $\psi^t$ with the flow generated
by $\partial/\partial t$ on $D$. Since $d\Sigma^{2n-1}$
is invariant under $\psi^t$, $\Psi^*(d\Sigma^{2n-1})$ is invariant
under the flow of $\partial/\partial t$. Then
$$
\Psi^*(d\Sigma^{2n-1})(x,\xi;t)=a(x,\xi)d\Sigma^{2n-2}\wedge dt
$$
for some function $a$ on $\partial_+SM$, so that \eqref{pullback}
takes the form
$$
\int_{SM}\varphi\,d\Sigma^{2n-1}
=\int_{\partial_+SM}a(x,\xi)\,d\Sigma^{2n-2}(x,\xi)
\int_0^{\ell(x,\xi)}\varphi(\gamma_{x,\xi}(t),\dot\gamma_{x,\xi}(t))\,dt.
$$

We are left with proving that $a(x,\xi)=\langle \xi,\nu(x)\rangle$.
To this end, it suffices to show that
$$
\Psi^*(d\Sigma^{2n-1})(x,\xi;0)=\langle \xi,\nu(x)\rangle\, d\Sigma^{2n-2}\wedge dt.
$$

Note that $d\Sigma^{2n-1}=-d\Sigma^{2n-2}\wedge dr$ on $\partial(SM)$,
where $r(x)=\dist(x,\partial M)$ is the distance function from $x$
to $\partial M$. Since the differential of $\Psi$ is the identity on
$T_{(x,\xi)}(\partial_+SM)$ and takes $\partial /\partial t$ to
the generator $\mathbf G_\mu$ of the flow $\psi^t$, we have
$$
\Psi^*(d\Sigma^{2n-1})(x,\xi;0)=(\mathbf G_\mu r)\,d\Sigma^{2n-2}(x,\xi)\wedge dt.
$$
As soon as $\mathbf G_\mu r=\langle \xi,\nabla r\rangle=-\langle \xi,\nu(x)\rangle$,
we are done.
\end{proof}

\subsection{Index form of a magnetic geodesic}
Let $(M,g,\alpha)$ be a simple magnetic system.
For every $x\in M$, $\exp_{x}^{\mu}:T_{x}M\to M$
is a diffeomorphism restricted to
a suitable set in $T_{x}M$ which is diffeomorphic to a closed ball.

Let $\pi:SM\to M$ be the canonical projection and let for $v\in SM$,
\[V(v):=\mbox{\rm ker}\,d_{v}\pi,\]
which is an $(n-1)$-dimensional subspace of $T_{v}SM$, and
\begin{equation}\label{def-e}
E(v):=V(v)\oplus\mathbb R\mathbf G_{\mu}(v).
\end{equation}

\begin{Lemma} If $\gamma:[0,T]\to M$ is a unit speed magnetic geodesic,
then
\[d_{\dot{\gamma}(0)}\psi^{t}(E)\cap V(\dot{\gamma}(t))=\{0\}\]
for every $t\in (0,T]$.
\label{trans}
\end{Lemma}

\begin{proof} Take $v\in SM$ and $t\in (0,T]$.
From the definition of $\exp_{x}^\mu$ one sees right away that
\[\mbox{\rm image}(d_{tv}\exp_{x}^{\mu})
=d_{\dot{\gamma}(t)}\pi(d_{\dot{\gamma}(0)}\psi^{t}(E)).\]
Since $d_{w}\exp_{x}^{\mu}$ is a linear isomorphism for every $w\in T_{x}M$
at which $\exp_{x}^{\mu}$ is defined, the lemma follows.
\end{proof}

Given a unit speed magnetic geodesic $\gamma:[0,T]\to M$,
let ${\mathcal A}$ and ${\mathcal C}$ be the operators
on smooth vector fields along $\gamma$ defined by
\begin{equation*}
{\mathcal A}(Z)=\ddot{Z}+R(\dot{\gamma},Z)\dot{\gamma}
-Y(\dot{Z})-(\nabla_{Z}Y)(\dot{\gamma}),
\end{equation*}
\begin{equation*}
{\mathcal C}(Z)=R(\dot{\gamma},Z)\dot{\gamma}- Y(\dot{Z})
-(\nabla_{Z}Y)(\dot{\gamma}).
\end{equation*}

A vector field $J$ along $\gamma$ is said to be a {\em magnetic Jacobi field}
if it satisfies the equation
\begin{equation}
{\mathcal A}(J)=0.  \label{jacobi}
\end{equation}

Let $\Lambda$ denote the $\mathbb R$-vector space of smooth vector fields
$Z$ along~$\gamma$ such that $Z(0)=Z(T)=0$.
Define the quadratic form
$\ind:\Lambda\to \mathbb R$ by
\begin{equation*}
\ind(Z,Z)=\int_0^T\left\{|\dot Z|^2
-\langle {\mathcal C}(Z),Z\rangle
-\langle Y(\dot\gamma),Z\rangle^2\right\}\,dt.
\end{equation*}
Note that
\begin{equation*}
\ind(Z,Z)=-\int_{0}^{T}\{\left\langle {\mathcal A}(Z),Z\right\rangle
+\langle Y(\dot{\gamma}),Z\rangle^{2}\}\,dt.  \label{G}
\end{equation*}

Clearly, $\ind(Z,Z)$ generalizes the index form of a geodesic
in a Riemannian manifold. The next lemma is an analog of
a well-known assertion from Riemannian geometry.

\begin{Lemma}[Index Lemma]\label{index-lemma}
If $Z\in\Lambda $ is orthogonal to $\dot{\gamma}$, then
\[\ind(Z,Z)\geq 0,\]
with equality if and only if $Z$ vanishes.
\end{Lemma}

\begin{proof}
Note that the subbundle $E$ defined by \eqref{def-e} is Lagrangian.

If $\xi\in E(v)$, then $J_{\xi}(t)=d\pi\circ d\psi^{t}(\xi)$
satisfies the Jacobi equation (\ref{jacobi}).
Since
\[\left. d_{\dot{\gamma}(t)}\pi\right|_{E(\dot{\gamma}(t))}:
E(\dot{\gamma}(t))\rightarrow T_{\gamma(t)}M\]
is an isomorphism for all $t\in (0,T]$,
there exists a basis $\{\xi_{1},\dots,\xi_{n}\}$ for $E(v)$
such that $\{J_{\xi_{1}}(t),\dots,J_{\xi_{n}}(t)\}$ is a basis
of $T_{\gamma(t)}M$ for all $t\in (0,T]$.
Without loss of generality we may assume that
$\xi_{1}=\mathbf G_{\mu}(v)$, $J_{\xi_{1}}=\dot{\gamma}$
and $J_{\xi_{i}}(0)=0$ for $i\geq 2$.

Let us set for brevity $J_{i}=J_{\xi_{i}}$.
Then if $Z$ is an element of $\Lambda$, we can write
for $t\in (0,T]$
\[ Z(t)=\sum_{i=1}^{n}f_{i}(t)J_{i}(t)\]
for some smooth functions $f_{1},\dots,f_{n}$.
The functions $f_{i}$ can in fact be smoothly
extended to $t=0$. Indeed, for $i\geq 2$, we can write $J_{i}(t)=t\,A_{i}(t)$
where $A_{i}$ is a smooth vector field such that $A_{i}(0)=\dot{J}_{i}(0)$.
Since $\{\dot{\gamma}(t),A_{2}(t),\dots,A_{n}(t)\}$ is
now a basis for all $t\in [0,T]$, there
exist smooth functions $g_{i}$ such that for all $t\in [0,T]$
\[ Z(t)=g_{1}(t)\dot{\gamma}(t) +\sum_{i=2}^{n}g_{i}(t)A_{i}(t).\]
Therefore for $t\in (0,T]$, $g_{1}(t)=f_{1}(t)$ and for $i\geq 2$,
$g_{i}(t)=t\,f_{i}(t)$.
Since $Z(0)=0$, $g_{i}(0)=0$ for all $i$
and the $f_{i}$'s smoothly extend to $t=0$.

Now we can write
\begin{equation}
\ind(Z,Z)=-\sum_{i,j}\int_{0}^{T}\left\langle{\mathcal A}(f_{i}J_{i}),f_{j}J_{j}\right\rangle\,dt
-\int_{0}^{T}\langle Y(\dot{\gamma}),Z\rangle^{2}\,dt.      \label{bilinear}
\end{equation}
An easy computation shows that
\[{\mathcal A}(f_{i}J_{i})
=\ddot{f}_{i}J_{i}+2\dot{f}_{i}\dot{J}_{i}
-\dot{f}_{i}Y(J_{i})+f_{i}{\mathcal A}(J_{i}).\]
Since $J_{i}$ satisfies equation (\ref{jacobi}),
then ${\mathcal A}(J_{i})=0$ and hence,
\[\left\langle{\mathcal A}(f_{i}J_{i}),J_{j}\right\rangle
=\ddot{f}_{i}\left\langle J_{i},J_{j}\right\rangle+
2\dot{f}_{i}\langle\dot{J}_{i},J_{j}\rangle
-\dot{f}_{i}\langle Y(J_{i}),J_{j}\rangle.\]
Observe that since $E$ is a Lagrangian subspace,
\[\langle J_{i},\dot{J}_{j}\rangle
-\langle\dot{J}_{i},J_{j}\rangle+\langle Y(J_{i}),J_{j}\rangle=0,\]
and then
\[\left\langle{\mathcal A}(f_{i}J_{i}),J_{j}\right\rangle
=\frac{d}{dt}(\dot{f}_{i}\left\langle J_{i},J_{j}\right\rangle).\]
Now we can write
\[\int_{0}^{T}\left\langle{\mathcal A}(f_{i}J_{i}),f_{j}J_{j}\right\rangle\,dt
=\left.\langle\dot{f}_{i}J_{i},f_{j}J_{j}\rangle\right|_{0}^{T}
-\int_{0}^{T}\langle\dot{f}_{i}J_{i},\dot{f}_{j}J_{j}\rangle\,dt.\]
Combining the last equality with (\ref{bilinear}) we obtain
\[\ind(Z,Z)=\int_{0}^{T}\Big\|\sum_{i=1}^{n}\dot{f}_{i}J_{i}\Big\| ^{2}\,dt
-\Big\langle\sum_{i=1}^{n}\dot{f}_{i}J_{i},Z\Big\rangle\Big|_{0}^{T}-
\int_{0}^{T}\langle Y(\dot{\gamma}),Z\rangle^{2}\,dt.\]
But $Z(0)=Z(T)=0$, therefore
\begin{equation}
I(V,V)=\int_{0}^{T}\Big\| \sum_{i=1}^{n}\dot{f}_{i}J_{i}\Big\| ^{2}\,dt
-\int_{0}^{T}\langle Y(\dot{\gamma}),Z\rangle^{2}\,dt.
\label{finvi}
\end{equation}
Now let
\[W:=\sum_{i=2}^{n}\dot{f}_{i}J_{i}.\]
Since $J_1=\dot{\gamma}$ we have:
\[\Big\langle \sum_{i=1}^{n}\dot{f}_{i}J_{i},
\sum_{i=1}^{n}\dot{f}_{i}J_{i}\Big\rangle=
\langle \dot{f}_{1}\dot{\gamma}+W, \dot{f}_{1}\dot{\gamma}+W\rangle
=\dot{f}_{1}^{2}+2\dot{f}_{1}\langle \dot{\gamma},W\rangle+\langle W,W\rangle.\]
Differentiating $\langle Z,\dot{\gamma}\rangle=0$ we get
\[\langle \dot{Z},\dot{\gamma}\rangle+\langle Z,Y(\dot{\gamma})\rangle=0.\]
But
\[\langle \dot{Z},\dot{\gamma}\rangle
=\Big\langle \sum_{i=1}^{n}\dot{f}_{i}J_{i},\dot{\gamma}\Big\rangle
=\dot{f}_{1}+\langle W,\dot{\gamma}\rangle\]
since $\langle \dot{J}_{i},\dot{\gamma}\rangle=0$ for all $i$. Therefore
\[\langle Y(\dot{\gamma}),Z\rangle^{2}=\dot{f}_{1}^{2}+2\dot{f}_{1}\langle W,\dot{\gamma}\rangle+
\langle W,\dot{\gamma}\rangle^{2}.\]
Thus
\[\Big\langle \sum_{i=1}^{n}\dot{f}_{i}J_{i},\sum_{i=1}^{n}\dot{f}_{i}J_{i}\Big\rangle
-\langle Y(\dot{\gamma}),Z\rangle^{2}=\langle W,W\rangle-\langle W,\dot{\gamma}\rangle^{2}.\]
If we let $W^{\perp}$ be the orthogonal projection of $W$
to $\dot{\gamma}^{\perp}$, the last equation and (\ref{finvi}) give:
\[\ind(Z,Z)=\int_{0}^{T}\| W^{\perp}\|^{2}\,dt\geq 0\]
with equality if and only if $W^{\perp}$ vanishes identically.
But if $W^{\perp}$ vanishes, then
\[-\langle W,\dot{\gamma}\rangle\dot{\gamma}+\sum_{i=2}^{n}\dot{f}_{i}J_{i}=0\]
which implies that the functions $f_{i}$ are constant for $i\geq 2$. Thus
$Z$ is of the form $f_{1}\dot{\gamma}+ J$ where $J$ is a magnetic Jacobi field.
But $Z(T)=0$ implies $J(T)=0$. Since the $J_{i}'s$ are linearly independent at $T$,
$J$ must vanish identically and since $Z$ is orthogonal to $\dot{\gamma}$, $Z$ must also vanish.
\end{proof}

\section{Study of a certain class of integral operators
with singular kernels}\label{singular}

As we mentioned before, the fact that the magnetic exponential map is smooth in polar coordinates only forces us to work in polar coordinates as well. In this appendix we study a  class of operators that naturally arise in our analysis.

Let $U\subset \mathbb R^n$ be open and $g$ be a smooth Riemannian metric in
a neighborhood of~$\bar U$.

\begin{Lemma}   \label{lemma_SA}
Let $A :C_0(U)\to C(U)$ be the operator
\begin{equation}  \label{S19a}
\mathcal Af(x) = \int_{S_xU}\int_{\mathbb{R}}  A(x,r,\omega) f(x+r\omega) \,dr\,d\sigma_x(\omega),
\end{equation}
with $A\in C^\infty(U\times \mathbb{R}\times S_xU)$.
Then $\mathcal A$ is a classical $\Psi$DO of order $-1$ with full symbol
\[
a(x,\xi) \sim \sum_0^\infty a_k(x,\xi),
\]
where
\[
a_k(x,\xi) = 2\pi \frac{i^k}{k!}\int_{S_xU} \partial_r^k A(x,0,\omega)\delta^{(k)}(\omega\cdot \xi) \, d\sigma_x(\omega).
\]
\end{Lemma}

\begin{proof}
Notice first that if $A$ is an odd function of $(r,\omega)$, then $\mathcal Af=0$.
Therefore, we can replace $A$ above by
$A_\text{even}(r,\omega) = (A(r,\omega) + A(-r,-\omega))/2$.
Next, it is easy to check that we can integrate over $r\ge0$ only and
double the result. Therefore,
\begin{equation}  \label{S19b}
\mathcal Af(x) = 2\int_{S_xM}\int_0^\infty  A_\text{even}(x,r,\omega) f(x+r\omega) \,dr\,d\sigma_x(\omega),
\end{equation}
Consider now $r$, $\omega$ as polar coordinates for $z=r\omega$,
and make also the change of variables $y=x+z$ to get
\begin{equation}  \label{S20}
\mathcal Af(x) = 2 (\det g(x))^{1/2}\int A_\text{even}\bigg(x,|y-x|_g  ,\frac{y-x}{|y-x|_g} \bigg) \frac{f(y)}{  |y-x|_g^{n-1}}\,dy,
\end{equation}
where the subscript $g$ refers to $g(x)$. Let
\begin{equation}  \label{S20a}
A_\text{even}(x,r,\omega) = \sum_{k=0}^{N-1} A_{\text{even},k} (x,\omega) r^k+r^{N}R_N(x,r,\omega)
\end{equation}
be a finite Taylor expansion of $A_\text{even}$ in $r$ near $r=0$ with $N>0$.
It follows easily that
$2A_{\text{even},k}(x,\omega) = A_k(x,\omega) +(-1)^k A_k(x,-\omega)$,
where $k!A_k = \partial^k_r|_{r=0}A$, and in particular,
$A_{\text{even},k}(x,\omega)r^k$ is even w.r.t.\ $(r,\omega)$.
The remainder term contributes to (\ref{S20}) an operator that
maps $L^2_\text{comp}(U)$ into $H^{N-N_0}(U)$   with some fixed $N_0$.
To study the contribution of the other terms, write
\begin{equation}  \label{S20b}
\mathcal A_{\text{even},k} f(x) = 2(\det g(x))^{1/2} \int A_{\text{even},k}\bigg(x,  \frac{y-x}{|y-x|_g} \bigg) |y-x|_g^{k-n+1} f(y) \, dy.
\end{equation}
The kernel of $\mathcal A_{\text{even},k}$ is therefore a function of $x$ and $z=y-x$, with a polynomial
singularity at $y-x=0$, and it is therefore a formal $\Psi$DO with symbol
that can be obtained by taking Fourier transform in the $z$ variable.
Motivated by this, apply the Plancherel theorem to the integral above to get
\[
\mathcal A_{\text{even},k}f(x) = (2\pi)^{-n} \int e^{ix\cdot \xi} a_k(x,\xi) \hat f(\xi)\, d\xi,
\]
where
\begin{align}  \nonumber
a_k(x,\xi) &= 2\int e^{-iy\cdot \xi} A_{\text{even},k}\bigg(x, \frac{y-x}{|y-x|_g} \bigg) |y-x|_g^{k-n+1} (\det g(x))^{1/2}  \,dy\\ \nonumber
           &= 2\int_{S_xU}\int_0^\infty e^{-ir\omega\cdot \xi} A_{\text{even},k} (x,\omega) r^k\, dr\, d\sigma_x(\omega)\\ \nonumber
           &= \int_{S_xU}\int_{-\infty}^\infty e^{-ir\omega\cdot \xi} A_k (x,\omega) r^k\, dr\, d\sigma_x(\omega)\\
           &= 2\pi i^k\int_{S_xU}  A_k(x,\omega)\delta^{(k)}(\omega\cdot \xi) \, d\sigma_x(\omega).     \label{S21}
\end{align}
In the third line, we used the fact that
$A_{\text{even},k}(x,\omega) r^k$ is even.
Note that $a_k(x,\xi)$ is homogeneous in $\xi$ of order $-k-1$
and smooth away from $\xi=0$ but a distribution (in $\mathcal{S}'$) near zero.
To deal with this, choose $\chi\in C_0^\infty$  supported in $|\xi|\le1$
and equal to $1$ near $\xi=0$.
Write $a(x,\xi) = \chi(\xi)a(x,\xi) + (1-\chi(\xi))a(x,\xi)$.
The second term is a classical amplitude,
while the first one contributes the term
\begin{equation}  \label{S21a}
\mathcal A_{\text{even},k} (\check{\chi}*f )
\end{equation}
to \eqref{S20b} that is smooth,  as can be easily seen by making
the change of variables $z=y-x$ in \eqref{S20b}, see also \eqref{S44a}.
\end{proof}

\begin{Remark}
If $A(x,r,\omega)$ and $g$ are smooth of class $C^k$ only,
then $\mathcal A$ is an $\Psi$DO with an amplitude of finite smoothness $l(k)$,
admitting a finite expansion. If $k\gg1$, then $l\gg1$,
and one can still construct a finite order parametrix of an elliptic $\Psi$DO
in this class and the usual $H^{s_1}\to H^{s_2}$ estimates still hold,
if $k\gg1$, depending on $s_1$, $s_2$.
This has been used  already in \cite{SU, SU1}.
\end{Remark}

We return to the analysis of the singular operator $\mathcal A$
introduced in Lemma~\ref{lemma_SA} under the assumption that $A$ and $g$ are analytic. Our reference for the calculus of analytic $\Psi$DOs is \cite{T}.

\begin{Lemma}   \label{lemma_SA_an}
Let $\mathcal A :C_0(U)\to C(U)$ be the operator (\ref{S19a})
with $A(x,r,\omega)$ analytic for $(x,\omega)\in U\times S_xU$,
and $r\in\mathbb{R}$ such that $x+r\omega\in U$.
Then $A$ is an analytic $\Psi$DO of order $-1$ with a symbol expansion as in Lemma~\ref{lemma_SA}.
\end{Lemma}

\begin{proof}

Notice first that $\mathcal A$ is analytic pseudolocal, see \cite[Theorem~V.2.1]{T}.

By performing the change of variables $\omega' = g^{1/2}(x)\omega$, we reduce the lemma to the case where $g$ is the Euclidean metric. Let $U'\subset\subset U$.

Let us estimate $a_k(x,\xi)$, see Lemma~\ref{lemma_SA}.
Since $A(x,r,\omega)$ is analytic, and $a_k$ is homogeneous of order $-k-1$,
we have
\[
|a_k(x,\xi)| \le C^{k+1} k! |\xi|^{-k-1 }
\]
with some $C>0$. Using the homogeneity, we get
\begin{equation}   \label{S39}
|\partial_\xi^\alpha a_k(x,\xi)| \le C^{k+|\alpha|+1}\alpha! k! |\xi|^{-k-|\alpha|-1 }
\end{equation}
for $\xi$ is in a complex neighborhood of $\mathbb{R}$ and $x$ in a complex
neighborhood of $\bar{U}'$. Therefore, there exists
a pseudoanalytic symbol $a\sim \sum a_k$, see \cite{T}.
This symbol is defined by
\begin{equation}   \label{S40}
a(x,\xi) = \sum_{k=0}^\infty\varphi_k(\xi) a_k(x,\xi),
\end{equation}
where $\varphi_k$ have the properties
(see \cite[V.2]{T}): $0\le\varphi_k\le1$, $\varphi_k(\xi)=0$
for $|\xi|< 2R\max(k,1)$, $\varphi_k(\xi)=1$
for $|\xi|>3R\max(k,1)$, $|D^\alpha\varphi_k|\le (C/R)^{|\alpha|}$
for $|\alpha|\le 2k$, where $R>1$ is a large parameter.
We will prove next that $a(x,D)$ differs from $A$ by an analytic regularizing
operator.

Let $u\in \mathcal{E}'(U')$. Let $s\ge0$ be such that $u$ can be represented as
a finite sum of derivatives of continuous functions of order not exceeding $s$.

By \eqref{S20a}, for any $N\ge1$,
\begin{equation}  \label{S40a}
\mathcal Au = \sum_{k=0}^{N-1} a_k(x,D)u + \tilde{\mathcal R}_Nu,
\end{equation}
where, despite the strong singularity of $a_k$ at $\xi=0$,
$a_k(x,D)$ have  regular (integrable) Schwartz kernels, and
\begin{align} \nonumber
\tilde{\mathcal R}_Nu(x) &= 2\int_{\mathbb{S}^{n-1}}\int_0^\infty r^{N} R_N(x,r,\omega) u(x+r\omega)\, dr\,d\sigma(\omega)\\
                &= 2\int |x-y|^{N-n+1}  R_N \Big(x, |x-y|, \frac{x-y}{|x-y|}\Big)u(y)\,dy.
\end{align}
We express  $R_N$ in its  Cauchy form  as
\[
R_N(x,r,\omega) = \frac{1}{(N-1)!} \int_0^1 \partial_r^N A_{\text{even},k}(x,tr,\omega) (1-t)^{N-1}\,dt.
\]
We have
\begin{equation}  \label{S41}
\big| D^\alpha \tilde{\mathcal R}_N u \big| \le C^{N}\alpha! \quad\text{in $U'$ for $|\alpha|\le N-s$}.
\end{equation}
Splitting the sum \eqref{S40} into two parts, we write
\begin{equation}   \label{S42}
a(x,D) = \text{Op}\Big(\sum_{k=0}^{N-1}\varphi_k(\xi) a_k(x,\xi)\Big)
+ \text{Op}\Big(\sum_{k=N}^{\infty}\varphi_k(\xi) a_k(x,\xi)\Big).
\end{equation}
For the second term we have (see (3.15) in chapter~V in \cite{T})
\begin{equation}   \label{S43}
\Big| D^\alpha \text{Op}\Big(\sum_{k=N}^{\infty}\varphi_k(\xi) a_k(x,\xi)\Big)u\Big| \le C^N\alpha!  \quad\text{in $U'$ for $|\alpha|\le N-s$}.
\end{equation}
We are left to compare the first sum in \eqref{S42} with the sum in \eqref{S40a}:
\[
\mathcal B_Nu := \text{Op}\Big(\sum_{k=0}^{N-1}(1-\varphi_k(\xi)) a_k(x,\xi) \Big) u
= \mathcal B'_Nu+\mathcal B_N''u,
\]
where
\begin{align*}
\mathcal B'_Nu &=  \text{Op}\Big(\sum_{k=0}^{N-1}(\varphi_N(\xi)-\varphi_k(\xi)) a_k(x,\xi) \Big) u,
\\
\mathcal  B''_Nu &=  \text{Op}\Big(\sum_{k=0}^{N-1}(1-\varphi_N(\xi)) a_k(x,\xi) \Big) u.
\end{align*}
On $\text{supp}\,(\varphi_N-\varphi_k)$, we have $2Rk\le|\xi|\le 3RN$
provided that $k<N$ and, as always, we assume that $R\gg1$.
Using this and \eqref{S39}, we get
\begin{equation}   \label{S44}
|D^\alpha \mathcal B_N'u| \le C(CRN)^{|\alpha|-1+s} \quad\text{in $U'$ for $|\alpha|\le N-1$},
\end{equation}
compare with (3.17) in chapter~V in \cite{T}.

We write $\mathcal B_N''u$ in the form (see \eqref{S21a})
\begin{equation}   \label{S44a}
\mathcal B_N''u(x) = 2  \sum_{k=0}^{N-1} \int A_{\text{even},k}\Big(x,\frac{z}{|z|}\Big) |z|^{k-n+1} \left(1-\varphi_N\right)\check{}*f(z+x)\,dz.
\end{equation}
This implies
\begin{equation}   \label{S45}
|D^\alpha \mathcal B_N''u| \le C^N(CRN)^{|\alpha|+s} \quad\text{in $U'$ for $|\alpha|\le N-1$}.
\end{equation}

Combining \eqref{S41}, \eqref{S43},  \eqref{S44}, and \eqref{S45}, we get
\[
\left| D^\alpha (\mathcal A- a(x,D))u  \right| \le C^N N! \quad\text{in $U'$ for $|\alpha|\le N-s$}.
\]
For $N\ge s$, choose $|\alpha|= N-s$ to conclude that $(\mathcal A- a(x,D))u$ is analytic in $U'$.
\end{proof}


\begin{thebibliography}{ABC}

\bibitem{A1} {I. Alexandrova}, Structure of the semi-classical amplitude for general scattering
relations, {\it Comm. PDE} {\bf 30} (2005), 1505-1535.


\bibitem{A2} {I. Alexandrova}, Structure of the short range amplitude for
general scattering relations, to appear in {\em Asymptotic Analysis}.

\bibitem{AS} D. V. Anosov, Y. G. Sinai,
Certain smooth ergodic systems [Russian],
{\em Uspekhi Mat. Nauk} {\bf 22} (1967), 107--172.

\bibitem{Ar} V. I. Arnold,
Some remarks on flows of line elements and frames,
{\em Sov. Math. Dokl.} {\bf 2} (1961), 562--564.

\bibitem{ArG} V. I. Arnold, A. B. Givental,
{\em Symplectic Geometry}, Dynamical Systems IV, Encyclopaedia
of Mathematical Sciences, Springer Verlag, Berlin, 1990.




\bibitem{CI2} G. Contreras, R. Iturriaga,
{\em Global Minimizers of Autonomous Lagrangians},
22 Colloqio Brasileiro de Matematica, 1999.



\bibitem{CIPP1} G. Contreras, R. Iturriaga, G. Paternain, M. Paternain,
Lagrangian graphs, minimizing measures and {M}a\~n\'e's critical values,
{\em Geom. Funct. Anal.} {\bf 8} (1998), no. 5, 788--809.

\bibitem{CIPP2} G. Contreras, R. Iturriaga, G. Paternain, M. Paternain,
The Palais--Smale condition and Ma\~n\'e's critical values,
{\em Ann. Henri Poincare} {\bf 1} (2000),  no. 4, 655--684.

\bibitem{Cr1} C. B. Croke,
Rigidity and distance between boundary points,
{\em J. Diff. Geom.} {\bf 33} (1991), 445--464.

\bibitem{Cr2} C. B. Croke,
Rigidity theorems in Riemannian geometry,
in {\em Geometric Methods in Inverse Problems and PDE Control},
IMA Vol. Math. Appl. {\bf 137}, Springer, New York, 2004, 47--72.

\bibitem{CDS} C. Croke, N. Dairbekov, V. Sharafutdinov,
Local boundary rigidity of a compact Riemannian manifold with
curvature bounded above,
{\em Trans. Amer. Math. Soc.} {\bf 352} (2000), no. 9, 3937--3956.

\bibitem{DP} N. S. Dairbekov, G. P. Paternain,
{\em Rigidity properties of Anosov optical hypersurfaces}, preprint 2005.

\bibitem{E} G. Eskin,
Inverse scattering problem in anisotropic media,
{\em Comm. Math. Phys.} {\bf 199} (1998), 471--491.

\bibitem{FM} A. Fathi, E. Maderna,
Weak KAM theorem on non compact manifolds,
to appear in {\em Nonlinear Differ. Equ. Appl.}



\bibitem{G} V. Guillemin, Sojourn times and asymptotic properties of
the scattering matrix. Proceedings of the
Oji Seminar on Algebraic Analysis and the RIMS
Symposium on Algebraic Analysis (Kyoto Univ., Kyoto,
1976). {\it Publ. Res. Inst. Math. Sci.} {\bf 12}(1976/77), supplement, 69--88.




\bibitem{Ha} P. Hartman,
{\em Ordinary Differential Equations}, Birkh\"auser, Boston, 1982.

\bibitem{Hor} L. H\"ormander,
{\em The Analysis of Liner Partial Differential
Operators. III}, Springer-Verlag, Berlin-Heildelberg-New York-Tokyo, 1985.


\bibitem{Ko} V. V. Kozlov,
Calculus of variations in the large and classical mechanics,
{\em Russian Math. Surveys} {\bf 40} (1985), no. 2, 37--71.

\bibitem{LaU} M. Lassas, G. Uhlmann,
On determining a Riemannian manifold from the Dirichlet-to-Neumann map,
{\em Annales Scientifiques de L'Ecole Normale Superieure,}
{\bf 34} (2001), 771--787.

\bibitem{LSU} M. Lassas, V. Sharafutdinov, G. Uhlmann,
Semiglobal boundary rigidity for Riemannian metrics,
{\em Math. Ann.} {\bf 325} (2003), 767-793.


\bibitem{LeU}
J. Lee, G. Uhlmann,
Determining anisotropic
real-analytic conductivities by boundary measurements,
{\em Comm. Pure Appl. Math.} {\bf 42} (1989), 1097--1112.

\bibitem{M} R. Michel,
Sur la rigite\'e impos\'ee par la longueur des g\'eod\'esiques,
{\em Invent. Math.} {\bf 65} (1981), 71--83.

\bibitem{Mu} R. G. Mukhometov,
On a problem of reconstructing Riemannian metrics,
{\em Siberian Math. J.} {\bf 22} (1982), no. 3, 420--433.

\bibitem{MR} R. G. Mukhometov, V.G. Romanov,
On the problem of finding an isotropic Riemannian metric
in $n$-dimensional space,
{\em Soviet Math. Dokl.} {\bf 19} (1979), no. 6, 1330--1333.

\bibitem{N1} S. P. Novikov,
Variational methods and periodic solutions of equations
of Kirchhoff type. II,
{\em Functional Anal. Appl.} {\bf 15} (1981), 263--274.

\bibitem{N2} S. P. Novikov,
Hamiltonian formalism and a multivalued analogue of Morse theory,
{\em Russian Math. Surveys} {\bf 37} (1982), no. 5, 1--56.

\bibitem{NS} S. P. Novikov, I. Shmel'tser,
Periodic solutions of the Kirchhoff equations for the free
motion of a rigid body in a liquid, and the extended
Lyusternik-Schnirelmann-Morse theory. I.
{\em J. Functional Anal. Appl.} {\bf 15} (1981), 197--207.


\bibitem{PP} G. P. Paternain, M. Paternain,
Anosov geodesic flows and  twisted symplectic structures,
in {\em International Congress on Dynamical Systems in Montevideo (a tribute to Ricardo Ma\~n\'e)},
F. Ledrappier, J. Lewowicz, S. Newhouse eds,
Pitman Research Notes in Math. {\bf 362} (1996), 132--145.

\bibitem{P} L. Pestov,
{\em Questions of Well-Posedness of the Ray Tomography Problems} [Russian],
Sib. Nauch. Izd., Novosibirsk, 2003.

\bibitem{PSh}
L. N. Pestov, V. A. Sharafutdinov,
Integral geometry of tensor fields on a manifold of negative curvature,
{\em Siberian Math. J.} {\bf 29} (1988), no. 3, 427--441.

\bibitem{PU2} L. Pestov, G. Uhlmann,
On the characterization of the range and inversion formulas for the geodesic
X-ray transform,
{\em International Math. Research Notices}, {\bf 80} (2004), 4331-4347.

\bibitem{PU} L. Pestov, G. Uhlmann,
Two dimensional compact simple Riemannian manifolds are
boundary distance rigid,
{\em Ann. of Math.} {\bf 161}, no. 2 (2005), 1089-1106.


\bibitem{Sh1} V. A. Sharafutdinov,
{\em Integral Geometry of Tensor Fields},
VSP, Utrecht, the Netherlands, 1994.

\bibitem{Sh2} V. A. Sharafutdinov,
Ray transform on Riemannian manifolds.
{\em Eight lectures on integral geometry},
\hbox{http://www.math.washington.edu/\~{}sharafut/Ray\_transform.dvi}.

\bibitem{Sh3} V. A. Sharafutdinov,
An integral geometry problem in a nonconvex domain,
{\em Siberian Math. J.} {\bf 46} (2002), no. 6, 1159--1168.

\bibitem{SSU}  V. Sharafutdinov. M. Skokan, G. Uhlmann,
{\em Regularity of ghosts in tensor tomo\-graphy},  {\em J. Geom. Anal.}   {\bf 15}(2005),  499--542.


\bibitem{SU} P. Stefanov, G. Uhlmann,
Stability estimates for the X-ray transform of tensor fields
and boundary rigidity,
{\em  Duke Math. J.} {\bf 123} (2004), 445-467.

\bibitem{SU1}
P. Stefanov, G.Uhlmann,
Boundary rigidity and stability for generic simple metrics,
{\em J. Amer. Math. Soc.}  {\bf 18} (2005), 975--1003.

\bibitem{SU2} P. Stefanov, G.Uhlmann,
Recent progress on the boundary rigidity problem,
{\em Electron. Res. Announc. Amer. Math. Soc.} {\bf 11} (2005),  64-70.

\bibitem{T}
F. Treves,
{\em Introduction to Pseudodifferential and Fourier Integral Operators, Vol. 1. Pseudodifferential Operators}.
The University Series in Mathematics,  Plenum Press, New York--London, 1980.

\bibitem{Tr}
H. Triebel,
{\em Interpolation Theory, Function Spaces, Differential Operators}.
North-Holland, 1978.


\end{thebibliography}
\end{document}